\documentclass[draft]{amsart}

\usepackage{tikz} \usetikzlibrary{arrows}
\usepackage{amssymb, esint}
\usepackage[draft=false, backref=page]{hyperref}

\numberwithin{equation}{section}
\newtheorem{proposition}[equation]{Proposition}
\newtheorem{theorem}[equation]{Theorem}

\newtheorem{corollary}[equation]{Corollary}
\theoremstyle{definition}
\newtheorem{definition}[equation]{Definition}
\theoremstyle{remark}

\newcommand{\RR}{{\mathbb{R}}}
\newcommand{\R}{{\mathbb{R}}}
\newcommand{\NN}{{\mathbb{N}}}

\newcommand{\CC}{{\mathbb{C}}}

\newcommand{\po}{\partial\Omega}

\newcommand\D{\mathcal{D}}
\newcommand\s{\mathcal{S}}
\newcommand\e{\mathbf{e}}
\newcommand\Tr{\mathop{\mathrm{Tr}}\nolimits}
\newcommand\re{\mathop{\mathrm{Re}}\nolimits}
\newcommand\dist{\mathop{\mathrm{dist}}}
\newcommand\supp{\mathop{\mathrm{supp}}}
\newcommand\diam{\mathop{\mathrm{diam}}}
\newcommand\CAP{\mathop{\mathrm{cap}}\nolimits}
\newcommand\Div{\mathop{\mathrm{div}}\nolimits}
\newcommand\grad{\mathop{\mathrm{grad}}\nolimits}
\newcommand\abs[1]{{\lvert#1\rvert}}

\newcommand\biggabs[1]{{\biggl\lvert#1\biggr\rvert}}
\newcommand\doublebar[1]{{\lVert#1\rVert}}

\newcommand\negphantom[1]{{\ifmmode \hskip -1em
	\setbox0=\hbox{$#1$}%
	\else \setbox0=\hbox{#1}\fi \hskip -\wd0 \relax}}
\newcommand\arr{\dot}
\newcommand\ring{\mathaccent"0017 }
\newcommand\mat{}

\begin{document}
%
%
%
%
%
%
%
%
%
\title[Higher-order elliptic equations]{Higher-order elliptic equations in non-smooth domains: history and recent results}
\author{Ariel Barton}
\address{Ariel Barton, Mathematics Department, University of Missouri, Columbia, Missouri 65211}
\email{bartonae@missouri.edu}
\author{Svitlana Mayboroda}
\address{Svitlana Mayboroda, Department of Mathematics, University of Minnesota, Minneapolis, Minnesota 55455}
\email{svitlana@math.umn.edu}
\thanks{Svitlana Mayboroda is partially supported by  the NSF grants DMS 1220089 (CAREER), DMS 1344235 (INSPIRE), DMR 0212302 (UMN MRSEC Seed grant), and the the Alfred P. Sloan Fellowship.}


\date{}

\begin{abstract} Recent years have brought significant advances in the theory of higher order elliptic equations in non-smooth domains. Sharp pointwise estimates on derivatives of polyharmonic functions in arbitrary domains were established, followed by the higher order Wiener test. Certain boundary value problems for higher order operators with variable non-smooth coefficients were addressed, both in divergence form and in composition form, the latter being adapted to the context of Lipschitz domains. These developments brought new estimates on the fundamental solutions and the Green function, allowing for the lack of smoothness of the boundary or of the coefficients of the equation. Building on our earlier account of history of the subject in \cite{BarM14}, this survey presents the current state of the art, emphasizing the most recent results  and emerging open problems.
\end{abstract}

\maketitle
\tableofcontents

\section{Introduction}

The theory of boundary value problems for second order elliptic operators on Lipschitz domains is a well-developed subject. It has received a great deal of study in the past decades and while some important open questions remain, well-posedness of the Dirichlet, Neumann, and regularity problems in $L^p$ and other function spaces has been extensively studied in the full generality of divergence form operators $-\Div A \nabla$ with bounded measurable coefficients.

The corresponding theory for elliptic equations of order greater than two is much less well developed. Such equations are common in physics and in engineering design, with applications ranging from standard models of elasticity \cite{Mel03} to cutting-edge research of Bose-Einstein condensation in graphene and similar materials \cite{SedGK14p}. They naturally appear in many areas of mathematics too, including conformal geometry (Paneitz operator, $Q$-curvature \cite{Cha05}, \cite{Cha02}), free boundary problems \cite{Ada92}, and non-linear elasticity \cite{Skr94},  \cite{Cia97}, \cite{Ant05}.

It was realized very early in the study of higher order equations that most of the methods developed for the second order scenario break down. Further investigation brought challenging hypotheses and surprising counterexamples, and few general positive results. For instance, Hadamard's 1908 conjecture regarding positivity of the biharmonic Green function \cite{Had08} was actually refuted in 1949 (see \cite{Duffin}, \cite{Garabedian}, \cite{Shapiro}), and later on the weak maximum principle was proved to fail as well, at least in high dimensions \cite{MazR91}, \cite{PipV95A}. Another curious feature is a paradox of passage to the limit for solutions under approximation of a smooth domain by polygons \cite{Bab63}, \cite{MazN86B}.

For the sake of concreteness, we will mention that the prototypical example of a higher-order elliptic operator, well known from the theory of elasticity, is the bilaplacian $\Delta^2=\Delta(\Delta)$; a more general example is the polyharmonic operator $\Delta^m$, $m\geq 2$. The biharmonic problem in a domain $\Omega\subset \RR^n$ with Dirichlet boundary data consists, roughly speaking, of finding a function $u$ such that for given $f$, $g$, $h$,
\begin{equation}
\label{eqn:intro:bih}
\Delta^2 u=h  \text{ in }\Omega, \quad u\big\vert_{\partial\Omega}=f, \quad \partial_\nu u\big\vert_{\partial\Omega}=g,
\end{equation}
subject to appropriate estimates on $u$ in terms of the data.
To make it precise, as usual, one needs to properly interpret restriction of solution to the boundary $u\big\vert_{\partial\Omega}$ and its normal derivative $\partial_\nu u\big\vert_{\partial\Omega}$, as well as specify the desired estimates.

This survey concentrates on three directions in the study of the higher order elliptic problems. First, we discuss the fundamental a priori estimates on solutions to biharmonic and other higher order differential equations in arbitrary bounded domains. For the Laplacian, these properties are described by the maximum principle and by the 1924 Wiener criterion. The case of the polyharmonic operator has been only settled in 2014--2015 \cite{MayMaz14}, \cite{MayMaz15}, and is one of the main subjects of the present review.
Then we turn to the known well-posedness results for higher order boundary problems on Lipschitz domains with data in $L^p$, still largely restricted to the constant coefficient operators and, in particular, to the polyharmonic case. Finally, we present some advancements of the past several years in the theory of variable coefficient higher order equations. In contrast to the second order operators, here the discussion splits according to severals forms of underlying operators.
Let us now outline some details.

On smooth domains the study of higher order differential equations  went hand-in-hand with the second order theory; in particular, the weak maximum principle was established  in 1960 (\cite{Agm60}; see also \cite{Mir48}, \cite{Mir58}). Roughly speaking, for a solution $u$ to the equation $Lu=0$ in~$\Omega$, where $L$ is a differential operator of order~$2m$ and $\Omega\subset \RR^n$ is smooth, the maximum principle guarantees
\begin{equation}\label{eqMax}
\max_{\abs{\alpha}\leq m-1} \doublebar{\partial^\alpha u}_{L^\infty(\Omega)}
\leq
C\max_{\abs{\beta}\leq m-1}\doublebar{\partial^\beta u}_{L^\infty(\partial\Omega)},\end{equation}
with the usual convention that the zeroth-order derivative of $u$ is simply $u$ itself. For the Laplacian ($m=1$), this formula is a slightly weakened formulation of the maximum principle.
In striking contrast with the case of harmonic functions, the maximum principle for an elliptic operator of order $2m\geq 4$ may fail, even in a Lipschitz domain. To be precise, in general, the derivatives of order $(m-1)$ of a solution to an elliptic equation of order $2m$ need not be bounded. However, in the special case of three dimensions, \eqref{eqMax} was proven for the $m$-Laplacian $(-\Delta)^m$ in  domains with Lipschitz boundary, (\cite{PipV93}, \cite{PipV95A}; see also \cite{DahJK84}, \cite{PipV92}, \cite{She06A}, \cite{She06B} for related work), and, by different methods, in three-dimensional domains diffeomorphic to a polyhedron (\cite{KozMR01}, \cite{MazR91}).

Quite recently, in 2014, the boundedness of the $(m-1)$-st derivatives of a solution to the polyharmonic equation $(-\Delta)^m u=0$ was established in arbitrary three-dimensional domains \cite{MayMaz14}. Moreover, the authors derived sharp bounds on the $k$-th derivatives of solutions in higher dimensions, with $k$ strictly less than $m-1$ when the dimension is bigger than 3. These results were accompanied by pointwise estimates on the polyharmonic Green function, also optimal in the class of arbitrary domains. Furthermore, introducing the new notion of polyharmonic capacity, in \cite{MayMaz15} the authors established an analogue of the Wiener test. In parallel with the celebrated 1924 Wiener criterion for the Laplacian, the higher order Wiener test describes necessary and sufficient capacitory conditions on the geometry of the domain corresponding to continuity of the derivatives of the solutions. Some earlier results were also available for boundedness and continuity of the solutions themselves (see, e.g., \cite{Maz79}, \cite{Maz99A}, \cite{Maz02}).

We shall extensively describe all these developments and their historical context in Section~\ref{sec:boundedness}.

Going further, in Sections~\ref{sec:BVP} and~\ref{sec:var} we consider boundary value problems in irregular media. Irregularity can manifest itself through lack of smoothness of the boundary of the domain and/or lack of smoothness of the coefficients of the underlying equation. Section~\ref{sec:BVP} largely concentrates on constant coefficient higher order operators, in particular, the polyharmonic equation, in domains with Lipschitz boundaries.
Large parts of this section are taken verbatim from our earlier survey \cite{BarM14}; we have added some recent results of Brewster, I. Mitrea and M. Mitrea. We have chosen to keep our description of the older results, for completeness, and also to provide background and motivation for study of the boundary problems for operators with variable coefficients---the subject of Section~\ref{sec:var}.

The simplest example is the $L^p$-Dirichlet problem for the bilaplacian
\begin{equation}\label{bih-intro}
\Delta^2 u=0 \quad \mbox{in}\,\,\Omega, \quad u\big\vert_{\partial\Omega}=f\in W_1^p(\po), \quad \partial_\nu u\big\vert_{\partial\Omega}=g \in L^p(\po),
\end{equation}
in which case  the expected sharp estimate on the solution is
\begin{equation}\label{bih-intro2}\doublebar{ N(\nabla u)}_{L^p(\partial\Omega)}
\leq C\doublebar{\nabla_\tau f}_{L^p(\partial\Omega)} + C\doublebar{g}_{L^p(\partial\Omega)},\end{equation}
where  $ N$ denotes the non-tangential maximal function and $W_1^{p}(\po)$ is the Sobolev space of functions with one tangential derivative in $L^p$ (cf.\ Section~\ref{sec:dfn} for precise definitions). In Sections~\ref{sec:L2dirichlet}--\ref{sec:maximumproof} we discuss \eqref{bih-intro} and~\eqref{bih-intro2},  and more general higher-order homogeneous Dirichlet and regularity boundary value problems with constant coefficients, with boundary data in $L^p$.   Section~\ref{sec:convex} describes the specific case of  convex domains. The Neumann problem for the bilaplacian is addressed in Section~\ref{sec:neumann}. In Section~\ref{sec:inhom}, we discuss inhomogeneous boundary value problems (such as problem~\eqref{eqn:intro:bih}, with $h\neq 0$); for such problems it is natural to consider boundary data $f$, $g$ in Besov spaces, which, in a sense, are intermediate between those with Dirichlet and regularity data.

Finally, in Section~\ref{sec:var} we discuss higher order operators with non-smooth coefficients. The results are still very scarce, but the developments of past several years promise to lay a foundation for a general theory.

To begin, let us mention that contrary to the second order scenario, there are several natural generalizations of higher order differential equations to the variable coefficient context.
Recall that the prototypical higher order operator is the biharmonic operator $\Delta^2$; there are two natural ways of writing the biharmonic operator, either as a \emph{composition} $\Delta^2 u= \Delta(\Delta u)$, or in higher-order divergence form
\begin{equation*}\Delta^2 u = \sum_{j=1}^n \sum_{k=1}^n \partial_j \partial_k (\partial_j \partial_k u) = \sum_{\abs{\alpha}=2} \frac{2}{\alpha!} \partial^\alpha(\partial^\alpha u).\end{equation*}
If we regard $\Delta^2$ as a composition of two copies of the Laplacian, then one generalization to variable coefficients is to replace each copy by a more general second-order variable coefficient operator $L_2=-\Div A\nabla$ for some matrix~$A$; this yields operators in  \emph{composition form}
\begin{equation}\label{eqn:intro:composition}
Lu(X) = \Div B(X)\nabla(a(X)\,\Div A(X)\nabla u(X))\end{equation}
for some scalar-valued function~$a$ and two matrices $A$ and~$B$. Conversely, if we regard $\Delta^2$ as a divergence-form operator, we may generalize to a variable-coefficient operator in \emph{divergence form}
\begin{equation}\label{eqn:intro:divergence}
Lu(X)=(-1)^m \Div_m A\nabla^m = (-1)^m \sum_{\abs{\alpha}=\abs{\beta}=m} \partial^\alpha (a_{\alpha\beta}(X)\partial^\beta u(X)).\end{equation}
Both classes of operators will be defined more precisely in Section~\ref{sec:dfn}; we will see that the composition form is closely connected to changes of variables, while operators in divergence form are directly associated to positive bilinear forms.

Sections~\ref{sec:kato} and~\ref{sec:variablediv} discuss higher order operators in divergence form \eqref{eqn:intro:divergence}: these sections discuss, respectively, the Kato problem and the known well-posedness results for higher-order operators, all of which at present require boundary data in fractional smoothness spaces (cf.\ Section~\ref{sec:inhom}).
Section~\ref{sec:variablecomp} addresses well-posedness of the Dirichlet boundary value problem for a fourth order operator in a composition form \eqref{eqn:intro:composition} with data in $L^p$, $p=2$, in particular, generalizing the corresponding results for the biharmonic problem~\eqref{bih-intro}. To date, this is the only result addressing well-posedness with $L^p$ boundary data for higher order elliptic operators with bounded measurable coefficients, and its extensions to more general operators, other values of $p$, and other types of boundary data are still open. Returning to divergence form operators~\eqref{eqn:intro:divergence}, one is bound to start with the very foundations of the theory---the estimates on the fundamental solutions. This is a subject of Section~\ref{sec:fundamental}, and the emerging results are new even in the second order case. Having those at hand, and relying on the Kato problem solution \cite{AusHMT01}, we plan to pass to the study of the corresponding layer potentials and eventually, to the well-posedness of boundary value problems with data in~$L^p$.  In this context, an interesting new challenge, unparalleled in the second order case, is a proper definition of ``natural'' Neumann boundary data. For higher order equations the choice of Neumann data is not unique. Depending on peculiarities of the Neumann operator, one can be led to well-posed and ill-posed problems even for the bilaplacian, and more general operators give rise to new issues related to the coercivity of the underlying form. We extensively discuss these issues in the body of the paper and present a certain functional analytic approach to definitions in Section~\ref{sec:neumann:var}. Finally, Section~\ref{sec:open} lays out open questions and some preliminary results.

To finish the introduction, let us point out a few directions of analysis of higher order operators on non-smooth domains not covered in this survey. First, an excellent expository paper \cite{Maz99B} by Vladimir Maz'ya  on the topic of the Wiener criterion and pointwise estimates details the state of the art in the end of the previous century and  provides a considerably more extended discussion of the questions we raised in Section~\ref{sec:boundedness}, surrounding results and open problems. Here, we have concentrated on recent developments for the polyharmonic equation and their historical context.  Secondly, we do not touch upon the methods and results of the part of elliptic theory studying the behavior of solutions in the domains with isolated singularities,  conical points, cuspidal points, etc. For a good first exposure to that theory, one can consult, e.g., \cite{KozMR01} and references therein. Instead, we have intentionally concentrated on the case of Lipschitz domains, which can display accumulating singularities---a feature drastically affecting both the available techniques and the actual properties of  solutions.

\section{Higher order operators: divergence form and composition form}\label{sec:dfn}

As we pointed out in the introduction, the prototypical higher-order elliptic equation is the biharmonic equation $\Delta^2u=0$, or, more generally, the polyharmonic equation $\Delta^m u=0$ for some integer $m\geq 2$. It naturally arises in numerous applications in physics and in engineering, and in mathematics it is a basic model for a higher-order partial differential equation. For second-order differential equations, the natural generalization of the Laplacian is a divergence-form elliptic operator.  However, it turns out that even \emph{defining} a suitable {general higher-order elliptic operator} with variable coefficients is already a challenging problem with multiple \emph{different} solutions, each of them important in its own right.

Recall that there are two important features possessed by the polyharmonic operator. First, it is a ``divergence form'' operator in the sense that there is an associated positive bilinear form, and this positive bilinear form can be used in a number of ways; in particular, it allows us to define weak solutions in appropriate Sobolev space. Secondly, it is a ``composition operator'', that is, it is defined by composition of several copies of the Laplacian.  Moreover, if one considers the differential equation obtained from the polyharmonic equation by change of variables, the result would again be a composition of second-order operators. Hence, both generalizations are interesting and important for applications, albeit leading to \emph{different} higher-order differential equations.

Let us discuss the details. To start, a general \emph{constant coefficient} elliptic operator is defined as follows.

\begin{definition} \label{dfn:ellipticop}
Let $L$ be an operator acting on functions $u:\RR^n\mapsto \CC^\ell$. Suppose that we may write
\begin{equation} \label{eqn:constsystem}
(Lu)_j =  \sum_{k=1}^\ell \sum_{\abs{\alpha}=\abs{\beta}=m} \partial^\alpha a_{\alpha\beta}^{jk} \partial^\beta u_k
\end{equation}
for some coefficients $a_{\alpha\beta}^{jk}$ defined for all $1\leq j,k\leq \ell$ and all multiindices $\alpha$, $\beta$ of length $n$ with $\abs{\alpha}=\abs{\beta}=m$. Then we say that $L$ is a \emph{differential operator of order~$2m$}.

Suppose the coefficients $a_{\alpha\beta}^{jk}$ are constant and satisfy the Legendre-Hadamard ellipticity condition
\begin{equation}
\label{eqn:constelliptic}
\re \sum_{j,k=1}^\ell \sum_{\abs{\alpha}=\abs{\beta}=m } a_{\alpha\beta}^{jk} \xi^\alpha\xi^\beta \zeta_j \bar\zeta_k
\geq
\lambda \abs{\xi}^{2m} \abs{\zeta}^2
\end{equation}
for all $\xi\in \RR^{n}$ and all $\zeta\in\CC^\ell$, where $\lambda>0$ is a real constant.
Then we say that $L$ is an \emph{elliptic operator of order $2m$}.

If $\ell=1$ we say that $L$ is a \emph{scalar operator} and refer to the equation $Lu=0$ as an elliptic equation; if $\ell>1$ we refer to $Lu=0$ as an elliptic system. If $a^{jk}=a^{kj}$, then we say the operator $L$ is \emph{symmetric}. If $a^{jk}_{\alpha\beta}$ is real for all $\alpha$, $\beta$, $j$, and~$k$, we say that $L$ has real coefficients.
\end{definition}
Here if $\alpha$ is a multiindex of length $n$, then $\partial^\alpha = \partial_{x_1}^{\alpha_1}\dots\partial_{x_n}^{\alpha_n}$.

Now let us discuss the case of variable coefficients.
A divergence-form higher-order elliptic operator is given by
\begin{equation}
\label{eqn:divform}
(Lu)_j(X) =  \sum_{k=1}^\ell \sum_{\abs{\alpha}=\abs{\beta}=m} \partial^\alpha (a_{\alpha\beta}^{jk}(X) \partial^\beta u_k(X)).
\end{equation}
If the coefficients $a_{\alpha\beta}^{jk}:\RR^n\to\CC$ are sufficiently smooth, we may rewrite \eqref{eqn:divform} in nondivergence form
\begin{equation}
\label{eqn:nondivform}
(Lu)_j(X) =  \sum_{k=1}^\ell \sum_{\abs{\alpha}\leq 2m} a_{\alpha}^{jk}(X)
\partial^{\alpha} u_k(X).
\end{equation}
This form is particularly convenient when we allow equations with lower-order terms (note their appearance in \eqref{eqn:nondivform}).

A simple criterion for ellipticity of the operators $L$ of \eqref{eqn:nondivform} is the condition that \eqref{eqn:constelliptic} holds with $a_{\alpha\beta}^{jk}$ replaced by $a_{\alpha}^{jk}(X)$ for any $X\in\RR^n$, that is, that
\begin{equation}
\label{eqn:varelliptic1}
\re \sum_{j,k=1}^\ell \sum_{\abs{\alpha}=2m} a_{\alpha}^{jk}(X) \xi^\alpha \zeta_j \bar\zeta_k
\geq
\lambda \abs{\xi}^{2m}\abs{\zeta}^2
\end{equation}
for any fixed $X\in\RR^n$ and for all $\xi\in\RR^n$, $\zeta\in\CC^2$.
This means in particular that ellipticity is only a property of the highest-order terms of \eqref{eqn:nondivform}; the value of $a^{jk}_{\alpha}$, for $\abs{\alpha}<m$, is not considered.

Returning to divergence form operators~\eqref{eqn:divform}, notice that in this case we have a notion of weak solution; we say that $Lu=h$ weakly if
\begin{equation}\label{eqn:weaksoln}
\sum_{j=1}^\ell \int_\Omega \varphi_j\, h_j
=
\sum_{\abs{\alpha}=\abs{\beta}=m} \sum_{j,k=1}^\ell (-1)^m\int_\Omega \partial^\alpha \varphi_j \, a_{\alpha\beta}^{jk}\, \partial^\beta u_k
\end{equation}
for any function $\varphi:\Omega\mapsto \CC^\ell$ smooth and compactly supported. The right-hand side $\langle \varphi, Lu\rangle$ may be regarded as a bilinear form $A[\varphi,u]$. If $L$ satisfies the ellipticity condition \eqref{eqn:varelliptic1}, then this bilinear form is positive definite. A more general ellipticity condition available in the divergence case is simply that $\langle \varphi, L\varphi\rangle \geq \lambda \doublebar{\nabla^m \varphi}_{L^2}^2$ for all appropriate test functions~$\varphi$ (see formula~\eqref{eqn:varelliptic} below); this condition is precisely that the form $A[\varphi,u]$ be positive definite.

\begin{figure}

\begin{tikzpicture}

\fill[white!85!black] (0,0)--(0.8,1.2)--(1.4,0.2)-- (2.5,0.3) -- (3,0.9) -- (4.2,0.2) -- (4.7,0.5) -- (4.7,-1)--(0,-1)--cycle;
\draw (0,0)--(0.8,1.2)--(1.4,0.2)-- (2.5,0.3) -- (3,0.9) -- (4.2,0.2) -- (4.7,0.5);

\node [below] at (2.35,1.5) {$\Delta u=h$};

\draw[->, >=stealth, thick] (4.4,1.5) to [out=40,in=140] node [above] {$x\mapsto\rho(x)=\tilde x$} (6.6,1.5);

\fill[white!85!black] (6.3,0)--(11,0)--(11,-1)--(6.3,-1)--cycle;
\draw (6.3,0)--(11,0);

\node [below] at (8.65,1.5) {$\displaystyle\abs{J_\rho}\Div \biggl(\frac{1}{\abs{J_\rho}} J_\rho \, J_\rho^t\biggr)\nabla\tilde u=\tilde h$};

\end{tikzpicture}

\caption{The behavior of Laplace's equation after change of variables. Here $\tilde u=u\circ\rho$ and $J_\rho$ is the Jacobean matrix for the change of variables~$\rho$.}
\label{fig:changevars}

\end{figure}
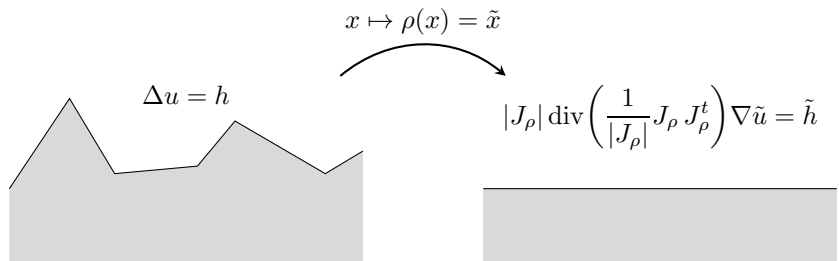

As mentioned above, there is another important form of higher-order operators.
Observe that second-order divergence-form equations arise from a change of variables as follows. If $\Delta u=h$, and $\tilde u=u\circ\rho$ for some  change of variables $\rho$, then
\begin{equation*}a\Div A\nabla \tilde u=\tilde h,\end{equation*}
where $a(X)$ is a real number and $A(X)$ is a real symmetric matrix  (both depending only on~$\rho$; see Figure~\ref{fig:changevars}). In particular, if $u$ is harmonic then $\tilde u$ satisfies the divergence-form equation
\begin{equation*}\sum_{\abs{\alpha}=\abs{\beta}=1}
\partial^\alpha(a_{\alpha\beta}(X)\,\partial^\beta \tilde u(X))=0\end{equation*}
and so the study of divergence-form equations in simple domains (such as the upper half-space) encompasses the study of harmonic functions in more complicated, not necessarily smooth, domains (such as the domain above a Lipschitz graph).

If $\Delta^2 u=0$, however, then after a change of variables $\tilde u$ does not satisfy a divergence-form equation (that is, an equation of the form \eqref{eqn:divform}). Instead, $\tilde u$ satisfies an equation of the following composition form:
\begin{equation}
\label{eqn:prodform}\Div A\nabla (a\Div A\nabla \tilde u)=0.
\end{equation}
In Section~\ref{sec:variablecomp} we shall discuss some new results pertaining to such operators.

Finally, let us mention that throughout we let $C$ and $\varepsilon$ denote positive constants whose value may change from line to line. We let $\fint$ denote the average integral, that is, $\fint_E f\,d\mu=\frac{1}{\mu(E)} \int_E f \,d\mu$. The only measures we will consider are the Lebesgue measure $dX$ (on $\RR^n$ or on domains in $\RR^n$) or the surface measure $d\sigma$ (on the boundaries of domains).

\section{Boundedness and continuity of derivatives of solutions, the maximum principle and the Wiener test}

\label{sec:boundedness}

\subsection{Miranda-Agmon maximum principle and related geometric restrictions on the boundary}
\label{sec:maximum}

The maximum principle for harmonic functions is one of the
fundamental results in the theory of elliptic equations. It holds
in arbitrary domains and guarantees that every solution to the
Dirichlet problem for the Laplace equation, with bounded data, is
bounded. Moreover, it remains valid for all second-order divergence-form elliptic equations with real coefficients.

In the case of equations of higher order, the maximum principle has been established only in relatively nice domains. It was proven to hold for operators with smooth coefficients in smooth domains of dimension two in \cite{Mir48} and \cite{Mir58}, and of arbitrary dimension  in \cite{Agm60}. In the early 1990s, it was extended to
three-dimensional domains diffeomorphic to a polyhedron (\cite{KozMR01,MazR91}) or having a Lipschitz boundary (\cite{PipV93,PipV95A}).  However, in general domains, no direct analog of the maximum principle exists  (see Problem~4.3, p.~275, in
Ne\v cas's book \cite{Nec67}). The increase of the order
leads to the failure of the methods which work for second order
equations, and the properties of the solutions themselves become more
involved.

To be more specific, the following theorem was proved by Agmon. \begin{theorem}[{\cite[Theorem~1]{Agm60}}]
\label{thm:Agmon}
Let $m\geq 1$ be an integer. Suppose that $\Omega$ is domain with $C^{2m}$ boundary. Let
\begin{equation*}
L =\sum_{\abs{\alpha}\leq 2m} a_\alpha(X)\partial^\alpha
\end{equation*}
be a scalar operator of order $2m$, where $a_\alpha\in C^{\abs{\alpha}}(\overline\Omega)$.
Suppose that $L$ is elliptic in the sense of \eqref{eqn:varelliptic1}.
Suppose further that solutions to the Dirichlet problem for $L$ are unique.

Then, for every $u\in C^{m-1}(\overline\Omega)\cap C^{2m}(\Omega)$ that satisfies $Lu=0$ in~$\Omega$, we have
\begin{equation}
\label{eqn:max:lower}
\max_{\abs{\alpha}\leq m-1} \doublebar{\partial^\alpha u}_{L^\infty(\Omega)}
\leq
C\max_{\abs{\beta}\leq m-1}\doublebar{\partial^\beta u}_{L^\infty(\partial\Omega)}.\end{equation}
\end{theorem}
We remark that the requirement that the Dirichlet problem have unique solutions is not automatically satisfied for elliptic equations with lower-order terms; for example, if $\lambda$ is an eigenvalue of the Laplacian then solutions to the Dirichlet problem for $\Delta u-\lambda u$ are not unique.

Equation \eqref{eqn:max:lower} is called the \emph{Agmon-Miranda maximum principle}.
In \cite{Sch75A}, {\v{S}}ul'ce generalized this to systems of the form \eqref{eqn:nondivform}, elliptic in the sense of \eqref{eqn:varelliptic1},
that satisfy a positivity condition (strong enough to imply Agmon's requirement that solutions to the Dirichlet problem be unique).

Thus the Agmon-Miranda maximum principle holds for sufficiently smooth operators and domains. Moreover, for some operators, the maximum principle is valid even in domains with Lipschitz boundary, provided the dimension is small enough. We postpone a more detailed discussion of the Lipschitz case to Section~\ref{sec:maximumproof}; here we simply state the main results.
In \cite{PipV93} and \cite{PipV95A}, Pipher and Verchota showed that the maximum principle holds for the biharmonic operator $\Delta^2$, and more generally for the polyharmonic operator $\Delta^m$, in bounded Lipschitz domains in $\RR^2$ or $\RR^3$. In \cite[Section~8]{Ver96}, Verchota extended this to symmetric, strongly elliptic systems with real constant coefficients in three-dimensional Lipschitz domains.

For Laplace's equation and more general second order elliptic operators, the maximum principle continues to hold in \emph{arbitrary} bounded domains. In contrast, the maximum principle for higher-order operators  in rough domains generally \emph{fails}.

In \cite{MazNP83}, Maz'ya, Nazarov and Plamenevskii studied the Dirichlet problem (with zero boundary data) for constant-coefficient elliptic systems in cones.
Counterexamples to \eqref{eqn:max:lower} for systems of order $2m$ in dimension $n\geq 2m+1$ immediately follow from their results. (See \cite[formulas (1.3), (1.18) and~(1.28)]{MazNP83}.) Furthermore, Pipher and Verchota constructed counterexamples to \eqref{eqn:max:lower} for the biharmonic operator $\Delta^2$ in dimension $n=4$ in \cite[Section~10]{PipV92}, and for the polyharmonic equation  $\Delta^m u=0$ in dimension $n$, $4\leq n<2m+1$, in \cite[Theorem~2.1]{PipV95A}.
Independently Maz'ya and Rossmann showed that \eqref{eqn:max:lower} fails in the exterior of a sufficiently thin cone in dimension~$n$, $n\geq 4$, where $L$ is any constant-coefficient elliptic scalar operator of order $2m\geq 4$ (without lower-order terms). See \cite[Theorem~8 and Remark~3]{MazR92}.

Moreover, with the exception of \cite[Theorem~8]{MazR92}, the aforementioned counterexamples actually provide a stronger negative result than simply the failure of the maximum principle: they show that the left-hand side of \eqref{eqn:max:lower} may be infinite even if the data of the elliptic problem is as nice as possible, that is, smooth and compactly supported.

The counterexamples, however, pertain to high dimensions. This phenomenon raises two fundamental questions: whether the boundedness of the $(m-1)$-st derivatives remains valid in dimensions $n\leq 3$, and whether there are some other, possibly lower-order, estimates that characterize the solutions when $n\geq 4$.  This issue has been completely settled in \cite{MayMaz09B} and \cite{MayMaz14} for the polyharmonic equation in arbitrary domains.

\subsection{Sharp pointwise estimates on the derivatives of solutions in arbitrary domains}
\label{sec:pointwise}

The main results addressing pointwise bounds for solutions to the polyharmonic equation in arbitrary domains are as follows.

\begin{theorem}[\cite{MayMaz14}]\label{tp1.1}
Let $\Omega$ be a bounded domain in $\RR^n$, $2\leq n \leq 2m+1$,
and
\begin{equation}\label{eqp1.2}
(-\Delta)^m u=f \,\,{\mbox{in}}\,\,\Omega, \quad f\in
C_0^{\infty}(\Omega),\quad u\in \ring
W^{m,2}(\Omega).
\end{equation}
Then  the
solution to the boundary value problem \eqref{eqp1.2}
satisfies
\begin{equation}\label{eqp1.3}
\nabla^{m-n/2+1/2} u\in L^\infty(\Omega)\,\,\mbox{when $n$ is odd, } \nabla^{m-n/2} u\in L^\infty(\Omega)\,\,\mbox{when $n$ is even}.
\end{equation}
In particular,
\begin{equation}\label{eqp1.4}
\nabla^{m-1} u\in L^\infty(\Omega)\,\,\mbox{when $n=2,3$}.
\end{equation}

\end{theorem}

Here the space $\ring W^{m,2}(\Omega)$, is, as usual, a
completion of $C_0^\infty(\Omega)$ in the norm given by $\|u\|_{\ring
W^{m,2}(\Omega)}=\|\nabla^m u\|_{L^2(\Omega)}$. We note that $\ring W^{m,2}(\Omega)$ embeds into $C^k(\Omega)$ only when $k$ is strictly smaller than $m-\frac n2$, $n<2m$. Thus, whether the dimension is even or odd, Theorem~\ref{tp1.1} gains one derivative over the outcome of Sobolev embedding.

The results of Theorem~\ref{tp1.1} are sharp, in the sense that the solutions  do not exhibit higher smoothness than warranted by \eqref{eqp1.3}--\eqref{eqp1.4} in general domains. Indeed, assume that $n\in [3,2m+1]\cap\NN$ is odd and let $\Omega\subset{\mathbb{R}}^n$ be the punctured unit ball
$B_1\setminus\{O\}$, where $B_r=\{x\in{\mathbb{R}}^n:\,|x|<r\}$. Consider
a function $\eta\in C_0^\infty(B_{1/2})$ such that $\eta=1$ on
$B_{1/4}$. Then let
\begin{equation}\label{int6}
u(x):=\eta(x)\,\partial_x^{m-\frac n2-\frac 12}(|x|^{2m-n}),\qquad x\in B_1\setminus\{O\},
\end{equation}
where $\partial_x$ stands for a derivative in the direction of $x_i$ for some $i=1,\dots,n$.
It is straightforward to check that $u\in \ring W^{m,2}(\Omega)$ and $(-\Delta)^m u\in
C_0^\infty (\Omega)$. While $\nabla^{m-\frac n2+\frac 12} u$ is bounded, the derivatives of order $m-\frac n2+\frac 32$ are not, and moreover, $\nabla^{m-\frac n2+\frac 12} u$
is not
continuous at the origin. Therefore, the estimates \eqref{eqp1.3}--\eqref{eqp1.4} are optimal in general domains.

As for the case when $n$ is even, the results in \cite[Section~10.4]{KozMR01} demonstrate  that in the exterior of a ray there is an $m$-harmonic function behaving as $|x|^{m-\frac n2+\frac 12}$. Thus, upon truncation by the aforementioned cut-off $\eta$, one obtains a solution to \eqref{eqp1.2} in $B_1\setminus\{x_1=0, \dots, x_{n-1}=0, 0\leq x_n<1\}$, whose derivatives of order $m-\frac n2+1$ are not bounded. More delicate examples can be obtained from our results for the Wiener test to be discussed in Section~\ref{sec:wiener}. Those show that the derivatives of order $m-\frac n2$ need not be continuous in even dimensions. Therefore, in even dimensions \eqref{eqp1.3} is a sharp property as well.

It is worth noting that the results above address also boundedness of solutions (rather than their derivatives) corresponding to the case the case when $m-\frac n2+\frac 12=0$ in odd dimensions, or, respectively, $m-\frac n2=0$ in the even case. In this respect, we would also like to mention higher dimensional results following from the Green function estimates in \cite{Maz99B}. As will be discussed in the next section, one can show that, in addition to our results above,  if $\Omega\subset\RR^n$ is bounded for $n\leq 2m+2$, and if $u$ is a solution to the polyharmonic equation \eqref{eqp1.2}, then $u\in L^\infty(\Omega)$. This result also holds if $\Omega\subset\RR^7$ and $m=2$.

If $\Omega\subset\RR^n$ is bounded and $n\geq 2m+3$, or if $m=2$ and $n\geq8$, then the question of whether solutions $u$ to \eqref{eqp1.2} are bounded is open. In particular, it is not known whether solutions $u$ to
\begin{equation*}
\Delta^2 u=h \text{ in }\Omega, \quad  u\in \mathring W^2_2(\Omega)
\end{equation*}
are bounded if $\Omega\subset\RR^n$ for $n\geq 8$. However, there exists another fourth-order operator whose solutions are \emph{not} bounded in higher-dimensional domains. In \cite{MazN86A}, Maz'ya and Nazarov showed that if $n\geq 8$ and if $a>0$ is large enough, then there exists an open cone $K\subset\RR^n$ and a function $h\in C^\infty_0(\overline K\setminus\{0\})$ such that the solution $u$ to
\begin{equation}
\label{eqn:counterexample}
\Delta^2 u + a \,\partial_{n}^4 u =h \text{ in }K, \quad  u\in \mathring W^2_2(K)
\end{equation}
is unbounded near the origin.

\subsection{Green function estimates}
\label{sec:Green}

Theorem~\ref{tp1.1} has several quantitative manifestations, providing specific estimates on the solutions to \eqref{eqp1.2}. Most importantly, the authors established sharp pointwise estimates on Green's function of the polyharmonic operator and its derivatives, once again without any restrictions on the geometry of the domain.

To start, let us recall the definition of the fundamental solution for the polyharmonic equation (see. e.g., \cite{PolyharmBook}). A fundamental solution for the $m$-Laplacian is a linear combination of the characteristic singular solution (defined below) and any $m$-harmonic function in $\RR^n$. The characteristic singular solution is
\begin{eqnarray}\label{eqp8.2}
&C_{m,n}|x|^{2m-n},\quad &\mbox{if $n$ is odd, or if $n$ is even with $n\geq 2m+2$},\\[4pt]
\label{eqp8.3}
&C_{m,n}|x|^{2m-n} \log |x|,\quad &\mbox{if $n$ is even with $n\leq 2m$}.
\end{eqnarray}
The exact expressions for constants  $C_{m,n}$ can be found in \cite{PolyharmBook}, p.~8.  Hereafter we will use the fundamental solution given by
\begin{equation}\label{eqp8.4}
\Gamma(x)=C_{m,n}  \left\{\begin{array}{l}
|x|^{2m-n},\quad \mbox{if $n$ is odd},\\[4pt]
|x|^{2m-n}\log \frac{{\rm diam}\,\Omega}{|x|},\quad \mbox{if $n$ is even and $n\leq 2m$},\\[4pt]
|x|^{2m-n},\quad \mbox{if $n$ is even and $n\geq 2m+2$}.
\end{array}
\right.
\end{equation}

As is customary, we denote the Green's
function for the polyharmonic equation by $G(x,y)$, $x$, $y\in\Omega$, and its regular part by $S(x,y)$, that is, $S(x,y)=G(x,y)-\Gamma(x-y)$.  By definition, for every
fixed $y\in\Omega$ the function $G(\,\cdot\, ,y)$ satisfies
\begin{equation}\label{eqp8.1}
(-\Delta_x)^m G(x,y)=\delta(x-y), \qquad x\in\Omega,
\end{equation}
in the space $\ring W^{m,2}(\Omega)$. Here  $\Delta_x$ stands for the Laplacian in the $x$ variable. Similarly, we use the notation $\Delta_y$, $\nabla_y$, $\nabla_x$ for
the Laplacian and gradient in $y$, and gradient in $x$,
respectively. By $d(x)$ we denote the distance from $x\in\Omega$ to
$\po$.

\begin{theorem}\label{tp8.1} Let $\Omega\subset \RR^n$ be an arbitrary bounded
domain, $m\in\NN$, $n\in [2,2m+1]\cap\NN$, and let
\begin{equation}\label{eqp7.5}
\lambda=\left\{\begin{array}{l} m-n/2+1/2 \quad \,\,\mbox{when $n$ is odd},\\[8pt] m-n/2\qquad\qquad\mbox{when $n$ is even}.\end{array}
\right.
\end{equation}
 Fix any number $N\geq 25$.
Then there exists a  constant $C$ depending only on $m$, $n$, $N$ such that for
every $x,y\in\Omega$ the following estimates hold.

If $n\in [3,2m+1] \cap\NN$ is odd then
\begin{equation}\label{eqp8.25}
|\nabla_x^i \nabla_y^{j}G(x,y)|\leq
C\,\frac{d(y)^{\lambda-j}}{|x-y|^{\lambda+n-2m+i}}, \quad \mbox{when } |x-y|\geq N\,d(y),\quad 0\leq i,j\leq \lambda,
\end{equation}
and
\begin{equation}\label{eqp8.26}
|\nabla_x^i \nabla_y^{j}G(x,y)|\leq
C\,\frac{d(x)^{\lambda-i}}{|x-y|^{\lambda+n-2m+j}}, \quad \mbox{when } |x-y|\geq N\,d(x),\quad 0\leq i,j\leq \lambda.
\end{equation}
Next,
\begin{equation}\label{eqp8.42}
|\nabla_x^i\nabla_y^jG(x,y)|\leq
\frac{C}{|x-y|^{n-2m+i+j}}, \quad \\[4pt]
\end{equation}
\mbox{when} $|x-y|\leq N^{-1}\max\{d(x),d(y)\}$, {and} $i+j\geq 2m-n$, $0\leq i,j\leq m-n/2+1/2$,
and
\begin{equation}\label{eqp8.43}
|\nabla_x^i\nabla_y^jG(x,y)|\leq
{C}\,{\min\{d(x),d(y)\}^{2m-n-i-j}}, \quad\\[4pt] \end{equation}
\mbox{when} $|x-y|\leq N^{-1}\max\{d(x),d(y)\}$, \mbox{and} $i+j\leq 2m-n$, $0\leq i,j\leq m-n/2+1/2$.
Finally,
\begin{align}\label{eqp8.53}
|\nabla_x^i\nabla_y^jG(x,y)|
&\leq  \frac{C}{\min \{d(x),d(y),|x-y|\}^{n-2m+i+j}}\\[4pt] \nonumber
&\approx \frac{C}{\max \{d(x),d(y),|x-y|\}^{n-2m+i+j}},
\end{align}
\mbox{when} $N^{-1}\,d(x)\leq |x-y|\leq Nd(x)$ \mbox{and} $
N^{-1}\,d(y)\leq |x-y|\leq Nd(y)$, $0\leq i,j\leq \lambda$.

Furthermore, if $n\in [3,2m+1] \cap\NN$ is odd, the estimates on the regular part of the Green function $S$ are as follows:
\begin{equation}\label{eqp8.30}
|\nabla_x^i\nabla_y^j S(x-y)|\leq \frac{C}{|x-y|^{n-2m+i+j}}\,\, \mbox{when}\,\, |x-y|\geq N \min\{d(x),d(y)\},\,\,  0\leq i,j\leq \lambda.\end{equation}
Next,
\begin{equation}\label{eqp8.40}
|\nabla_x^i\nabla_y^jS(x,y)|\leq
\frac{C}{\max\{d(x),d(y)\}^{n-2m+i+j}},
\end{equation}
\mbox{when} $|x-y|\leq N^{-1}\max\{d(x),d(y)\}$, and $i+j\geq 2m-n$, $0\leq i,j\leq m-n/2+1/2$,
and
\begin{equation}\label{eqp8.41}
|\nabla_x^i\nabla_y^jS(x,y)|\leq
{C}\,{\min\{d(x),d(y)\}^{2m-n-i-j}},
\end{equation}
\mbox{when} $|x-y|\leq N^{-1}\max\{d(x),d(y)\}$, \mbox{and} $i+j\leq 2m-n$, $0\leq i,j\leq m-n/2+1/2$.
Finally,
\begin{align}\label{eqp8.54}
|\nabla_x^i\nabla_y^jS(x,y)|
&\leq  \frac{C}{\min \{d(x),d(y),|x-y|\}^{n-2m+i+j}}\\[4pt]\nonumber
&\approx \frac{C}{\max \{d(x),d(y),|x-y|\}^{n-2m+i+j}},
\end{align}
when $N^{-1}\,d(x)\leq |x-y|\leq Nd(x)$ and $
N^{-1}\,d(y)\leq |x-y|\leq Nd(y)$, $0\leq i,j\leq m-n/2+1/2$.

If $n\in [2,2m]\cap\NN$ is even, then \eqref{eqp8.25}--\eqref{eqp8.26} and \eqref{eqp8.53} are valid with $\lambda=m-\frac n2$, and
\begin{equation}\label{eqp8.49}
|\nabla_x^i\nabla_y^jG(x,y)|\leq
{C}\,{\min\{d(x),d(y)\}^{2m-n-i-j}}\,\left( C'+\log \frac{\min\{d(x),d(y)\}}{|x-y|}\right),\end{equation}
when $|x-y|\leq N^{-1}\max\{d(x),d(y)\}$ and $ 0\leq i,j\leq m-n/2$.

Furthermore, if $n\in [2,2m]\cap\NN$ is even, the estimates on the regular part of the Green function $S$ are as follows:
\begin{equation}\label{eqp8.31}
|\nabla_x^i\nabla_y^j S(x-y)|\leq C\,|x-y|^{-n+2m-i-j}\left(C'+\log \frac{{\rm diam}\,{(\Omega)}}{|x-y|}\right)
\end{equation}
when $ |x-y|\geq N \min\{d(x),d(y)\}$, $ 0\leq i,j\leq m-n/2$.
Next,
\begin{equation}\label{eqp8.46}
|\nabla_x^i\nabla_y^jS(x,y)|\leq
{C}\,{\min\{d(x),d(y)\}^{2m-n-i-j}}\left( C'+\log \frac{{\rm diam}\,\Omega}{{\max\{d(x),d(y)\}}}\right),
\end{equation}
when $|x-y|\leq N^{-1}\max\{d(x),d(y)\}$, $0\leq i,j\leq m-n/2$.
Finally,
\begin{multline}\label{eqp8.55}|\nabla_x^i\nabla_y^jS(x,y)|\leq C\,{\min \{d(x),d(y),|x-y|\}^{2m-n-i-j}}\times\\[4pt]\times\left(C'+\log \frac{{\rm diam}\,\Omega}{\max \{d(x),d(y),|x-y|\}^{n-2m+i+j}}\right)
\end{multline}
\mbox{when} $N^{-1}\,d(x)\leq |x-y|\leq Nd(x)$ \mbox{and} $N^{-1}\,d(y)\leq |x-y|\leq Nd(y)$, $ 0\leq i,j\leq m-n/2$.
\end{theorem}

We would like to highlight the most important case of the estimates above, pertaining to the highest order derivatives.

\begin{corollary}\label{cp1.2} Let $\Omega\subset \RR^n$ be an arbitrary bounded
domain.  If $n\in [3,2m+1] \cap\NN$ is odd, then for all $x$, $y\in \Omega$,
\begin{equation}\label{eqp8.8.1}
\left|\nabla_x^{m-\frac n2+\frac 12}\nabla_y^{m-\frac n2+\frac 12} (G(x,y)-\Gamma(x-y))\right|\leq \frac{C}{\max\{d(x),d(y),|x-y|\}},  \end{equation}
and, in particular,
\begin{equation}\label{eqp8.5.1}
\left|\nabla_x^{m-\frac n2+\frac 12}\nabla_y^{m-\frac n2+\frac 12} G(x,y)\right|\leq  \frac{C}{|x-y|}.
\end{equation}

If $n\in [2,2m]\cap\NN$ is even, then for all $x$, $y\in \Omega$,
\begin{multline}\label{eqp8.10.1}
\left|\nabla_x^{m-\frac n2}\nabla_y^{m-\frac n2} (G(x,y)-\Gamma(x-y))\right| \\[4pt]\leq C\, \log \left(1+\frac{{\rm diam}\,\Omega}{\max\{d(x),d(y),|x-y|\}}\right),
\end{multline}
and
\begin{equation}
\label{eqp8.7.1}
\left|\nabla_x^{m-\frac n2}\nabla_y^{m-\frac n2} G(x,y)\right| \leq  C  \log \left(1+\frac{\min\{d(x),d(y)\}}{|x-y|}\right)
.\end{equation}

The constant $C$ in \eqref{eqp8.8.1}--\eqref{eqp8.7.1}  depends on $m$ and $n$ only. In particular, it does not depend on the size or the geometry of the domain $\Omega$.

\end{corollary}

We mention that the pointwise bounds  on the absolute value of Green's function itself have been treated previously in dimensions $2m+1$ and $2m+2$ for  $m>2$ and dimensions $5,6,7$ for $m=2$ in \cite[Section~10]{Maz99B} (see also \cite{Maz79}).
In particular, in \cite[Section~10]{Maz99B}, Maz'ya showed that the Green's function $G_m(x,y)$ for $\Delta^m$ in an arbitrary bounded domain $\Omega\subset\RR^n$ satisfies
\begin{equation}\label{eqn:Greenbound}
\abs{G_m(x,y)}\leq \frac{C(n)}{\abs{x-y}^{n-2m}}
\end{equation}
if $n=2m+1$ or $n=2m+2$. If $m=2$, then \eqref{eqn:Greenbound} also holds in dimension $n=7=2m+3$ (cf.~\cite{Maz79}). Whether \eqref{eqn:Greenbound} holds in dimension $n\geq 8$ (for $m=2$) or $n\geq 2m+3$ (for $m>2$) is an open problem; see \cite[Problem~2]{Maz99B}.
Also, similarly to the case of general solutions discussed above, there exist results for Green functions in smooth domains \cite{DASw}, \cite{Kras}, \cite{Sol1}, \cite{Sol2}, in conical domains \cite{MazPlam}, \cite{KozMR01}, and in polyhedra  \cite{MazR91}.

Furthermore, using standard techniques, the Green's function estimates can be employed to establish the bounds on the solution to \eqref{eqp1.2} for general classes of data $f$, such as $L^p$ for a certain range of $p$, Lorentz spaces etc. A sample statement to this effect is as follows.

\begin{proposition}\label{p9.1} Let $\Omega\subset \RR^n$ be an arbitrary bounded
domain, $m\in\NN$, $n\in [2,2m+1]\cap\NN$, and let $\lambda$ retain the significance of \eqref{eqp7.5}. Consider the
boundary value problem
\begin{equation}\label{eq9.1}
(-\Delta)^m u=\sum_{|\alpha|\leq\lambda}c_\alpha \partial^{\alpha}f_{\alpha}, \quad u\in \ring W^{m,2}(\Omega).
\end{equation}
Then the solution satisfies the following estimates.

If  $n\in [3,2m+1] \cap\NN$ is odd, then for all $x\in\Omega$,
\begin{equation}\label{eq9.2}|\nabla^{m-\frac n2+\frac 12} u(x)|
\leq C_{m,n}  \sum_{|\alpha|\leq  m-\frac n2+\frac 12}
\int_\Omega \frac{d(y)^{m-\frac n2+\frac 12-|\alpha|}}{|x-y|}\,\left|f_\alpha(y)\right|\,dy,
\end{equation}
whenever the integrals on the right-hand side of \eqref{eq9.2} are finite. In particular,
\begin{multline}\label{eq9.3}\|\nabla^{m-\frac n2+\frac 12} u\|_{L^\infty(\Omega)}\\[4pt]\leq C_{m,n,\Omega}    \sum_{|\alpha|\leq  m-\frac n2+\frac 12}
\| d(\,\cdot\,)^{m-\frac n2-\frac 12-|\alpha|}f_{\alpha}\|_{L^p(\Omega)},\quad p>\frac{n}{n-1},
\end{multline}
provided that the norms on the right-hand side of \eqref{eq9.3} are finite.

If  $n\in [2,2m] \cap\NN$ is even, then  for all $x\in\Omega$,
\begin{multline}\label{eq9.12}|\nabla^{m-\frac n2} u(x)|\\[4pt]\leq  C_{m,n}  \sum_{|\alpha|\leq  m-\frac n2}
\int_\Omega d(y)^{m-\frac n2-|\alpha|}\,\log\left(1+\frac{d(y)}{|x-y|}\right)\,\left|f_\alpha(y)\right|\,dy,
\end{multline}
whenever the integrals on the right-hand side of \eqref{eq9.12} are finite. In particular,
\begin{align}\label{eq9.13}\|\nabla^{m-\frac n2} u\|_{L^\infty(\Omega)}&\leq C_{m,n,\Omega}    \sum_{|\alpha|\leq  m-\frac n2}
\| d(\cdot)^{m-\frac n2-|\alpha|}f_{\alpha}\|_{L^p(\Omega)},\,\, p>1,
\end{align}
provided that the norms on the right-hand side of \eqref{eq9.13} are finite.

The constants $C_{m,n}$ above depend on $m$ and $n$ only, while the constants denoted by $C_{m,n,\Omega}$ depend on $m$, $n$, and the diameter of the domain $\Omega$.
\end{proposition}

By the same token, if \eqref{eqn:Greenbound} holds, then solutions to \eqref{eqp1.2} satisfy
\begin{equation*}\doublebar{u}_{L^\infty(\Omega)}\leq C(m,n,p)\diam(\Omega)^{2m-n/p}\doublebar{f}_{L^p(\partial\Omega)}\end{equation*}
provided $p>n/2m$ (see, e.g., \cite[Section~2]{Maz99B}). Thus, e.g., if $\Omega\subset\RR^n$ is bounded for $n=2m+2$,
and if $u$ satisfies \eqref{eqp1.2} for a reasonably nice function~$f$, then $u\in L^\infty(\Omega)$. This result also holds if $\Omega\subset\RR^7$ and $m=2$. This complements the results in Theorem~\ref{tp1.1}, as discussed in Section~~\ref{sec:pointwise}.

To conclude our discussion of Green's functions, we mention two results from \cite{MitM11}; these results are restricted to relatively well-behaved domains. In \cite{MitM11}, D.~Mitrea and I.~Mitrea showed that, if $\Omega$ is a bounded Lipschitz domain in $\RR^3$, and $G$ denotes the Green's function for the bilaplacian~$\Delta^2$, then the estimates
\begin{equation*}\nabla^2 G(x,\,\cdot\,)\in L^3(\Omega),\quad
\dist(\,\cdot\,,\partial\Omega)^{-\alpha} \nabla G(x,\,\cdot\,) \in L^{3/\alpha,\infty}\end{equation*}
hold, uniformly in $x\in\Omega$, for all $0<\alpha\leq 1$.

Moreover, they considered more general elliptic systems. Suppose that $L$ is an arbitrary elliptic operator of order $2m$ with constant coefficients, as defined by Definition~\ref{dfn:ellipticop}, and that $G$ denotes the Green's function for~$L$. Suppose that $\Omega\subset\RR^n$, for $n>m$, is a Lipschitz domain, and that the unit outward normal $\nu$ to $\Omega$ lies in the Sarason space $VMO(\partial\Omega)$ of functions of vanishing mean oscillations on $\partial\Omega$. Then the estimates
\begin{gather}
\label{eqn:MM11}
\nabla^m G(x,\,\cdot\,)\in L^{\frac{n}{n-m},\infty}(\Omega),\\\nonumber
\dist(\,\cdot\,,\partial\Omega)^{-\alpha}\nabla^{m-1} G(x,\,\cdot\,)\in L^{\frac{n}{n-m-1+\alpha},\infty}(\Omega)
\end{gather}
hold, uniformly in $x\in\Omega$, for any $0\leq\alpha\leq 1$.

\subsection{The Wiener test: continuity of solutions}
\label{sec:wiener}

In this section, we discuss conditions that ensure that solutions (or appropriate gradients of solutions) are continuous up to the boundary.
These conditions parallel the famous result of Wiener, who in 1924 formulated a criterion that ensured continuity of \emph{harmonic} functions at boundary points \cite{Wie24}. Wiener's criterion has been
extended to a variety of second-order elliptic and parabolic
equations (\cite{LitSW63,FabJK82B,FabGL89,DalM86, MalZ97, AdaH96, TruW02, Lab02, EvaG82}; see also
the review papers \cite{Maz97, Ada97}). However, as with the maximum principle, extending this criterion to higher-order elliptic equations is a subtle matter, and many open questions remain.

We begin by stating the classical Wiener criterion for the Laplacian. If $\Omega\subset\RR^n$ is a domain and $Q\in\partial\Omega$, then $Q$ is called \emph{regular} for the Laplacian if every solution $u$ to
\begin{equation*}\Delta u=h \text{ in }\Omega,\quad u\in \mathring W^2_1(\Omega)\end{equation*}
for $h\in C^\infty_0(\Omega)$ satisfies $\lim_{X\to Q} u(X)=0$. According to Wiener's theorem \cite{Wie24}, the boundary point $Q\in\po$ is regular if and only if the equation
\begin{equation}
\label{eqn:Wiener2}
\int_0^1 \CAP_{2}(\overline{B(Q,s)}\setminus \Omega) s^{1-n}\,ds=\infty
\end{equation}
holds, where
\begin{equation}\label{Cap2o}
\CAP_{2}(K)=\inf\Bigl\{\doublebar{u}_{L^2(\RR^n)}^2+
\doublebar{\nabla u}_{L^2(\RR^n)}^2:u\in C^\infty_0(\RR^n),\>u\geq 1 \text{ on }K\Bigr\}.\end{equation}

For example, suppose $\Omega$ satisfies the exterior cone condition at~$Q$. That is, suppose there is some open cone $K$ with vertex at $Q$ and some $\varepsilon>0$ such that $K\cap B(Q,\varepsilon)\subset\Omega^C$.
It is elementary to show that $\CAP_2(\overline{B(Q,s)}\setminus\Omega)\geq C(K) s^{n-2}$ for all $0<s<\varepsilon$, and so \eqref{eqn:Wiener2} holds and $Q$ is regular. Regularity of such points was known prior to Wiener (see \cite{Poi90}, \cite{Zar09}, and \cite{Leb13}) and provided inspiration for the formulation of the Wiener test.

By \cite{LitSW63}, if $L=-\Div A\nabla$ is a second-order divergence-form operator, where the matrix $A(X)$ is bounded, measurable, real, symmetric and elliptic, then $Q\in\partial\Omega$ is regular for $L$ if and only if $Q$ and $\Omega$ satisfy \eqref{eqn:Wiener2}. In other words, $Q\in\partial\Omega$ is regular for the Laplacian if and only if it is regular for all such operators.
Similar results hold for some other classes of second-order equations; see, for example,
\cite{FabJK82B},
\cite{DalM86},
or \cite{EvaG82}.

One would like to consider the Wiener criterion for higher-order elliptic equations, and that immediately gives rise to the question of natural generalization of the concept of a regular point. The Wiener criterion for the second order PDEs ensures, in particular, that weak $\mathring W^{2}_1$ solutions are \emph{classical}. That is, the solution approaches its boundary values in the pointwise sense (continuously). From that point of view, one would extend the concept of regularity of a boundary point as continuity of derivatives of order $m-1$ of the solution to an equation of order $2m$ up to the boundary. On the other hand, as we discussed in the previous section, even the boundedness of solutions cannot be guaranteed in general, and thus, in some dimensions the study of the continuity up to the boundary for solutions themselves is also very natural. We begin with the latter question, as it is better understood.

Let us first define a regular point for an arbitrary differential operator $L$  of order $2m$ analogously to the case of the Laplacian, by requiring that every solution $u$ to
\begin{equation}
\label{eqn:regular}
L u=h \text{ in }\Omega,\quad u\in \mathring W^2_m(\Omega)
\end{equation}
for $h\in C^\infty_0(\Omega)$ satisfy $\lim_{X\to Q} u(X)=0$. Note that by the Sobolev embedding theorem, if $\Omega\subset\RR^n$ for $n\leq 2m-1$, then every  $u\in \mathring W^2_m(\Omega)$ is H\"older continuous on $\overline\Omega$ and so satisfies $\lim_{X\to Q} u(X)=0$ at every point $Q\in\partial\Omega$.
Thus, we are only interested in continuity of the solutions at the boundary when $n\geq 2m$.

In this context, the appropriate concept of capacity is the potential-theoretic Riesz capacity of order $2m$, given by
\begin{equation}\label{eqRieszCap}\CAP_{2m}(K)=\inf\Bigl\{\sum_{0\leq\abs{\alpha}\leq m}
\doublebar{\partial^\alpha u}_{L^2(\RR^n)}^2:u\in C^\infty_0(\RR^n),\>u\geq 1 \text{ on }K\Bigr\}.\end{equation}

The following is known. If $m\geq 3$, and if $\Omega\subset\RR^n$ for $n=2m$, $2m+1$ or $2m+2$, or if $m=2$ and $n=4$, $5$, $6$ or~$7$, then
$Q\in\partial\Omega$ is regular for $\Delta^m$ if and only if
\begin{equation}
\label{eqn:Wiener}
\int_0^1 \CAP_{2m}(\overline{B(Q,s)}\setminus \Omega) s^{2m-n-1}\,ds=\infty.
\end{equation}
The biharmonic case was treated in  \cite{Maz77} and \cite{Maz79}, and the polyharmonic case for $m\geq 3$ in  \cite{DonM83} and \cite{Maz99A}.

Let us briefly discuss the method of the proof in order to explain the restrictions on the dimension. Let $L$ be an arbitrary elliptic operator, and let $F$ be the fundamental solution for $L$ in $\RR^n$ with pole at $Q$. We say that $L$ is positive with weight $F$ if, for all $u\in C^\infty_0(\RR^n\setminus\{Q\})$, we have that
\begin{equation}\label{eqn:poswt}
\int_{\RR^n} Lu(X)\cdot u(X)\,F(X)\,dX \geq c\sum_{k=1}^m \int_{\RR^n} \abs{\nabla^k u(X)}^2\abs{X}^{2k-n}\,dX.\end{equation}

The biharmonic operator is positive with weight $F$ in dimension $n$ if $4\leq n\leq 7$, and the polyharmonic operator $\Delta^m$, $m\geq 3$, is positive with weight $F$ in dimension $2m\leq n\leq 2m+2$. (The Laplacian $\Delta$ is positive with weight $F$ in any dimension.) The biharmonic operator $\Delta^2$ is not positive with weight $F$ in dimensions $n\geq 8$, and $\Delta^m$ is not positive with weight $F$ in dimension $n\geq 2m+3$. See \cite[Propositions 1 and~2]{Maz99A}.

The proof of the Wiener criterion for the polyharmonic operator required positivity with weight $F$. In fact, it turns out that positivity with weight $F$ suffices to provide a Wiener criterion for an \emph{arbitrary} scalar elliptic operator with constant coefficients.

\begin{theorem}[{\cite[Theorems 1 and~2]{Maz02}}]
\label{thm:posweightF}
Suppose $\Omega\subset\RR^n$ and that $L$ is a scalar elliptic operator of order $2m$ with constant real coefficients, as defined by Definition~\ref{dfn:ellipticop}.

If $n=2m$, then  $Q\in\partial\Omega$ is regular for $L$ if and only if
\eqref{eqn:Wiener} holds.

If $n\geq 2m+1$, and if the condition \eqref{eqn:poswt} holds, then again $Q\in\partial\Omega$ is regular for $L$ if and only if
\eqref{eqn:Wiener} holds.
\end{theorem}
This theorem is also valid for certain variable-coefficient operators in divergence form; see the remark at the end of \cite[Section~5]{Maz99A}.

Similar results have been proven for some second-order elliptic \emph{systems}. In particular, for the Lam\'e system $Lu=\Delta u+\alpha\grad\Div u$, $\alpha>-1$, positivity with weight $F$ and Wiener criterion have been established for a range of $\alpha$ close to zero, that is, when the underlying operator is close to the Laplacian  (\cite{LuoM10}). It was also shown that positivity with weight $F$ may in general fail for the Lam\'e system. Since the present review is restricted to the higher order operators, we shall not elaborate on this point and instead refer the reader to   \cite{LuoM10} for more detailed discussion.

In the absense of the positivity condition \eqref{eqn:poswt}, the situation is much more involved. Let us point out first that the condition \eqref{eqn:poswt} is \emph{not} necessary for regularity of a boundary point, that is, the continuity of the solutions. There exist fourth-order elliptic operators that are not positive with weight $F$  whose solutions exhibit nice behavior near the boundary; there exist other such operators whose solutions exhibit very bad behavior near the boundary.

Specifically, recall that \eqref{eqn:poswt} fails for $L=\Delta^2$ in dimension $n\geq 8$. Nonetheless, solutions to $\Delta^2u=h$ are often well-behaved near the boundary. By \cite{MazP81}, the vertex of a cone is regular for the bilaplacian in any dimension. Furthermore, if the capacity condition~\eqref{eqn:Wiener} holds with $m=2$, then  by \cite[Section 10]{Maz02}, any solution $u$ to
\begin{equation*}\Delta^2 u=h \text{ in }\Omega,\quad u\in \mathring W^{2}_2(\Omega)\end{equation*}
for $h\in C^\infty_0(\Omega)$ satisfies $\lim_{X\to Q} u(X)=0$ provided the limit is taken along a \emph{nontangential} direction.

Conversely, if $n\geq 8$ and $L=\Delta^2+a\partial_n^4$, then by \cite{MazN86A}, there exists a cone $K$ and a function $h\in  C^\infty_0(\overline K\setminus\{0\})$ such that the solution $u$ to \eqref{eqn:counterexample} is not only discontinuous but \emph{unbounded} near the vertex of the cone.
We remark that a careful examination of the proof in \cite{MazN86A} implies that solutions to \eqref{eqn:counterexample} are unbounded even along some nontangential directions.

Thus, conical points in dimension eight are regular for the bilaplacian and irregular for the operator $\Delta^2 + a\, \partial_{n}^4$. Hence, a relevant Wiener condition \emph{must} use different capacities for these two operators. This is a striking contrast with the second-order case, where the same capacity condition implies regularity for all divergence-form operators, even with variable coefficients.

This concludes the discussion of regularity in terms of continuity of the solution.  We now turn to regularity in terms of continuity of the $(m-1)$-st derivatives. Unfortunately, much less is known in this case.

\subsection{The higher order Wiener test: continuity of derivatives of polyharmonic functions}
\label{sec:wiener2}

The most natural generalization of the Wiener test to the higher order scenario concerns the continuity of the derivatives of the solutions, rather than solutions themselves, as derivatives constitute part of the boundary data. However, a necessary prerequisite for such results is boundedness of the corresponding derivatives of the solutions---an extremely delicate matter in its own right as detailed in Section~\ref{sec:pointwise}. In the context of the polyharmonic equation, Theorem~\ref{tp1.1} has set the stage for an extensive investigation of the Wiener criterion and, following earlier results in \cite{MayMaz09A}, the second author of this paper and V.~Maz'ya have recently obtained a full extension of the Wiener test to the polyharmonic context  in \cite{MayMaz15}. One of the most intricate issues is the proper definition of the polyharmonic capacity, and we start by addressing it.

At this point Theorem~\ref{tp1.1} finally sets the stage for a discussion of the {\it Wiener test} for continuity of the corresponding derivatives of the solution, which brings us to the main results of the present paper.

Assume that $m\in\NN$ and $n\in [2,2m+1]\cap \NN$. Let us denote by $Z$ the following set of indices:
\begin{align}
\label{z1} Z=&\{0,1,\dots,m-n/2+1/2\} \quad \text{if $n$ is odd}, \\[4pt]
\label{z2} Z=&\{-n/2+2, -n/2+4,\dots, m-n/2-2, m-n/2\}\cap(\NN\cup\{0\}) \\\nonumber&\qquad \text{if $n$ is even, $m$ is even}, \\[4pt]
\label{z3} Z=&\{-n/2+1,-n/2+3,\dots, m-n/2-2, m-n/2\}\cap(\NN\cup\{0\})  \\\nonumber&\qquad \text{if $n$ is even, $m$ is odd}.
\end{align}
Now let $\Pi$ be the space of linear combinations of spherical harmonics
\begin{equation}\label{cap6.1}
P(x)=\sum_{p\in Z} \sum_{l=-p}^{p} b_{pl}
Y_l^p(x/|x|), \qquad b_{pl}\in\RR, \quad x\in\RR^n\setminus \{O\},
\end{equation}
with the norm
\begin{equation}\label{cap6.2}
\|P\|_{\Pi}:=\left(\sum_{p\in Z} \sum_{l=-p}^{p} b_{pl}^2\right)^{\frac 12} \quad \mbox{and}\quad \Pi_1:=\{P\in\Pi:\,\|P\|_{\Pi}=1\}.
\end{equation}

Then, given $P\in\Pi_1$, an open set $D$ in $\RR^n$ such that $O\in
\RR^n\setminus D$, and a compactum $K$ in $D$, we define
\begin{multline}\label{cap8}
{\rm Cap}_P\,(K,D)\\[4pt]:=\inf\Bigg\{\int_D |\nabla^m u(x)|^2\,dx:\,\,u\in \ring W^{m,2}(D),\,\,u=P\mbox{ in a
neighborhood of }K\Bigg\},
\end{multline}
with
\begin{equation}\label{cap9}
{\rm Cap}\,(K,D):=\inf_{P\in \Pi_1} {\rm Cap}_P\,(K,D).
\end{equation}
In the context of the Wiener test, we will be working extensively with the capacity of the complement of a domain $\Omega\subset \RR^n$ in the balls $B_{2^{-j}}$, $j\in\NN$, and even more so, in dyadic annuli, $C_{2^{-j}, 2^{-j+2}}$, $j\in\NN$, where $C_{s,as}:=\{x\in\RR^n:\,s<|x|<as\}$, $s$, $a>0$. As is customary, we will drop the reference to the ``ambient'' set
\begin{equation}\label{drop_amb}{{\rm
Cap}_P\,(\overline{C_{2^{-j},2^{-j+2}}}\setminus\Omega)}:= {{\rm
Cap}_P\,(\overline{C_{2^{-j},2^{-j+2}}}\setminus\Omega, C_{2^{-j-2},2^{-j+4}})},\quad j\in\NN,
\end{equation}
and will drop the similar reference for~${\rm Cap}$. In fact, it will be proven below that there are several equivalent definitions of capacity, in particular,
for any $n\in [2, 2m+1]$  and for any $s>0$, $a>0$, $K\subset \overline{C_{s,as}}$, we have
\begin{multline}\label{eq5.12-2-intro} {\rm Cap}_P(K,
C_{s/2,2as})
\\[4pt]
\approx\inf\Bigg\{\sum_{k=0}^m\int_{\RR^n} \frac{|\nabla^ku(x)|^2}{|x|^{2m-2k}}\,dx:\,\,u\in \ring W^{m,2}(\RR^n\setminus\{O\}),\\[-4pt]
u=P\text{ in a neighborhood of }K\Bigg\}.
\end{multline}
In the case when the dimension is odd, also
$${{\rm
Cap}_P\,(\overline{C_{s,as}}\setminus\Omega, C_{s/2, 2as}})\approx {{\rm
Cap}_P\,(\overline{C_{s,as}}\setminus\Omega, \RR^n\setminus \{O\})}.$$
Thus, either of the above can be used in \eqref{drop_amb}, as convenient.

Let $\Omega$ be a domain in $\RR^n$, $n\geq 2$. The point $Q\in\po$ is \emph{$k$-regular} with
respect to the domain $\Omega$ and the operator $(-\Delta)^m$, $m\in\NN$, if the solution to the boundary problem
\begin{equation}\label{eqw6.1}
(-\Delta)^m u=f \,\,{\mbox{in}}\,\,\Omega, \quad f\in
C_0^{\infty}(\Omega),\quad u\in \ring W^{m,2}(\Omega),
\end{equation}
satisfies the condition
\begin{equation}\label{eqw6.2}
\nabla^k u(x)\to 0\mbox{ as } x\to Q,\,x\in\Omega,
\end{equation}
that is, all partial derivatives of $u$ of order $k$ are continuous.
Otherwise, we say that $Q\in\po$ is $k$-irregular.

\begin{theorem}[\cite{MayMaz09A}, \cite{MayMaz15}]\label{tw1.2} Let $\Omega$ be an arbitrary open set in $\RR^n$, $m\in\NN$, $2\leq n \leq 2m+1$. Let $\lambda$ be given by
\begin{equation}\label{eqw7.5}
\lambda=\left\{\begin{array}{l} m-n/2+1/2 \quad \,\,\mbox{when $n$ is odd},\\[8pt] m-n/2\qquad\qquad\mbox{when $n$ is even}.\end{array}
\right.
\end{equation}
If
\begin{equation}\label{eqw1.9}
\sum_{j=0}^\infty 2^{-j(2m-n)} \,\inf_{P\in\Pi_1}{{\rm
Cap}_P\,(\overline{C_{2^{-j},2^{-j+2}}}\setminus\Omega)} =+\infty, \quad \mbox{when $n$ is odd},
\end{equation}
and
\begin{equation}\label{eqw1.9-even}
\sum_{j=0}^\infty j\,2^{-j(2m-n)} \,\inf_{P\in\Pi_1}{{\rm
Cap}_P\,(\overline{C_{2^{-j},2^{-j+2}}}\setminus\Omega)}=+\infty, \quad \mbox{when $n$ is even},
\end{equation}
then the point $O$ is $\lambda$-regular with
respect to the domain $\Omega$ and the operator $(-\Delta)^m$.

Conversely, if the point $O\in\po$ is $\lambda$-regular with
respect to the domain $\Omega$ and the operator $(-\Delta)^m$ then
\begin{equation}\label{eqw1.10}
\inf_{P\in\Pi_1} \sum_{j=0}^\infty 2^{-j(2m-n)} \,{{\rm
Cap}_P\,(\overline{C_{2^{-j},2^{-j+2}}}\setminus\Omega)}=+\infty, \quad \mbox{when $n$ is odd},
\end{equation}
and
\begin{equation}\label{eqw1.10-even}
\inf_{P\in\,\Pi_1} \sum_{j=0}^\infty j\,2^{-j(2m-n)} \,{{\rm
Cap}_P\,(\overline{C_{2^{-j},2^{-j+2}}}\setminus\Omega)}=+\infty, \quad \mbox{when $n$ is even}. \end{equation}
Here, as before,  $C_{2^{-j},2^{-j+2}}$ is the annulus $\{x\in\RR^n:\,2^{-j}<|x|<2^{-j+2}\}$, $j\in\NN\cup \{0\}$.
\end{theorem}

Let us now discuss the results of Theorem~\ref{tw1.2} in more detail. This was the first treatment of the continuity of derivatives of an elliptic equation of order $m\geq 2$ at the boundary, and the first time the capacity \eqref{cap8} appeared in the literature. When applied to the case $m=1$, $n=3$, it yields the classical Wiener criterion for continuity of a harmonic function (cf.\ \eqref{eqn:Wiener2}). Furthermore, as discussed in the previous section, continuity of the solution itself (rather than  its derivatives) has been previously treated for the polyharmonic equation, and for $(-\Delta)^m$ the resulting criterion also follows from Theorem~\ref{tw1.2}, in particular, when $m=2n$, the new notion of capacity \eqref{z1}--\eqref{cap6.2} coincides with the potential-theoretical Bessel capacity used in \cite{Maz02}. In the case $\lambda=0$, covering both of the above, necessary and sufficient condition in Theorem~\ref{tw1.2} are trivially the same, as $P\equiv 1$ when $n=2m$ in even dimensions and $n=2m+1$ in odd ones.
For lower dimensions $n$ the discrepancy is not artificial, for, e.g., \eqref{eqw1.9} may fail to be necessary as was shown in \cite{MayMaz09A}.

It is not difficult to verify that we also recover known bounds in Lipschitz and in smooth domains, as the capacity of a cone and hence, capacity of an intersection with a complement of a Lipschitz domains, assures divergence of the series in \eqref{eqw1.9}--\eqref{eqw1.9-even}. On the other hand, given Theorem~\ref{tw1.2} and following considerations traditional in this context (choosing sufficiently small balls in the consecutive annuli to constitute a complement of the domain), we can build a set with a convergent capacitory integral and, respectively, an irregular solution with discontinuous derivatives of order $\lambda$ at the point $O$. Note that this yields further sharpness of the results of Theorem~\ref{tp1.1}. In particular, in even dimensions, it is a stronger counterexample than that of a continuum  (not only $m-n/2+1$ derivatives are not bounded, but $m-n/2$ derivatives might be discontinuous). We refer the reader back to Section~~\ref{sec:pointwise} for more details.

One of the most difficult aspects of proof of Theorem~\ref{tw1.2} is finding a correct notion of polyharmonic capacity and understanding its key properties. A peculiar choice of linear combinations of spherical harmonics (see \eqref{z1}--\eqref{z3} and \eqref{cap6.1}) is crucial at several stages of the argument, specific to the problem at hand, and no alterations would lead to reasonable necessary and sufficient conditions. At the same time, the new capacity and the notion of higher-order regularity sometimes exhibit surprising properties, such as for instance sensitivity to the affine changes of coordinates \cite{MayMaz09A}, or the aforementioned fact that in sharp contrast with the second order case \cite{LitSW63}, one does not expect the same geometric conditions to be responsible for regularity of solutions to all higher order elliptic equations.

It is interesting to point out that despite fairly involved definitions, capacitory conditions may reduce to a simple and concise criterion, e.g., in a case of a graph. To be precise, let $\Omega\subset\RR^3$ be a domain whose boundary is the graph of a function~$\varphi$, and let $\omega$ be its modulus of continuity.  If
\begin{equation}\label{eqGraph} \int_0^1 \frac{t\,dt}{\omega^2(t)}=\infty,
\end{equation}
then every solution to the biharmonic equation satisfies $\nabla  u\in C(\overline\Omega)$.
Conversely, for every $\omega$ such that the integral in \eqref{eqGraph} is convergent,
there exists a $C^{0,\omega}$ domain and a solution $u$ of the
biharmonic equation such that $\nabla u\notin
C(\overline\Omega)$. In particular, as expected, the gradient of a solution to the biharmonic equation is always bounded in Lipschitz domains and is not necessarily bounded in a H\"older domain. Moreover, one can deduce from \eqref{eqGraph} that the gradient of a solution is always bounded, e.g.,  in a domain with $\omega(t)\approx t\log^{1/2}t$, which is not Lipschitz, and might fail to be bounded in a domain with $\omega(t)\approx t\log t$. More properties of the new capacity and examples can be found in \cite{MayMaz09A}, \cite{MayMaz15}.

\section[Boundary value problems with constant coefficients]{Boundary value problems in Lipschitz domains for elliptic operators with constant coefficients}
\label{sec:BVP}

The maximum principle \eqref{eqn:max:lower} provides estimates on solutions whose boundary data lies in~$L^\infty$. Recall that for second-order partial differential equations with real coefficients, the maximum principle is valid in arbitrary bounded domains. The corresponding sharp estimates for boundary data in~$L^p$, $1<p<\infty$, are much more delicate. They are \emph{not} valid in arbitrary domains, even for harmonic functions, and they depend in a delicate way on the geometry of the boundary. At present, boundary-value problems for the Laplacian and for general real symmetric elliptic operators of the second order are fairly well understood on Lipschitz domains. See, in particular, \cite{Ken94}.

We consider biharmonic functions and more general higher-order elliptic equations. The question of estimates on biharmonic functions with data in $L^p$  was raised by Rivi\`ere in the 1970s
(\cite{CalFS79}), and later Kenig redirected it towards Lipschitz domains in \cite{Ken90, Ken94}. The sharp range of well-posedness in~$L^p$, even for biharmonic functions, remains an open problem (see \cite[Problem 3.2.30]{Ken94}). In this section we shall review the current state of the art in the subject, the main techniques that have been successfully implemented, and their limitations in the higher-order case.

Most of the results we will discuss are valid in Lipschitz domains, defined as follows.

\begin{definition}\label{dfn:Lipschitz} A domain $\Omega\subset\RR^n$ is called a \emph{Lipschitz domain} if, for every $Q\in\partial\Omega$, there is a number $r>0$, a Lipschitz function $\varphi:\RR^{n-1}\mapsto\RR$ with $\doublebar{\nabla \varphi}_{L^\infty}\leq M$, and a rectangular coordinate system for $\RR^n$ such that
\begin{equation*}B(Q,r)\cap \Omega=\{(x,s):x\in\RR^{n-1},\>s\in\RR,\>\abs{(x,s)-Q)}<r,\text{ and } s>\varphi(x)\}.\end{equation*}

If we may take the functions $\varphi$ to be $C^k$ (that is, to possess $k$ continuous derivatives), we say that $\Omega$ is a \emph{$C^k$ domain.}

The {outward normal vector} to $\Omega$ will be denoted~$\nu$.
The {surface measure} will be denoted~$\sigma$, and the tangential derivative along $\partial\Omega$ will be denoted $\nabla_\tau$.
\end{definition}
In this paper, we will assume that all domains under consideration have connected boundary. Furthermore, if $\partial\Omega$ is unbounded, we assume that there is a single Lipschitz function $\varphi$ and coordinate system that satisfies the conditions given above; that is, we assume that $\Omega$ is the domain above (in some coordinate system) the graph of a Lipschitz function.

In order to properly state boundary-value problems on Lipschitz domains, we will need the notions of non-tangential convergence and non-tangential maximal function.

In this and subsequent sections we say that $u\big\vert_{\partial\Omega} =f$ if $f$ is the \emph{nontangential limit} of $u$, that is, if
\begin{equation*}\lim_{X\to Q,\>X\in\Gamma(Q)} u(X)=f(Q)\end{equation*}
for almost every ($d\sigma$) $Q\in\partial\Omega$, where $\Gamma(Q)$ is the \emph{nontangential cone}
\begin{equation}
\label{eqn:nontangcone}
\Gamma(Q)=\{Y\in\Omega:\dist(Y,\partial\Omega)<(1+a)\abs{X-Y}\}.
\end{equation}
Here $a>0$ is a positive parameter; the exact value of $a$ is usually irrelevant to applications.
The \emph{nontangential maximal function} is given by
\begin{equation}
\label{eqn:nontangmax}
NF(Q)=\sup\{\abs{F(X)}:X\in\Gamma(Q)\}.
\end{equation}

The normal derivative of $u$ of order $m$ is defined as
\begin{equation*}\partial_\nu^m u(Q)=\sum_{\abs{\alpha}=m} \nu(Q)^\alpha \frac{m!}{\alpha!}\partial^\alpha u(Q),\end{equation*}
where $\partial^\alpha u(Q)$ is taken in the sense of nontangential limits as usual.

\subsection{The Dirichlet problem: definitions, layer potentials, and some well-posedness results}
\label{sec:L2dirichlet}

We say that the $L^p$-Dirichlet problem for the biharmonic operator $\Delta^2$ in a domain $\Omega$ is well-posed if there exists a constant $C>0$ such that, for every $f\in W^p_1(\partial\Omega)$ and every $g\in L^p(\partial\Omega)$, there exists a unique function $u$ that satisfies
\begin{equation}
\label{eqn:bidirichlet}
\left\{\begin{aligned}
\Delta^2 u &= 0 && \text{in } \Omega,\\
u&=f && \text{on } \partial\Omega,\\
\partial_\nu u&=g && \text{on } \partial\Omega,\\
\doublebar{N(\nabla u)}_{L^p(\partial\Omega)}
&\leq C \doublebar{g}_{L^p(\partial\Omega)}
+C \doublebar{\nabla_\tau f}_{L^p(\partial\Omega)}.
\negphantom{\text{on } \partial\Omega,}
\end{aligned}\right.
\end{equation}

The $L^p$-Dirichlet problem for the polyharmonic operator  $\Delta^m$ is somewhat more involved, because the notion of boundary data is necessarily more subtle.
We say that the $L^p$-Dirichlet problem for $\Delta^m$ in a domain $\Omega$ is well-posed if there exists a constant $C>0$ such that, for every $g\in L^p(\partial\Omega)$ and every $\dot f$ in the Whitney-Sobolev space $W\!A^p_{m-1}(\partial\Omega)$, there exists a unique function $u$ that satisfies
\begin{equation}
\label{eqn:polydirichlet}
\left\{\begin{aligned}
\Delta^m u &= 0 && \text{in } \Omega,\\
\partial^\alpha u\big\vert_{\partial\Omega}&=f_\alpha &&\text{for all } 0\leq\abs{\alpha}\leq m-2,\\
\partial_\nu^{m-1} u&=g \qquad &&\text{on }\partial\Omega,\\
\doublebar{N(\nabla^{m-1} u)}_{L^p(\partial\Omega)}
&\leq C \doublebar{g}_{L^p(\partial\Omega}
+C\sum_{\abs{\alpha}=m-2} \doublebar{\nabla_\tau f_\alpha}_{L^p(\partial\Omega)}
.
\negphantom{\text{for all } 0\leq\abs{\alpha}\leq m-2,}
\end{aligned}\right.
\end{equation}

The space $W\!A^p_m(\partial\Omega)$ is defined as follows.
\begin{definition}\label{dfn:WhitneySobolev} Suppose that $\Omega\subset\RR^n$ is a Lipschitz domain, and  consider arrays of functions $\dot f=\{f_\alpha:\abs{\alpha}\leq m-1\}$ indexed by multiindices $\alpha$ of length~$n$, where $f_\alpha:\partial\Omega\mapsto \CC$.
We let $W\!A^{p}_{m}(\partial\Omega)$ be the completion of the set of arrays $\dot\psi=\{\partial^\alpha\psi: \abs{\alpha}\leq m-1\}$, for $\psi\in C^\infty_0(\RR^n)$, under the norm
\begin{equation}
\label{eqn:whitneynorm}
\sum_{\abs{\alpha}\leq m-1} \doublebar{\partial^\alpha \psi}_{L^p(\partial\Omega)}
+\sum_{\abs{\alpha}=m -1} \doublebar{\nabla_\tau \partial^\alpha\psi}_{L^p(\partial\Omega)}.
\end{equation}
\end{definition}

If we prescribe $\partial^\alpha u=f_\alpha$ on $\partial\Omega$ for some $f\in W\!A^p_m(\partial\Omega)$, then we are prescribing the values of $u$, $\nabla u,\dots,\nabla^{m-1}u$ on $\partial\Omega$, and requiring that (the prescribed part of) $\nabla^{m} u\big\vert_{\partial\Omega}$ lie in $L^p(\partial\Omega)$.

The study of these problems began with biharmonic functions in $C^1$ domains.
In \cite{SelS81}, Selvaggi and Sisto  proved that, if $\Omega$ is the domain above the graph of a compactly supported $C^1$ function $\varphi$, with $\doublebar{\nabla \varphi}_{L^\infty}$ small enough, then solutions to the Dirichlet problem exist provided $1<p<\infty$. Their method used certain biharmonic layer potentials composed with the Riesz transforms.

In \cite{CohG83}, Cohen and Gosselin proved that, if $\Omega$ is a bounded, simply connected $C^1$ domain contained in the plane $\RR^2$, then the $L^p$-Dirichlet problem is well-posed in $\Omega$ for any $1<p<\infty$. In \cite{CohG85}, they extended this result to the complements of such domains.
Their proof used multiple layer potentials introduced by Agmon in \cite{Agm57} in order to solve the Dirichlet problem with continuous boundary data.
The general outline of their proof parallelled that of the proof of the corresponding result \cite{FabJR78} for Laplace's equation.
We remark that by \cite[Theorem~6.30]{MitM13A}, we may weaken the condition that $\Omega$ be $C^1$ to the condition that the unit outward normal $\nu$ to $\Omega$ lies in $VMO(\partial\Omega)$. (Recall that this condition has been used in \cite{MitM11}; see formula~\eqref{eqn:MM11} above and preceding remarks. This condition was also used in \cite{MazMS10}; see Section~\ref{sec:variablediv}.)

As in the case of Laplace's equation, a result in Lipschitz domains soon followed. In \cite{DahKV86}, Dahlberg, Kenig and Verchota showed that the $L^p$-Dirichlet problem for the biharmonic equation is well-posed in any bounded simply connected Lipschitz domain $\Omega\subset\RR^n$, provided $2-\varepsilon<p<2+\varepsilon$ for some $\varepsilon>0$ depending on the domain $\Omega$.

In \cite{Ver87}, Verchota used the construction of \cite{DahKV86} to extend Cohen and Gosselin's results from planar $C^1$ domains to $C^1$ domains of arbitrary dimension. Thus, the $L^p$-Dirichlet problem for the bilaplacian is well-posed for $1<p<\infty$ in $C^1$ domains.

In \cite{Ver90}, Verchota showed that the $L^p$-Dirichlet problem for the polyharmonic operator $\Delta^m$ could be solved for $2-\varepsilon<p<2+\varepsilon$ in starlike Lipschitz domains by induction on the exponent~$m$. He simultaneously proved results for the $L^p$-regularity problem in the same range; we will thus delay discussion of his methods to Section~\ref{sec:regularity}.

All three of the papers \cite{SelS81}, \cite{CohG83} and \cite{DahKV86} constructed biharmonic functions as potentials. However, the potentials used differ. \cite{SelS81} constructed their solutions as
\begin{equation}
\label{eqn:SelS81}
u(X)=\int_{\partial\Omega} \partial_{n}^2 F(X-Y)\,f(Y)\,d\sigma(Y)
+\sum_{i=1}^{n-1}\int_{\partial\Omega} \partial_{i}\partial_{n} F(X-Y) R_i g(Y)\,d\sigma(Y)\end{equation}
where $R_i$ are the Riesz transforms. Here $F(X)$ is the fundamental solution to the biharmonic equation; thus, $u$ is biharmonic in $\RR^n\setminus\partial\Omega$.  As in the case of Laplace's equation,  well-posedness of the Dirichlet problem follows from the boundedness relation  $\doublebar{N(\nabla u)}_{L^p(\partial\Omega)} \leq C\doublebar{f}_{L^p(\partial\Omega)} +C\doublebar{g}_{L^p(\partial\Omega)}$ and from invertibility of the mapping $(f,g)\mapsto (u\big\vert_{\partial\Omega},\partial_\nu u)$ on $L^p(\partial\Omega)\times L^p(\partial\Omega)\mapsto W^p_1(\partial\Omega)\times L^p(\partial\Omega)$.

The multiple layer potential of \cite{CohG83} is an operator of the form
\begin{equation}
\label{eqn:CG83}
\mathcal{L} \dot f(P) = \mathop{\mathrm{p.v.}} \int_{\partial\Omega} \mathcal{L}(P,Q)\dot f(Q)\,d\sigma(Q)\end{equation}
where $\mathcal{L}(P,Q)$ is a $3\times3$ matrix of kernels, also composed of derivatives of the fundamental solution to the biharmonic equation, and $\dot f=(f,f_x,f_y)$ is a ``compatible triple'' of boundary data, that is, an element of $W^{1,p}(\partial\Omega)\times L^p(\partial\Omega)\times L^p(\partial\Omega)$ that satisfies $\partial_\tau f = f_x \tau_x + f_y \tau_y$.
Thus, the input is essentially a function and its gradient, rather than two functions, and the Riesz transforms are not involved.

The method of \cite{DahKV86} is to compose two potentials. First, the function $f\in L^2(\partial\Omega)$ is mapped to its Poisson extension $v$. Next, $u$ is taken to be the solution of the inhomogeneous equation $\Delta u(Y)= (n+2Y\cdot\nabla) v(Y)$ with $u=0$ on $\partial\Omega$. If $G(X,Y)$ is the Green's function for $\Delta$ in $\Omega$ and $k^Y$ is the harmonic measure density at~$Y$, we may write the map $f\mapsto u$ as
\begin{equation}
\label{eqn:DKV}
u(X)=\int_{\Omega} G(X,Y) (n+2Y\cdot \nabla)\int_{\partial\Omega} k^Y(Q)\,f(Q)\,d\sigma(Q)\,dY.
\end{equation}
Since $(n+2Y\cdot\nabla) v(Y)$ is harmonic, $u$ is biharmonic, and so $u$ solves the Dirichlet problem.

\subsection{The \texorpdfstring{$L^p$}{Lp}-Dirichlet problem: the summary of known results on well-posedness and ill-posedness}

Recall that by \cite{Ver90}, the $L^p$-Dirichlet problem is well-posed in Lipschitz domains provided $2-\varepsilon<p<2+\varepsilon$.
As in the case of Laplace's equation (see \cite{FabJL77}), the range $p>2-\varepsilon$ is sharp. That is, for any $p<2$ and any integers $m\geq 2$, $n\geq 2$, there exists a bounded Lipschitz domain $\Omega\subset\RR^n$
such that the $L^p$-Dirichlet problem for $\Delta^m$ is ill-posed in $\Omega$.
See \cite[Section~5]{DahKV86} for the case of the biharmonic operator $\Delta^2$, and the proof of Theorem~2.1 in \cite{PipV95A} for the polyharmonic operator $\Delta^m$.

The range $p<2+\varepsilon$ is not sharp and has been studied extensively. Proving or disproving well-posedness of the $L^p$-Dirichlet problem for $p>2$ in general Lipschitz domains has been an open question since \cite{DahKV86}, and was formally stated as such in \cite[Problem 3.2.30]{Ken94}. (Earlier in \cite[Question~7]{CalFS79}, the authors had posed the more general question of what classes of boundary data give existence and uniqueness of solutions.)

In \cite[Theorem~10.7]{PipV92}, Pipher and Verchota constructed Lipschitz domains $\Omega$ such that the $L^p$-Dirichlet problem for $\Delta^2$ was ill-posed in~$\Omega$, for any given $p>6$ (in four dimensions) or any given $p>4$ (in five or more dimensions).
Their counterexamples built on the study of solutions near a singular point, in particular upon \cite{MazNP83} and \cite{MazP81}. In \cite{PipV95A}, they provided other counterexamples to show that the $L^p$-Dirichlet problem for $\Delta^m$ is ill-posed, provided $p>2(n-1)/(n-3)$ and $4\leq n< 2m+1$. They remarked that if $n\geq 2m+1$, then ill-posedness follows from the results of \cite{MazNP83} provided $p>2m/(m-1)$.

The endpoint result at $p=\infty$ is the Agmon-Miranda maximum principle \eqref{eqn:max:lower} discussed above. We remark that if $2<p_0\leq \infty$, and the $L^{p_0}$-Dirichlet problem is well-posed (or \eqref{eqn:max:lower} holds) then by interpolation, the $L^p$-Dirichlet problem is well-posed for any $2<p<p_0$.

We shall adopt the following definition (justified by the discussion above).
\begin{definition} \label{dfn:prange}
Suppose that $m\geq 2$ and $n\geq 4$. Then $p_{m,n}$ is defined to be the extended real number that satisfies the following properties. If $2\leq p\leq p_{m,n}$, then the $L^p$-Dirichlet problem for $\Delta^m$ is well-posed in any bounded Lipschitz domain $\Omega\subset\RR^n$. Conversely, if $p>p_{m,n}$, then there exists a bounded Lipschitz domain $\Omega\subset\RR^n$ such that the $L^p$-Dirichlet problem for $\Delta^m$ is ill-posed in~$\Omega$. Here, well-posedness for $1<p<\infty$ is meant in the sense of~\eqref{eqn:polydirichlet}, and well-posedness for $p=\infty$ is meant in the sense of the maximum principle (see \eqref{eqn:maximum1} below).
\end{definition}

As in \cite{DahKV86}, we expect the range of solvability for any \emph{particular} Lipschitz domain $\Omega$ to be $2-\varepsilon<p<p_{m,n}+\varepsilon$ for some $\varepsilon$ depending on the Lipschitz character of~$\Omega$.

Let us summarize here the results currently known for $p_{m,n}$. More details will follow in Section~\ref{sec:regularity}.

For any $m\geq 2$, we have that
\begin{itemize}
\item If $n=2$ or $n=3$, then the $L^p$-Dirichlet problem for $\Delta^m$ is well-posed in any Lipschitz domain $\Omega$ for any $2\leq p<\infty$. (\cite{PipV92, PipV95A})
\item If $4\leq n\leq 2m+1$, then $p_{m,n}=2(n-1)/(n-3)$. (\cite{She06A, PipV95A}.)
\item If $n=2m+2$, then $p_{m,n}=2m/(m-1)=2(n-2)/(n-4)$. (\cite{She06B, MazNP83}.)
\item If $n\geq 2m+3$, then $2(n-1)/(n-3)\leq p_{m,n}\leq 2m/(m-1)$. (\cite{She06A, MazNP83}.)
\end{itemize}
The value of $p_{m,n}$, for $n\geq 2m+3$, is open.

In the special case of biharmonic functions ($m=2$), more is known.
\begin{itemize}
\item $p_{2,4}=6$, $p_{2,5}=4$, $p_{2,6}=4$, and  $p_{2,7}=4$. (\cite{She06A} and \cite{She06B})
\item If $n\geq 8$, then
\begin{equation*}2+\frac{4}{n-\lambda_n} <p_{2,n}\leq 4\end{equation*}
where
\begin{equation*}
\lambda_n=\frac{n+10+2\sqrt{2(n^2-n+2)}}{7}.
\end{equation*}
(\cite{She06C})
\item If $\Omega$ is a $C^1$ or convex domain of arbitrary dimension, then the $L^p$-Dirichlet problem for $\Delta^2$ is well-posed in $\Omega$ for any $1<p<\infty$. (\cite{Ver90, She06C, KilS11A}.)
\end{itemize}

We comment on the nature of ill-posedness. The counterexamples of \cite{DahKV86} and \cite{PipV95A} for $p<2$ are failures of uniqueness. That is, those counterexamples are nonzero functions $u$, satisfying $\Delta^m u=0$ in $\Omega$, such that $\partial_\nu^k u=0$ on $\partial\Omega$ for $0\leq k\leq m-1$, and such that $N(\nabla^{m-1} u)\in L^p(\partial\Omega)$.

Observe that if $\Omega$ is bounded and $p>2$, then $L^p(\partial\Omega)\subset L^2(\partial\Omega)$. Because the $L^2$-Dirichlet problem is well-posed, the failure of well-posedness for $p>2$ can only be a failure of the optimal estimate $N(\nabla^{m-1} u)\in L^p(\partial\Omega)$. That is, if the $L^p$-Dirichlet problem for $\Delta^m$ is ill-posed in $\Omega$, then for some Whitney array $\dot f\in W\!A^p_{m-1}(\partial\Omega)$ and some $g\in L^p(\partial\Omega)$, the unique function $u$ that satisfies $\Delta^m u=0$ in $\Omega$, $\partial^\alpha u=f_\alpha$, $\partial_\nu^{m-1} u=g$ and $N(\nabla^{m-1} u)\in L^2(\partial\Omega)$ does not satisfy $N(\nabla^{m-1} u)\in L^p(\partial\Omega)$.

\subsection{The regularity problem and the \texorpdfstring{$L^p$}{Lp}-Dirichlet problem}
\label{sec:regularity}

In this section we elaborate on some of the methods used to prove the Dirichlet well-posedness results listed above, as well as their historical context. This naturally brings up a consideration of a different boundary value problem, the $L^q$-regularity problem for higher order operators.

Recall that for second-order equations the regularity problem corresponds to finding a solution with prescribed tangential gradient along the boundary. In analogy, we say that the $L^q$-regularity problem for $\Delta^m$ is well-posed in $\Omega$ if there exists a constant $C>0$ such that, whenever $\dot f\in W\!A^q_m(\partial\Omega)$, there exists a unique function $u$ that satisfies
\begin{equation}
\label{eqn:polyregularity}
\left\{\begin{aligned}
\Delta^m u &= 0 && \text{in } \Omega,\\
\partial^\alpha u\big\vert_{\partial\Omega}&=f_\alpha &&\text{for all } 0\leq \abs{\alpha}\leq m-1,\\
\doublebar{N(\nabla^{m} u)}_{L^q(\partial\Omega)}
&\leq C \sum_{\abs{\alpha}=m-1} \doublebar{\nabla_\tau f_\alpha}_{L^q(\partial\Omega)}
.
\negphantom{\text{for all } 0\leq \abs{\alpha}\leq m-1,}
\end{aligned}\right.
\end{equation}

There is an important endpoint formulation at $q=1$ for the regularity problem. We say that the $H^1$-regularity problem is well-posed if there exists a constant $C>0$ such that, whenever $\dot f$ lies in the Whitney-Hardy space $H^1_m(\partial\Omega)$, there exists a unique function $u$ that satisfies
\begin{equation*}
\left\{\begin{aligned}
\Delta^m u &= 0 && \text{in } \Omega,\\
\partial^\alpha u\big\vert_{\partial\Omega}&=f_\alpha &&\text{for all } 0\leq \abs{\alpha}\leq m-1,\\
\doublebar{N(\nabla^{m} u)}_{L^1(\partial\Omega)}
&\leq C \sum_{\abs{\alpha}=m-1} \doublebar{\nabla_\tau f_\alpha}_{H^1(\partial\Omega)}.
\negphantom{\text{for all } 0\leq \abs{\alpha}\leq m-1,}
\end{aligned}\right.
\end{equation*}

The space $H^1_{m}(\partial\Omega)$ is defined as follows.
\begin{definition}\label{dfn:WhitneyHardy}
We say that $\dot a\in W\!A^q_m(\partial\Omega)$ is a {$H^1_{m}(\partial\Omega)$-$L^q$ atom} if $\dot a$ is supported  in a ball $B(Q,r)\cap\partial\Omega$ and if
\begin{equation*}\sum_{\abs{\alpha}=m -1} \doublebar{\nabla_\tau a_\alpha}_{L^q(\partial\Omega)}
\leq \sigma(B(Q,r)\cap\partial\Omega)^{1/q-1}
.\end{equation*}
If $\dot f\in W\!A^1_m(\partial\Omega)$ and there are $H^1_m$-$L^2$ atoms $\dot a_k$ and constants $\lambda_k\in\CC$ such that
\begin{equation*}\nabla_\tau f_\alpha=\sum_{k=1}^\infty \lambda_k \nabla_\tau (a_k)_\alpha \text{ for all }\abs{\alpha}={m-1}\end{equation*}
and such that $\sum\abs{\lambda_k}<\infty$, we say that $\dot f\in H^1_m(\partial\Omega)$, with $\doublebar{\dot f}_{H^1_m(\partial\Omega)}$ being the smallest $\sum\abs{\lambda_k}$ among all such representations.
\end{definition}

In \cite{Ver90}, Verchota proved well-posedness of the $L^2$-Dirichlet problem and the $L^2$-regularity problem for the polyharmonic operator $\Delta^m$ in any bounded starlike Lipschitz domain by simultaneous induction.

The base case $m=1$ is valid in all bounded Lipschitz domains by \cite{Dah79} and \cite{JerK81B}. The inductive step is to show that well-posedness for the Dirichlet problem for $\Delta^{m+1}$ follows from well-posedness of the lower-order problems.
In particular, solutions with $\partial^\alpha u=f_\alpha$ may be constructed using the regularity problem for~$\Delta^m$, and the boundary term $\partial_\nu^m u=g$, missing from the regularity data, may be attained using the \emph{inhomogeneous}
Dirichlet problem for $\Delta^m$. On the other hand, it was shown that the well-posedness for the regularity problem for $\Delta^{m+1}$ follows from well-posedness of the lower-order problems and from the Dirichlet problem for~$\Delta^{m+1}$, in some sense, by realizing the solution to the regularity problem as an integral of the solution to the Dirichlet problem.

As regards a broader range of $p$ and $q$,
Pipher and Verchota showed in \cite{PipV92} that the $L^p$-Dirichlet and $L^q$-regularity problems for $\Delta^2$ are well-posed in all bounded Lipschitz domains $\Omega\subset\RR^3$, provided $2\leq p<\infty$ and $1<q\leq 2$. Their method relied on duality. Using potentials similar to those of \cite{DahKV86}, they constructed solutions to the $L^2$-Dirichlet problem in domains above Lipschitz graphs.
The core of their proof was the invertibility on $L^2(\partial\Omega)$ of a certain potential operator~$T$. They were able to show that the invertibility of its adjoint $T^*$ on $L^2(\partial\Omega)$ implies that the $L^2$-regularity problem for $\Delta^2$ is well-posed. Then, using the atomic decomposition of Hardy spaces, they analyzed the $H^1$-regularity problem. Applying interpolation and duality for $T^*$ once again, now in the reverse regularity-to-Dirichlet direction, the full range for both regularity and Dirichlet problems was recovered in domains above graphs. Localization arguments then completed the argument in bounded Lipschitz domains.

In four or more dimensions, further progress relied on the following theorem of Shen.

\begin{theorem}[\cite{She06B}]\label{thm:Shen}
Suppose that $\Omega\subset\RR^n$ is a Lipschitz domain. The following conditions are equivalent.
\begin{itemize}
\item The $L^p$-Dirichlet problem for $L$ is well-posed, where $L$ is a symmetric elliptic system of order $2m$ with real constant coefficients.
\item There exists some constant $C>0$ and some $p>2$ such that\\
\null\hfill\llap{\parbox{\textwidth}
{%
\begin{equation}
\label{eqn:revholder}
\left(\fint_{B(Q,r)\cap\partial\Omega} N(\nabla^{m-1} u)^p\,d\sigma\right)^{1/p}
\leq C\left(\fint_{B(Q,2r)\cap\partial\Omega} N(\nabla^{m-1} u)^2\,d\sigma\right)^{1/2}
\end{equation}}}\\
holds whenever $u$ is a solution to the $L^2$-Dirichlet problem for $L$ in~$\Omega$, with $\nabla u\equiv 0$ on $B(Q,3r)\cap\partial\Omega$.
\end{itemize}
\end{theorem}
For the polyharmonic operator $\Delta^m$, this theorem was essentially proven in \cite{She06A}. Furthermore, the reverse H\"older estimate \eqref{eqn:revholder} with $p=2(n-1)/(n-3)$ was shown to follow from well-posedness of the $L^2$-regularity problem. Thus the $L^p$-Dirichlet problem is well-posed in bounded Lipschitz domains in $\RR^n$ for $p=2(n-1)/(n-3)$. By interpolation, and because reverse H\"older estimates have self-improving properties, well-posedness in the range $2\leq p\leq 2(n-1)/(n-3)+\varepsilon$ for any particular Lipschtiz domain follows automatically.

Using regularity estimates and square-function estimates, Shen was able to further improve this range of $p$. He showed that with $p=2+4/(n-\lambda)$, $0<\lambda<n$,   the reverse H\"older estimate \eqref{eqn:revholder} is true, provided that
\begin{equation}
\label{eqn:bdryholder}
\int_{B(Q,r)\cap\Omega}\abs{\nabla^{m-1} u}^2 \leq
C\left(\frac{r}{R}\right)^\lambda \int_{B(Q,R)\cap\Omega}\abs{\nabla^{m-1} u}^2
\end{equation}
holds whenever $u$ is a solution to the $L^2$-Dirichlet problem in $\Omega$ with $N(\nabla^{m-1}u)\in L^2(\partial\Omega)$ and  $\nabla^{k} u\big\vert_{B(Q,R)\cap\Omega}\equiv 0$  for all $0\leq k\leq m-1$.

It is illuminating to observe that the estimates arising in connection with the pointwise bounds on the solutions in arbitrary domains (cf.~Section~\ref{sec:maximum}) and the Wiener test (cf.~Section~\ref{sec:wiener}), take essentially the form \eqref{eqn:bdryholder}. Thus, Theorem~\ref{thm:Shen} and its relation to \eqref{eqn:bdryholder} provide a  direct way to transform results regarding local boundary regularity of solutions, obtained via the methods underlined in Sections~\ref{sec:maximum} and \ref{sec:wiener},
into well-posedness of the $L^p$-Dirichlet problem.

In particular, consider \cite[Lemma~5]{Maz02}. If $u$ is a solution to $\Delta^m u=0$ in $B(Q,R)\cap\Omega$, where $\Omega$ is a Lipschitz domain, then by \cite[Lemma~5]{Maz02} there is some constant $\lambda_0>0$ such that
\begin{equation}
\label{eqn:wienerholderlipschitz}
\sup_{B(Q,r)\cap\Omega} \abs{u}^2
\leq
\left(\frac{r}{R}\right)^{\lambda_0}
\frac{C}{R^n}\int_{B(Q,R)\cap\Omega} \abs{u(X)}^2\,dX
\end{equation}
provided that $r/R$ is small enough, that $u$ has zero boundary data on $B(Q,R)\cap\partial\Omega$, and where $\Omega\subset\RR^n$ has dimension $n=2m+1$ or $n=2m+2$, or where $m=2$ and $n=7=2m+3$. (The bound on dimension comes from the requirement that $\Delta^m$ be positive with weight~$F$; see equation~\eqref{eqn:poswt}.)

It is not difficult to see (cf., e.g.,  \cite[Theorem~2.6]{She06B}), that \eqref{eqn:wienerholderlipschitz} implies \eqref{eqn:bdryholder} for some $\lambda > n-2m+2$, and thus implies well-posedness of the $L^p$-Dirichlet problem for a certain range of~$p$. This provides an improvement on the results of \cite{She06A} in the case $m=2$ and $n=6$ or $n=7$, and in the case $m\geq 3$ and $n=2m+2$. Shen has stated this improvement in \cite[Theorems 1.4 and~1.5]{She06B}: the  $L^p$-Dirichlet problem for $\Delta^2$ is well-posed for $2\leq p<4+\varepsilon$ in dimensions $n=6$ or $n=7$, and the $L^p$-Dirichlet problem for $\Delta^m$ is well-posed if $2\leq p<2m/(m-1)+\varepsilon$ in dimension $n=2m+2$.

The method of weighted integral identities, related to positivity with weight~$F$ (cf.~\eqref{eqn:poswt}), can be further finessed in a particular case of the biharmonic equation.
\cite{She06C} uses this method (extending the ideas from \cite{Maz79}) to show that if $n\geq 8$, then \eqref{eqn:bdryholder} is valid for solutions to $\Delta^2$ with $\lambda=\lambda_n$, where
\begin{equation}
\label{eqn:lambdan}
\lambda_n=\frac{n+10+2\sqrt{2(n^2-n+2)}}{7}.
\end{equation}

We now return to the $L^q$-regularity problem. Recall that in \cite{PipV92}, Pipher and Verchota showed that if $2<p<\infty$ and $1/p+1/q<1$, then the $L^p$-Dirichlet problem and the $L^q$-regularity problem for $\Delta^2$ are both well-posed in three-dimensional Lipschitz domains. They proved this by showing that, in the special case of a domain above a Lipschitz graph, there is duality between the $L^p$-Dirichlet and $L^q$-regularity problems.
Such duality results are common. See \cite{KenP93}, \cite{She07A}, and \cite{KenR09} for duality results in the second-order case; although even in that case, duality is not always guaranteed. (See \cite{May10A}.)
Many of the known results concerning the regularity problem for the polyharmonic operator $\Delta^m$ are results relating the $L^p$-Dirichlet problem to the $L^q$-regularity problem.

In \cite{MitM10}, I.~Mitrea and M.~Mitrea showed that  if $1<p<\infty$ and $1/p+1/q=1$, and if the $L^q$-regularity problem for $\Delta^2$ and the $L^p$-regularity problem for $\Delta$ were both well-posed in a {particular}  bounded Lipschitz domain~$\Omega$, then the $L^p$-Dirichlet problem for $\Delta^2$ was also well-posed in~$\Omega$. They proved this result (in arbitrary dimensions) using layer potentials and a Green representation formula for biharmonic equations. Observe that the extra requirement of well-posedness for the Laplacian is extremely unfortunate, since in bad domains it essentially restricts consideration to $p<2+\varepsilon$ and thus does not shed new light on well-posedness in the general class of Lipschitz domains. As will be discussed below, later Kilty and Shen established an optimal duality result for biharmonic Dirichlet and regularity problems.

Recall that the formula~\eqref{eqn:revholder} provides a necessary and sufficient condition for well-posedness of the $L^p$-Dirichlet problem.
In \cite{KilS11B}, Kilty and Shen provided a similar condition for the regularity problem. To be precise, they demonstrated that if $q>2$ and $L$ is a symmetric elliptic system of order $2m$ with real constant coefficients, then the $L^q$-regularity problem for $L$ is well-posed if and only if the estimate
\begin{equation}
\label{eqn:revholderreg}
\left(\fint_{B(Q,r)\cap\Omega} N(\nabla^{m} u)^q\,d\sigma\right)^{1/q}\leq C\left(\fint_{B(Q,2r)\cap\Omega} N(\nabla^{m} u)^2\,d\sigma\right)^{1/2}
\end{equation}
holds for all points $Q\in\partial\Omega$, all $r>0$ small enough, and all solutions $u$ to the $L^2$-regularity problem with $\nabla^k u\big\vert_{B(Q,3r)\cap\partial\Omega}=0$ for $0\leq k\leq m-1$.
Observe that \eqref{eqn:revholderreg} is identical to \eqref{eqn:revholder} with $p$ replaced by $q$ and $m-1$ replaced by~$m$.

As a consequence, well-posedness of the $L^q$-regularity problem in $\Omega$ for certain values of $q$ implies well-posedness of the $L^p$-Dirichlet problem for some values of~$p$. Specifically, arguments using interior regularity and fractional integral estimates (given in \cite[Section 5]{KilS11B}) show that \eqref{eqn:revholderreg} implies \eqref{eqn:revholder} with $1/p=1/q-1/(n-1)$.
But recall from \cite{She06B} that \eqref{eqn:revholder} holds if and only if the $L^p$-Dirichlet problem for $L$ is well-posed in~$\Omega$.
Thus, if $2<q<n-1$, and if the $L^q$-regularity problem for a symmetric elliptic system is well-posed in a Lipschitz domain $\Omega$, then the $L^p$-Dirichlet problem for the same system and domain is also well-posed, provided $2<p<p_0+\varepsilon$ where $1/p_0=1/q-1/(n-1)$.

For the bilaplacian, a full duality result is known.
In \cite{KilS11A}, Kilty and Shen showed that, if $1<p<\infty$ and $1/p+1/q=1$, then well-posedness of the $L^p$-Dirichlet problem for $\Delta^2$ in a Lipschitz domain~$\Omega$, and well-posedness of the $L^q$-regularity problem for $\Delta^2$ in~$\Omega$, were both equivalent to the bilinear estimate
\begin{align}
\label{eqn:KSbilinear}
\biggabs{\int_\Omega \Delta u\,\Delta v}
&\leq C
\left(\doublebar{\nabla_\tau\nabla f}_{L^p}+\abs{\partial\Omega}^{-1/(n-1)} \doublebar{\nabla f}_{L^p}+\abs{\partial\Omega}^{-2/(n-1)} \doublebar{f}_{L^p}\right)
\\&\quad\nonumber
\times
\left(\doublebar{\nabla g}_{L^q}+\abs{\partial\Omega}^{-1/(n-1)} \doublebar{g}_{L^q}\right)
\end{align}
for all $f$, $g\in C^\infty_0(\RR^n)$, where $u$ and $v$ are solutions of the $L^2$-regularity problem with boundary data $\partial^\alpha u=\partial^\alpha f$ and $\partial^\alpha v=\partial^\alpha g$. Thus, if $\Omega\subset\RR^n$ is a bounded Lipschitz domain, and if $1/p+1/q=1$, then the $L^p$-Dirichlet problem is well-posed in $\Omega$ if and only if the $L^q$-regularity problem is well-posed in~$\Omega$.

All in all, we see that the $L^p$-regularity problem for $\Delta^2$ is well-posed in $\Omega\subset\RR^n$ if
\begin{itemize}
\item $\Omega$ is $C^1$ or convex, and $1<p<\infty$.
\item $n=2$ or $n=3$ and $1<p<2+\varepsilon$.
\item $n=4$ and $6/5-\varepsilon<p<2+\varepsilon$.
\item $n=5$, $6$ or $7$, and $4/3-\varepsilon < p<2+\varepsilon$.
\item $n\geq 8$, and $2-\frac{4}{4+n-\lambda_n}<p<2+\varepsilon$, where $\lambda_n$ is given by \eqref{eqn:lambdan}. The above ranges of $p$ are sharp, but this range is still open.
\end{itemize}

\subsection{Higher-order elliptic systems}
\label{sec:system}

The polyharmonic operator $\Delta^m$ is part of a larger class of elliptic higher-order operators. Some study has been made of boundary-value problems for such operators and systems.

The $L^p$-Dirichlet problem for a strongly elliptic system $L$ of order $2m$, as defined in Definition~\ref{dfn:ellipticop}, is well-posed in $\Omega$ if there exists a constant $C$ such that, for every $\dot f\in W\!A^p_{m-1}(\partial\Omega\mapsto \CC^\ell)$ and every $\vec g\in L^p(\partial\Omega\mapsto \CC^\ell)$, there exists a unique vector-valued function $\vec{u}:\Omega\mapsto\CC^\ell$ such that
\begin{equation}
\label{eqn:systemdirichlet}
\left\{\begin{aligned}
(L \vec{u})_j = \sum_{k=1}^\ell \sum_{\abs{\alpha}=\abs{\beta}=m } \partial^\alpha a_{\alpha\beta}^{jk} \partial^\beta u_k
&= 0 && \text{in } \Omega \text{ for each }1\leq j\leq \ell,\\
\partial^\alpha \vec{u} &= f_\alpha && \text{on } \partial\Omega \text{ for } \abs{\alpha}\leq m-2,\\
\partial_{\nu}^{m-1} \vec{u} &= \vec g && \text{on } \partial\Omega,\\
\doublebar{N(\nabla^{m-1} u)}_{L^q(\partial\Omega)}
\leq C &\sum_{\abs{\alpha}=m-2} \doublebar{\nabla_\tau f_\alpha}_{L^q(\partial\Omega)}
+C\doublebar{\vec g}_{L^p(\partial\Omega)}
.
\negphantom{\text{in } \Omega \text{ for each }1\leq j\leq \ell,}
\end{aligned}\right.
\end{equation}
The $L^q$-regularity problem is well-posed in $\Omega$ if there is some constant $C$ such that, for every $\dot f\in W\!A^p_{m}(\partial\Omega\mapsto \CC^\ell)$, there exists a unique $\vec u$ such that
\begin{equation}
\label{eqn:systemregularity}
\left\{\begin{aligned}
(L \vec{u})_j = \sum_{k=1}^\ell \sum_{\abs{\alpha}=\abs{\beta}=m } \partial^\alpha a_{\alpha\beta}^{jk} \partial^\beta u_k &= 0 && \text{in } \Omega \text{ for each }1\leq j\leq \ell,\\
\partial^\alpha \vec{u} &= f_\alpha && \text{on } \partial\Omega \text{ for } \abs{\alpha}\leq m-1,\\
\doublebar{N(\nabla^{m} u)}_{L^q(\partial\Omega)}
&\leq C \sum_{\abs{\alpha}=m-1} \doublebar{\nabla_\tau f_\alpha}_{L^q(\partial\Omega)}
.
\negphantom{\text{in } \Omega \text{ for each }1\leq j\leq \ell,}
\end{aligned}\right.
\end{equation}

In \cite{PipV95B}, Pipher and Verchota showed that the $L^p$-Dirichlet and $L^p$-reg\-u\-lar\-ity problems were well-posed for $2-\varepsilon<p<2+\varepsilon$, for any higher-order elliptic partial differential equation with real constant coefficients, in Lipschitz domains of arbitrary dimension. This was extended to symmetric elliptic systems in \cite{Ver96}. A key ingredient of the proof was the boundary G\aa rding inequality
\begin{multline*}
\frac{\lambda}{4}\int_{\partial\Omega} \abs{\nabla^m u} (-\nu_{n})\,d\sigma
\\\leq
\sum_{j,k=1}^\ell \sum_{\abs{\alpha}=\abs{\beta}=m}
\int_{\partial\Omega} \partial^\alpha a_{\alpha\beta}^{jk} \partial^\beta u_k (-\nu_n)\,d\sigma
+C\int_{\partial\Omega}\abs{\nabla^{m-1}\partial_n u}^2\,d\sigma
\end{multline*}
valid if  $u\in C^\infty_0(\RR^n)^\ell$, if $L=\partial^\alpha a_{\alpha\beta}^{jk} \partial^\beta$ is a symmetric elliptic system with real constant coefficients, and if $\Omega$ is the domain above the graph of a Lipschitz function. We observe that in this case, $(-\nu_n)$ is a positive number bounded from below. Pipher and Verchota then used this G\r{a}rding inequality and a Green's formula to construct the nontangential maximal estimate. See \cite{PipV95A} and \cite[Sections 4 and~6]{Ver96}.

As in the case of the polyharmonic operator $\Delta^m$, this first result concerned the $L^p$-Dirichlet problem and $L^q$-regularity problem only for $2-\varepsilon<p<2+\varepsilon$ and for $2-\varepsilon<q<2+\varepsilon$. The polyharmonic operator $\Delta^m$ is an elliptic system, and so we cannot in general improve upon the requirement that $2-\varepsilon<p$ for well-posedness of the $L^p$-Dirichlet problem.

However, we can improve on the requirement $p<2+\varepsilon$.
Recall that Theorem~\ref{thm:Shen} from \cite{She06B}, and its equivalence to \eqref{eqn:bdryholder}, were proven in the general case of strongly elliptic systems with real symmetric constant coefficients. As in the case of the polyharmonic operator $\Delta^m$, \eqref{eqn:revholder} follows from well-posedness of the $L^2$-regularity problem provided $p=2(n-1)/(n-3)$, and so if $L$ is such a system, the $L^p$-Dirichlet problem for $L$ is well-posed in $\Omega$ provided $2-\varepsilon<p<2(n-1)/(n-3)+\varepsilon$. This is \cite[Corollary~1.3]{She06B}.
Again, by the counterexamples of \cite{PipV95A}, this range cannot be improved if $m\geq 2$ and $4\leq n\leq 2m+1$; the question of whether this range can be improved for general operators~$L$ if $n\geq 2m+2$ is still open.

Little is known concerning the regularity problem in a broader range of~$p$. Recall that \eqref{eqn:revholderreg} from \cite{KilS11B} was proven in the general case of strongly elliptic systems with real symmetric constant coefficients. Thus, we known that for such systems, well-posedness of the $L^q$-regularity problem for $2<q<n-1$ implies well-posedness of the $L^p$-Dirichlet problem for appropriate~$p$. The question of whether the reverse implication holds, or whether this result can be extended to a broader range of~$q$, is open.

\subsection{The area integral}

\label{sec:square}

One of major tools in the theory of second-order elliptic differential equations is the Lusin area integral, defined as follows.
If $w$ lies in $W^2_{1,loc}(\Omega)$ for some domain $\Omega\subset\RR^n$, then the area integral (or square function) of $w$ is defined for $Q\in\partial\Omega$ as
\begin{equation*}
S w(Q) = \left(\int_{\Gamma(Q)} \abs{\nabla w(X)}^2 \dist(X,\partial\Omega)^{2-n}dX\right)^{1/2}.
\end{equation*}
In \cite{Dah80A}, Dahlberg showed that if $u$ is harmonic in a bounded Lipschitz domain $\Omega$, if $P_0\in\Omega$ and $u(P_0)=0$, then for any $0<p<\infty$,
\begin{equation}
\label{eqn:square2}
\frac{1}{C}\int_{\partial\Omega} Su^p\,d\sigma
\leq \int_{\partial\Omega} (Nu)^p\,d\sigma
\leq C\int_{\partial\Omega} (Su)^p\,d\sigma
\end{equation}
for some constants $C$ depending only on $p$, $\Omega$ and~$P_0$. Thus, the Lusin area integral bears deep connections to the $L^p$-Dirichlet problem. In \cite{DahJK84}, Dahlberg, Jerison and Kenig generalized this result to solutions to second-order divergence-form elliptic equations with real coefficients for which the $L^r$-Dirichlet problem is well-posed for at least one~$r$.

If $L$ is an operator of order $2m$, then the appropriate estimate is
\begin{equation}
\label{eqn:squarehigher}
\frac{1}{C}
\int_{\partial\Omega} N(\nabla u^{m-1})^p\,d\sigma
\leq
\int_{\partial\Omega} S(\nabla u^{m-1})^p\,d\sigma
\leq
C \int_{\partial\Omega} N(\nabla u^{m-1})^p\,d\sigma.
\end{equation}
Before discussing their validity for particular operators,  let us point out that such square-function estimates are very useful in the study of higher-order equations. In \cite{She06B}, Shen used \eqref{eqn:squarehigher} to prove the equivalence of \eqref{eqn:bdryholder} and \eqref{eqn:revholder}, above. In \cite{KilS11A}, Kilty and Shen used \eqref{eqn:squarehigher} to prove that well-posedness of the $L^p$-Dirichlet problem for $\Delta^2$ implies the bilinear estimate \eqref{eqn:KSbilinear}. The proof of the maximum principle \eqref{eqn:max:lower} in \cite[Section~8]{Ver96} (to be discussed in Section~\ref{sec:maximumproof}) also exploited \eqref{eqn:squarehigher}. Estimates on square functions can be used to derive estimates on Besov space norms; see \cite[Proposition~S]{AdoP98}.

In \cite{PipV91}, Pipher and Verchota proved that \eqref{eqn:squarehigher} (with $m=2$) holds for solutions $u$ to $\Delta^2 u=0$, provided $\Omega$ is a bounded Lipschitz domain, $0<p<\infty$, and $\nabla u(P_0)=0$ for some fixed $P_0\in\Omega$.
Their proof was an adaptation of Dahlberg's proof \cite{Dah80A} of the corresponding result for harmonic functions. They used the $L^2$-theory for the biharmonic operator \cite{DahKV86}, the representation formula  \eqref{eqn:DKV}, and the $L^2$-theory for harmonic functions to prove  good-$\lambda$ inequalities, which, in turn, imply $L^p$ estimates for $0<p<\infty$.

In  \cite{DahKPV97}, Dahlberg, Kenig, Pipher and Verchota proved that \eqref{eqn:squarehigher} held for solutions $u$ to $Lu=0$, for a symmetric elliptic system $L$ of order $2m$ with real constant coefficients, provided as usual that $\Omega$ is a bounded Lipschitz domain, $0<p<\infty$, and $\nabla^{m-1} u(P_0)=0$ for some fixed $P_0\in\Omega$. The argument is necessarily considerably more involved than the argument of \cite{PipV91} or \cite{Dah80A}. In particular, the bound $\doublebar{S(\nabla^{m-1}u)}_{L^2(\partial\Omega)}\leq C\doublebar{N(\nabla^{m-1}u)}_{L^2(\partial\Omega)}$
was proven in three steps.

The first step was to reduce from the elliptic system $L$ of order~$2m$ to the scalar elliptic operator $M=\det L$ of order $2\ell m$, where $\ell$ is as in formula~\eqref{eqn:constsystem}.
The second step was to reduce to elliptic equations of the form $\sum_{\abs{\alpha}=m} a_\alpha \partial^{2\alpha} u=0$, where $\abs{a_\alpha}>0$ for all $\abs{\alpha}=m$. Finally, it was shown that for operators of this form
\begin{equation*}
\sum_{\abs{\alpha}=m} \int_\Omega a_\alpha \,\partial^\alpha u(X)^2\,\dist(X,\partial\Omega)\,dX
\leq C
\int_{\partial\Omega} N(\nabla^{m-1} u)^2\,d\sigma.
\end{equation*}
The passage to $0<p<\infty$ in \eqref{eqn:squarehigher} was done, as usual,  using good-$\lambda$ inequalities.
We remark that these arguments used the result of \cite{PipV95B} that the $L^2$-Dirichlet problem is well-posed for such operators $L$ in Lipschitz domains.

It is quite interesting that for second-order elliptic systems, the only currently known approach to the square-function estimate \eqref{eqn:square2} is this reduction to a higher-order operator.

\subsection{The maximum principle in Lipschitz domains}
\label{sec:maximumproof}

We are now in a position to discuss the maximum principle \eqref{eqn:max:lower} for higher-order equations in Lipschitz domains.

We say that the maximum principle for an operator $L$ of order $2m$ holds in the bounded Lipschitz domain $\Omega$ if there exists a constant $C>0$ such that, whenever $f\in W\!A^\infty_{m-1}(\partial\Omega)\subset W\!A^2_{m-1}(\partial\Omega)$ and $g\in L^\infty(\partial\Omega)\subset L^2(\partial\Omega)$, the solution $u$ to the Dirichlet problem \eqref{eqn:systemdirichlet} with boundary data $f$ and $g$ satisfies
\begin{equation}
\label{eqn:maximum1}
\doublebar{\nabla^{m-1} u}_{L^\infty}
\leq C \doublebar{g}_{L^\infty(\partial\Omega)}
+C\sum_{\abs{\alpha}=m-2} \doublebar{\nabla_\tau f_\alpha}_{L^\infty(\partial\Omega)}.
\end{equation}

The maximum principle \eqref{eqn:maximum1} was proven to hold in three-dimensional Lipschitz domains by Pipher and Verchota in \cite{PipV93} (for biharmonic functions), in \cite{PipV95A} (for polyharmonic functions), and by Verchota in \cite[Section~8]{Ver96} (for solutions to symmetric systems with real constant coefficients). Pipher and Verchota also proved in \cite{PipV93} that the maximum principle was valid for biharmonic functions in $C^1$ domains of arbitrary dimension. In \cite[Theorem~1.5]{KilS11A}, Kilty and Shen observed that the same techinque gives validity of the maximum principle for biharmonic functions in convex domains of arbitrary dimension.

The proof of \cite{PipV93} uses the $L^2$-regularity problem in the domain $\Omega$ to construct the Green's function $G(X,Y)$ for $\Delta^2$ in~$\Omega$.
Then if $u$ is biharmonic in $\Omega$ with $N(\nabla u)\in L^2(\partial\Omega)$, we have that
\begin{equation}
\label{eqn:PipV93}
u(X)=\int_{\partial\Omega} u(Q)\,\partial_\nu \Delta G(X,Q)\,d\sigma(Q) +\int_{\partial\Omega} \partial_\nu u(Q)\,\Delta G(X,Q)\,d\sigma(Q)
\end{equation}
where all derivatives of $G$ are taken in the second variable~$Q$.
If the $H^1$-regularity problem is well-posed in appropriate subdomains of~$\Omega$, then $\nabla^2 \nabla_X G(X,\,\cdot\,)$ is in $L^1(\partial\Omega)$ with $L^1$-norm independent of~$X$, and so the second integral is at most $C\doublebar{\partial_\nu u}_{L^\infty(\partial\Omega)}$. By taking Riesz transforms, the normal derivative $\partial_\nu \Delta G(X,Q)$ may be transformed to tangential derivatives $\nabla_\tau\Delta G(X,Q)$; integrating by parts transfers these derivatives to~$u$. The square-function estimate \eqref{eqn:squarehigher} implies that the Riesz transforms of $\nabla_X \Delta_Q G(X,Q)$ are bounded on $L^1(\partial\Omega)$. This completes the proof of the maximum principle.

Similar arguments show that the maximum principle is valid for more general operators. See \cite{PipV95A} for the polyharmonic operator, or  \cite[Section~8]{Ver96} for arbitrary symmmetric operators with real constant coefficients.

An important transitional step is the  well-posedness of the $H^1$-regularity problem. It was established in three-dimensional (or $C^1$) domains in \cite[Theorem~4.2]{PipV93} and  \cite[Theorem~1.2]{PipV95A} and discussed in \cite[Section~7]{Ver96}. In each case, well-posedness was proven by analyzing solutions with atomic data $\dot f$ using a technique from \cite{DahK90}. A crucial ingredient in this technique is the well-posedness of the $L^p$-Dirichlet problem for some $p<(n-1)/(n-2)$; the latter is valid if $n=3$ by \cite{DahKV86}, and (for $\Delta^2$) in $C^1$ and convex domains by \cite{Ver90} and \cite{KilS11A}, but fails in general Lipschitz domains for $n\geq 4$.

\subsection{Biharmonic functions in convex domains}
\label{sec:convex}

We say that a domain $\Omega$ is \emph{convex} if, whenever $X$, $Y\in\Omega$, the line segment connecting $X$ and $Y$ lies in $\Omega$. Observe that all convex domains are necessarily Lipschitz domains but the converse does not hold. Moreover, while convex domains are in general no smoother than Lipschitz domains, the extra geometrical structure often allows for considerably stronger results.

Recall that in \cite{MayMaz09A}, the second author of this paper and Maz'ya showed that the gradient of a biharmonic function is bounded in a three-dimensional domain. This is a sharp property in dimension three, and in higher dimensional domains the solutions can be even less regular (cf.\ Section~\ref{sec:maximum}). However, using some intricate linear combination of weighted integrals, the same authors showed in
 \cite{MayMaz08} that \emph{second} derivatives to biharmonic functions were locally bounded when the domain was convex. To be precise, they showed that if $\Omega$ is convex, and  $u\in \mathring W^2_2(\Omega)$ is a solution to $\Delta^2 u=h$ for some $h\in C^\infty_0(\Omega\setminus B(Q,10R))$, $R>0$, $Q\in\po$, then
\begin{equation}
\label{eqn:hessian}
\sup_{B(Q,R/5)\cap\Omega} \abs{\nabla^2 u}
\leq
\frac{C}{R^2}\left(\fint_{\Omega\cap B(Q,5R)\setminus B(Q,R/2)} \abs{u}^2\right)^{1/2}.
\end{equation}
In particular, not only are all boundary points of convex domains $1$-regular, but the gradient $\nabla u$ is Lipschitz continuous near such points.

Kilty and Shen noted in \cite{KilS11A} that \eqref{eqn:hessian} implies that \eqref{eqn:revholderreg} holds in convex domains for any $q$; thus, the $L^q$-regularity problem for the bilaplacian is well-posed for any $2<q<\infty$ in a convex domain. Well-posedness of the $L^p$-Dirichlet problem for $2<p<\infty$ has been established by Shen in \cite{She06C}. By the duality result \eqref{eqn:KSbilinear}, again from \cite{KilS11A}, this implies that both the $L^p$-Dirichlet and $L^q$-regularity problems are well-posed, for any $1<p<\infty$ and any $1<q<\infty$, in a convex domain of arbitrary dimension. They also observed that, by the techniques of \cite{PipV93} (discussed in Section~\ref{sec:maximumproof} above),  the maximum principle \eqref{eqn:maximum1} is valid in arbitrary convex domains.

It is interesting to note how, once again, the methods and
results related to pointwise estimates, the Wiener criterion,
and local regularity estimates near the boundary are intertwined
with the well-posedness of boundary problems in~$L^p$.

\subsection{The Neumann problem for the biharmonic equation}
\label{sec:neumann}

So far we have only discussed the Dirichlet and regularity problems for higher order operators. Another common and important boundary-value problem that arises in applications is the Neumann problem.
Indeed, the principal physical motivation for the inhomogeneous biharmonic equation $\Delta^2 u=h$ is that it describes the equilibrium position of a thin elastic plate subject to a vertical force~$h$. The Dirichlet problem $u\big\vert_{\partial\Omega}=f$, $\nabla u\big\vert_{\partial\Omega}=g$ describes an elastic plate whose edges are clamped, that is, held at a fixed position in a fixed orientation. The Neumann problem, on the other hand, corresponds to the case of a free boundary.
Guido Sweers has written an excellent short paper \cite{Swe09} discussing the boundary conditions that correspond to these and other physical situations.

More precisely, if a thin two-dimensional plate is subject to a force $h$ and the edges are free to move, then its displacement $u$ satisfies the boundary value problem
\begin{equation*}
\left\{
\begin{aligned}
\Delta^2 u &= h && \text{in } \Omega,\\
\rho \Delta u+(1-\rho)\partial_\nu^2 u&=0 &&\text{on }\partial\Omega,\\
\partial_\nu \Delta u+(1-\rho) \partial_{\tau\tau\nu}u&=0 &&\text{on }\partial\Omega.\\
\end{aligned}\right.
\end{equation*}
Here $\rho$ is a physical constant, called the Poisson ratio. This formulation goes back to Kirchoff and is well known in the theory of elasticity; see, for example, Section~3.1 and Chapter~8 of the classic engineering text \cite{Nad63}.
We remark that by \cite[formula \hbox{(8-10)}]{Nad63},
\begin{equation*}
\partial_\nu \Delta u + (1-\rho) \partial_{\tau\tau\nu}u=
\partial_\nu \Delta u + (1-\rho) \partial_{\tau}
\left(\partial_{\nu\tau} u\right).\end{equation*}

This suggests the following homogeneous boundary value problem in a Lipschitz domain $\Omega$ of arbitrary dimension. We say that the $L^p$-Neumann problem is well-posed if there exists a constant $C>0$ such that, for every $f_0\in L^p(\partial\Omega)$ and $\Lambda_0\in W^p_{-1}(\partial\Omega)$, there exists a function $u$ such that
\begin{equation}
\label{eqn:neumann}
\left\{
\begin{aligned}
\Delta^2 u ={}& 0 && \text{in } \Omega,\\
M_\rho u:={}&\rho \Delta u + (1-\rho) \partial_\nu^2 u = f_0&& \text{on }\partial\Omega,\\
K_{\rho} u:={}&\partial_\nu \Delta u + (1-\rho) \frac{1}{2}
\partial_{\tau_{ij}}
\left(\partial_{\nu\tau_{ij}} u \right) = \Lambda_0&& \text{on }\partial\Omega,\\
\doublebar{N(\nabla^2 u)}_{L^p(\partial\Omega)}
\leq{}& C
\doublebar{f_0}_{W^p_1(\partial\Omega)}
+C\doublebar{\Lambda_0}_{W^p_{-1}(\partial\Omega)}.
\negphantom{\text{on }\partial\Omega,}
\end{aligned}\right.
\end{equation}
Here $\tau_{ij} = \nu_i \e_j-\nu_j \e_i$ is a vector orthogonal to the outward normal $\nu$ and lying in the $x_ix_j$-plane.

In addition to the connection to the theory of elasticity, this problem is of interest because it is in some sense adjoint to the Dirichlet problem \eqref{eqn:bidirichlet}. That is, if $\Delta^2 u=\Delta^2 w=0$ in~$\Omega$, then $\int_{\partial\Omega}\partial_\nu w\,M_\rho u- w\, K_\rho u\,d\sigma=\int_{\partial\Omega}\partial_\nu u\,M_\rho w- u\, K_\rho w\,d\sigma$, where $M_\rho$ and $K_\rho$ are as in~\eqref{eqn:neumann}; this follows from the more general formula
\begin{align}
\label{eqn:biharmonicDN}
\int_\Omega w\,\Delta^2 u =
\int_\Omega \left(\rho \Delta u\,\Delta w
+(1-\rho)\partial_{jk} u\,\partial_{jk}w\right)
+\int_{\partial\Omega} w\, K_\rho u-\partial_\nu w\,M_\rho u\,d\sigma
\end{align}
valid for arbitrary smooth functions. This formula is analogous to the classical Green's identity for the Laplacian
\begin{equation}\label{eqn:GrFla}
\int_{\Omega} w\,\Delta u
=
-\int_\Omega \nabla u\cdot \nabla w +\int_{\partial\Omega} w\,\nu\cdot\nabla u\,d\sigma.
\end{equation}

Observe that, contrary to the Laplacian or more general second order operators, there is a \emph{family} of relevant Neumann data for the biharmonic equation. Moreover, different values (or, rather,   ranges) of $\rho$ correspond to different natural physical situations. We refer the reader to  \cite{Ver05} for a detailed discussion.

In \cite{CohG85}, Cohen and Gosselin showed that the $L^p$-Neumann problem \eqref{eqn:neumann} was well-posed in $C^1$ domains contained in $\RR^2$ for for $1<p<\infty$, provided in addition that $\rho=-1$.
The method of proof was as follows.
Recall from \eqref{eqn:CG83} that Cohen and Gosselin showed that the $L^p$-Dirichlet problem was well-posed by constructing a multiple layer potential $\mathcal{L}\dot f$ with boundary values $(I+\mathcal{K})\dot f$, and showing that $I+\mathcal{K}$ is invertible.
We remark that because Cohen and Gosselin preferred to work with Dirichlet boundary data of the form $(u,\partial_x u,\partial_y u)\big\vert_{\partial\Omega}$ rather than of the form $(u,\partial_\nu u)\big\vert_{\partial\Omega}$, the notation of \cite{CohG85} is somewhat different from that of the present paper.
In the notation of the present paper, the method of proof of \cite{CohG85} was to observe that $(I+\mathcal{K})^* \dot\theta$ is equivalent to $(K_{-1}v\dot\theta,M_{-1}v\dot \theta)_{\partial\Omega^C}$, where $v$ is another biharmonic layer potential and $(I+\mathcal{K})^*$ is the adjoint to $(I+\mathcal{K})$. Well-posedness of the Neumann problem then follows from invertibility of $I+\mathcal{K}$ on~$\partial\Omega^C$.

In \cite{Ver05}, Verchota investigated the Neumann problem \eqref{eqn:neumann} in full generality. He considered Lipschitz domains with compact, connected boundary contained in $\RR^n$, $n\geq 2$.
He showed that if $-1/(n-1)\leq \rho<1$, then the Neumann problem is well-posed provided $2-\varepsilon<p<2+\varepsilon$. That is, the solutions exist, satisfy the desired estimates, and are unique either modulo functions of an appropriate class, or (in the case where $\Omega$ is unbounded) when subject to an appropriate growth condition. See \cite[Theorems 13.2 and~15.4]{Ver05}. Verchota's proof also used boundedness and invertibility of certain potentials on $L^p(\partial\Omega)$; a crucial step was a coercivity estimate $\doublebar{\nabla^2 u}_{L^2(\partial\Omega)}\leq C\doublebar{K_\rho u}_{W^2_{-1}(\partial\Omega)}+C\doublebar{M_\rho u}_{L^2(\partial\Omega)}$. (This estimate is valid provided $u$ is biharmonic and satisfies some mean-value hypotheses; see \cite[Theorem~7.6]{Ver05}).

More recently, in \cite{She07B}, Shen improved upon Verchota's results by extending the range on $p$ (in bounded simply connected Lipschitz domains) to $2(n-1)/(n+1)-\varepsilon<p<2+\varepsilon$ if $n\geq 4$, and $1<p<2+\varepsilon$ if $n=2$ or $n=3$. This result again was proven by inverting layer potentials. Observe that the $L^p$-regularity problem is also known to be well-posed for $p$ in this range, and (if $n\geq 6$) in a broader range of~$p$; see Section~\ref{sec:regularity}. The question of the sharp range of $p$ for which the $L^p$-Neumann problem is well-posed in a Lipschitz domain is still open.

Finally, in \cite[Section~6.5]{MitM13A}, I. Mitrea and M. Mitrea showed that if $\Omega\subset\R^n$ is a simply connected domain whose unit outward normal $\nu$ lies in~$VMO(\partial\Omega)$ (for example, if $\Omega$ is a $C^1$ domain), then the acceptable range of~$p$ is $1<p<\infty$; this may be seen as a generalization of the result of Cohen and Gosselin to higher dimensions, to other values of~$\rho$, and to slightly rougher domains.

It turns out that extending the well-posedness results for the Neumann problem beyond the case of the bilaplacian is an excruciatingly difficult problem, even if one considers only fourth-order operators with constant coefficients. Even defining Neumann boundary values for more general operators is a difficult problem (see Section~\ref{sec:neumann:var}), and while some progress has been made (see \cite{Agr07,MitM13A,BarHM15p,Bar15p}, or Section~\ref{sec:neumann:var} below), at present there are no well-posedness results for the Neumann problem with $L^p$ boundary data.

In analogy to \eqref{eqn:biharmonicDN} and \eqref{eqn:GrFla},
one can write
\begin{equation}
\label{eqn:generalDN}
\int_\Omega w\, Lu =A[u,w]+\int_{\partial\Omega} w\,K_A u-\partial_\nu w\,M_A u \,d\sigma,
\end{equation}
where $A[u,w]=\sum_{\abs{\alpha}=\abs{\beta}=2} a_{\alpha \beta} \int_{\Omega} D^\beta u\, D^{\alpha} w$ is an energy form associated to the operator $L=\sum_{\abs{\alpha}=\abs{\beta}=2} a_{\alpha \beta} D^{\alpha} D^{\beta}$. Note that in the context of fourth-order operators, the pair $(w, \partial_\nu w)$ constitutes the Dirichlet data for $w$ on the boundary, and so one can say that the operators $K_A u$ and $M_A u$ define the Neumann data for~$u$. One immediately faces the problem that the same higher-order operator $L$ can be rewritten in many different ways and gives rise to different energy forms. The corresponding Neumann data will be different. (This is the reason why there is a family of Neumann data for the biharmonic operator.)

Furthermore, whatever the choice of the form, in order to establish well-posedness of the Neumann problem, one needs to be able to estimate all second derivatives of a solution on the boundary in terms of the Neumann data. In the analogous second-order case, such an estimate is provided by the Rellich identity, which shows that the tangential derivatives are equivalent to the normal derivative in $L^2$  for solutions of elliptic PDEs. In the higher-order scenario, such a result calls for certain coercivity estimates which are still rather poorly understood. We refer the reader to \cite{Ver10} for a detailed discussion of related results and problems.

\subsection[Inhomogeneous problems for the biharmonic equation]{Inhomogeneous problems and other classes of boundary data}
\label{sec:inhom}

In \cite{AdoP98}, Adolfsson and Pipher investigated the inhomogeneous Dirichlet problem for the biharmonic equation with data in Besov and Sobolev spaces. While resting on the results for homogeneous boundary value problems discussed in Sections~\ref{sec:L2dirichlet} and~\ref{sec:regularity},
such a framework presents a completely new setting, allowing for the inhomogeneous problem and for consideration of classes of boundary data which are, in some sense, intermediate between the Dirichlet and the regularity problems.

They showed that if $\dot f\in  W\!A^p_{1+s}(\partial\Omega)$ and $h\in L^p_{s+1/p-3}(\Omega)$, then there exists a unique function $u$ that satisfies
\begin{equation}
\label{eqn:biinhomDirich}
\left\{
\begin{aligned}
\Delta^2 u&=h &&\text{in }\Omega,\\
\Tr \partial^\alpha u &=f_\alpha, &&\text{for }0\leq\abs{\alpha}\leq 1
\end{aligned}\right.
\end{equation}
subject to the estimate
\begin{equation}
\label{eqn:estimate:AdoP98}
\doublebar{u}_{L^p_{s+1/p+1}(\Omega)} \leq C\doublebar{h}_{L^p_{s+1/p-3}(\Omega)}
+C\doublebar{\dot f}_{W\!A^p_{1+s}(\partial\Omega)}
\end{equation}
provided $2-\varepsilon<p<2+\varepsilon$ and $0<s<1$.
Here
$\Tr w$ denotes the trace of $w$ in the sense of Sobolev spaces; that these may be extended to functions $u\in L^p_{s+1+1/p}$, $s>0$, was proven in \cite[Theorem~1.12]{AdoP98}.

In Lipschitz domains contained in~$\RR^3$, they proved these results for a broader range of $p$ and~$s$, namely for $0<s<1$ and for
\begin{equation}
\label{eqn:psrange}
\max\left(1,\frac{2}{s+1+\varepsilon}\right)<p<
\begin{cases}\infty,&s<\varepsilon,\\
\frac{2}{s-\varepsilon},&\varepsilon\leq s<1.\end{cases}
\end{equation}
Finally, in $C^1$ domains, they proved these results for any $p$ and $s$ with $1<p<\infty$ and $0<s<1$.

In \cite{MitMW11}, I.~Mitrea, M.~Mitrea and Wright extended the three-dimen\-sional results to $p=\infty$ (for $0<s<\varepsilon$) or $2/(s+1+\varepsilon)<p\leq 1$ (for $1-\varepsilon<s<1$). They also extended these results to data $h$ and $\dot f$ in more general Besov or Triebel-Lizorkin spaces.

In \cite{MitM13B}, I.~Mitrea and M.~Mitrea extended the results of \cite{AdoP98} to higher dimensions. That is, they showed that if $\Omega\subset\R^n$ is a Lipschitz domain and $n\geq 4$, then there is a unique solution to the problem~\eqref{eqn:biinhomDirich} subject to the estimate \eqref{eqn:estimate:AdoP98} provided that $0<s<1$ and that
\begin{equation*}
\max\left(1,\frac{n-1}{s+(n-1)/2+\varepsilon}\right) < p <
\begin{cases}\infty,&(n-3)/2+s<\varepsilon,\\
\frac{n-1}{(n-3)/2+s-\varepsilon},&\varepsilon\leq s<1.\end{cases}
\end{equation*}
As in \cite{MitMW11}, their results extend to more general function spaces.

I.~Mitrea and M.~Mitrea also showed that, for the same values of $p$ and~$s$, there exist unique solutions to the inhomogeneous \emph{Neumann} problem
\begin{equation*}
\left\{
\begin{aligned}
\Delta^2 u&=h &&\text{in }\Omega,\\
M_\rho u &=f &&\text{on }\partial\Omega,\\
K_\rho u &=\Lambda &&\text{on }\partial\Omega\\
\end{aligned}\right.
\end{equation*}
where $M_\rho$ and $K_\rho$ are as in Section~\ref{sec:neumann}, subject to the estimate
\begin{equation}
\doublebar{u}_{L^p_{s+1/p+1}(\Omega)}
\leq C\doublebar{h}_{L^p_{s+1/p-3}(\Omega)}
+C\doublebar{f}_{B^{p,p}_{s-1}(\partial\Omega)}
+C\doublebar{\Lambda}_{B^{p,p}_{s-2}(\partial\Omega)}.
\end{equation}

Finally, in \cite[Section~6.4]{MitM13A}, I.~Mitrea and M.~Mitrea proved similar results for more general constant-coefficient elliptic operators. That is, if $L$ is an operator given by formula~\eqref{eqn:constsystem} whose coefficients satisfy the ellipticity condition~\eqref{eqn:constelliptic}, and if $\Omega\subset\R^n$ is a bounded Lipschitz domain, then there exists a unique solution to the Dirichlet problem
\begin{equation}
\label{eqn:inhomogeneous}
\left\{
\begin{aligned}
L u&=h &&\text{in }\Omega,\\
\Tr \partial^\alpha u &=f_\alpha, &&\text{for }0\leq\abs{\alpha}\leq 1
\end{aligned}\right.
\end{equation}
subject to the estimate
\begin{equation}
\label{eqn:inhomogeneous:estimate}
\doublebar{u}_{B^{p,q}_{m-1+s+1/p}(\Omega)} \leq C\doublebar{h}_{B^{p,q}_{-m-1+s+1/p}(\Omega)}
+C\doublebar{\dot f}_{W\!A^{p,q}_{m-1+s}(\partial\Omega)}
\end{equation}
provided $2-\varepsilon<p<2+\varepsilon$, $2-\varepsilon<q<2+\varepsilon$, and $1/2-\varepsilon<s<1/2+\varepsilon$. (A very similar result is valid for variable-coefficient divergence-form operators; see \cite{BreM13}, discussed as formula~\eqref{eqn:BreM13} below.) Furthermore, with a slightly stronger ellipticity condition
\begin{equation*}\re \sum_{\abs{\alpha}=\abs{\beta}=m}\sum_{j,k=1}^\ell a_{\alpha\beta}^{jk} \zeta_j^\alpha \overline{\zeta_k^\beta} \geq \lambda \abs{\zeta}^2,\end{equation*}
they established well-posedness of the inhomogeneous Neumann problem for arbitrary constant-coefficient operators.
See Section~\ref{sec:neumann:var} below for a formulation of Neumann boundary data in the case of arbitrary operators.
Finally, if $L$ is self-adjoint and $n>2m$, and if the unit outward normal $\nu$ to $\partial\Omega$ lies in $VMO(\partial\Omega)$, then the Dirichlet problem \eqref{eqn:inhomogeneous} has a unique solution satisfying the estimate~\eqref{eqn:inhomogeneous:estimate} for any $0<s<1$ and any $1<p<\infty$, $1<q<\infty$.

Let us define the function spaces appearing above. $L^p_\alpha(\RR^n)$ is defined to be $\{g:(I-\Delta)^{\alpha/2} g\in L^p(\RR^n)\}$; we say $g\in L^p_\alpha(\Omega)$ if $g=h\big\vert_{\Omega}$ for some $h\in L^p_\alpha(\RR^n)$. If $k$ is a nonnegative integer, then $L^p_k=W^p_k$.
If $m$ is an integer and $0<s<1$, then the Whitney-Besov space $W\!A^p_{m-1+s}=W\!A^{p,p}_{m-1+s}$ or $W\!A^{p,q}_{m-1+s}$ is defined analogously to $W\!A^p_{m}$ (see Definition~\ref{dfn:WhitneySobolev}), except that we take the completion with respect to the Whitney-Besov norm
\begin{equation}
\label{eqn:whitneybesovnorm}
\sum_{\abs{\alpha}\leq m-1} \doublebar{\partial^\alpha \psi}_{L^p(\partial\Omega)}
+\sum_{\abs{\alpha}=m-1} \doublebar{\partial^\alpha\psi}_{B^{p,q}_s(\partial\Omega)}
\end{equation}
rather than the Whitney-Sobolev norm
\begin{equation*}
\sum_{\abs{\alpha}\leq m-1} \doublebar{\partial^\alpha \psi}_{L^p(\partial\Omega)}
+\sum_{\abs{\alpha}=m-1} \doublebar{\nabla_\tau \partial^\alpha\psi}_{L^p(\partial\Omega)}.
\end{equation*}

In \cite{AdoP98}, the general problem \eqref{eqn:biinhomDirich} for $\Delta^2$ was first reduced to the case $h=0$ (that is, to a homogeneous problem) by means of trace/extension theorems, that is, subtracting $w(X)=\int_{\RR^n} F(X,Y)\,\tilde h(Y)\,dY$, and showing that if $h\in L^p_{s+1/p-3}(\Omega)$ then $(\Tr w,\Tr \nabla w)\in W\!A^p_{1+s}(\partial\Omega)$. Next, the well-posedness of Dirichlet and regularity problems discussed in Sections~\ref{sec:L2dirichlet} and~\ref{sec:regularity} provide the endpoint cases $s=0$ and $s=1$, respectively. The core of the matter is to show that, if $u$ is biharmonic, $k$ is an integer and $0\leq \alpha\leq 1$, then $u\in L^p_{k+\alpha}(\Omega)$ if and only if
\begin{equation}\label{eqn:equivestimates}
\int_\Omega \abs{\nabla^{k+1} u(X)}^p \dist(X,\partial\Omega)^{p-p\alpha} +\abs{\nabla^k u(X)}^p +\abs{u(X)}^p\,dX<\infty,
\end{equation}
(cf.\ \cite[Proposition~S]{AdoP98}). With this at hand, one can use square-function estimates to justify the aforementioned endpoint results. Indeed, observe that for $p=2$ the first integral on the left-hand side of \eqref{eqn:equivestimates}  is exactly the $L^2$ norm of $S(\nabla^k u)$. The latter, by \cite{PipV91} (discussed in Section~\ref{sec:square}), is equivalent to the $L^2$ norm of the corresponding non-tangential maximal function, connecting the estimate \eqref{eqn:estimate:AdoP98} to the nontangential estimates in the Dirichlet problem~\eqref{eqn:bidirichlet} and the regularity problem~\ref{eqn:polyregularity}. Finally, one can build an interpolation-type scheme to pass to well-posedness in intermediate Besov and Sobolev spaces.

The solution in \cite{MitM13A} to the problem~\eqref{eqn:inhomogeneous}, at least in the case of general Lipschitz domains, was constructed in the opposite way, by first reducing to the case where the boundary data $\arr f=0$. Using duality it is straightforward to establish well-posedness in the case $p=q=2$, $s=1/2$; perturbative results then suffice to extend to $p$, $q$ near $2$ and $s$ near $1/2$.

\section{Boundary value problems with variable coefficients}
\label{sec:var}

Results for higher order differential equations with variable coefficients are very scarce. As we discussed in Section~\ref{sec:dfn}, there are two natural manifestations of higher-order operators with variable coefficients. Operators in divergence form arise via the weak formulation framework.  Conversely, operators in composition form generalize the bilaplacian under a pull-back of a Lipschitz domain to the upper half-space.

Both classes of operators have been investigated. However, operators in divergence form have received somewhat more study; thus, we begin this section by reviewing the definition of divergence-form operator.
A divergence-form operator $L$, acting on $W^{2}_{m,loc}(\Omega\mapsto \CC^\ell)$, may be defined weakly via~\eqref{eqn:weaksoln}; we say that $Lu=h$ if
\begin{equation}\label{eqn:weaksoln2}
\sum_{j=1}^\ell \int_\Omega \varphi_j \,h_j
=
(-1)^m\sum_{j,k=1}^\ell \sum_{\abs{\alpha}=\abs{\beta}=m} \int_\Omega \partial^\alpha \varphi_j \,a_{\alpha\beta}^{jk}\, \partial^\beta u_k
\end{equation}
for all $\varphi$ smooth and compactly supported in~$\Omega$.

\subsection{The Kato problem and the Riesz transforms}
\label{sec:kato}

We begin with the Kato problem and the properties of the Riesz transform; this is an important topic in elliptic theory, which formally stands somewhat apart from the well-posedness issues.

Suppose that $L$ is a variable-coefficient operator in divergence form, that is, an operator defined by~\eqref{eqn:weaksoln2}. Suppose that $L$ satisfies the bound
\begin{equation}
\label{eqn:weakbdd}
\biggabs{\sum_{\abs{\alpha}=\abs{\beta}=m} \sum_{j,k=1}^\ell\int_{\RR^n} a_{\alpha\beta}^{jk} \,\partial^\beta\! f_k\,\partial^\alpha g_j}
\leq C
\doublebar{\nabla^m f}_{L^2(\RR^n)}
\doublebar{\nabla^m g}_{L^2(\RR^n)}
\end{equation}
for all $f$ and $g$ in $\dot W^2_m(\RR^n)$,
and the ellipticity estimate
\begin{equation}\label{eqn:varelliptic}
\re\sum_{\abs{\alpha}=\abs{\beta}=m } \sum_{j,k=1}^\ell
\int_\Omega
 a_{\alpha\beta}^{jk}(X) \partial^\beta \varphi_k(X) \partial^\alpha \overline{\varphi_j}(X)
\,dX\geq \lambda\sum_{\abs{\alpha}=m} \sum_{k=1}^\ell
\int_\Omega \abs{\partial^\alpha \varphi_k}^2
\end{equation}
for all functions $\varphi\in C^\infty_0(\Omega\mapsto \CC^\ell)$. Notice that this is a weaker requirement than the pointwise ellipticity condition~\eqref{eqn:varelliptic1}.

Auscher, Hofmann, McIntosh and Tchamitchian \cite{AusHMT01} proved that under these conditions, the Kato estimate
\begin{equation}
\label{eqn:Kato}
\frac{1}{C}\doublebar{\nabla^m f}_{L^2(\RR^n)}\leq\doublebar{\sqrt{L} f}_{L^2(\RR^n)}\leq C\doublebar{\nabla^m f}_{L^2(\RR^n)}
\end{equation}
is valid for some constant~$C$.
They also proved similar results for operators with lower-order terms.

It was later observed in \cite{Aus04} that by the methods of \cite{AusT98}, if $1\leq n\leq 2m$, then the bound on the Riesz transform $\nabla^m L^{-1/2}$ in $L^p$ (that is, the first inequality in formula~\eqref{eqn:Kato}) extends to the range $1<p<2+\varepsilon$, and the reverse Riesz transform bound (that is, the second inequality in formula~\eqref{eqn:Kato}) extends to the range $1<p<\infty$. This also holds if the Schwartz kernel $W_{t}(X,Y)$ of the operator $e^{-tL}$ satisfies certain pointwise bounds (e.g., if $L$ is second-order and the coefficients of $A$ are real).

In the case where $n>2m$, the inequality
$\doublebar{\nabla^m L^{-1/2} f}_{L^p(\RR^n)}\leq C\doublebar{f}_{L^p(\RR^n)}$ holds for $2n/(n+2m)-\varepsilon<p\leq 2$; see \cite{BluK04,Aus04}. The reverse inequality holds for $\max(2n/(n+4m)-\varepsilon,1)<p<2$ by \cite[Theorem~18]{Aus04}, and for $2<p<{2n}/({n-2m})+\varepsilon$ by duality (see \cite[Section~7.2]{Aus07}).

In the case of second-order operators, the Kato estimate implies well-posedness of boundary-value problems with $L^2$ data in the upper half space for certain coefficients in a special (``block'') form. We conjecture that the same is true in the case of higher-order operators; see Section~\ref{sec:open}.

\subsection{The Dirichlet problem for operators in divergence form}
\label{sec:variablediv}

In this section we discuss boundary-value problems for divergence-form operators with variable coefficients. At the moment, well-posedness results for such operators are restricted in that the boundary problems treated fall {\it strictly} between the range of $L^p$-Dirichlet and $L^p$-regularity, in the sense of Section~\ref{sec:inhom}. That is, there are at present no well-posedness results for the $L^p$-Dirichlet, regularity, and Neumann problems on Lipschitz domains with the usual sharp estimates in terms of the non-tangential maximal function for these divergence-form operators. (Such problems are now being considered; see Section~\ref{sec:open}.)

To be more precise, recall from the discussion in Section~\ref{sec:inhom} that the classical Dirichlet and regularity problems, with boundary data in~$L^p$, can be viewed as the $s=0$, $1$ endpoints of the boundary problem studied in \cite{AdoP98}, \cite{MitMW11}  and \cite{MitM13B}
\begin{align*}
\Delta^2 u&=h\text{ in }\Omega,\quad\partial^\alpha u\big\vert_{\partial\Omega}
=f_\alpha\text{ for all }\abs{\alpha}\leq 1
\end{align*}
with $\dot f$ lying in an \emph{intermediate} smoothness space $W\!A^p_{1+s}(\partial\Omega)$, $0\leq s\leq 1$.
In the context of divergence-form higher-order operators with variable coefficients, essentially the known results pertain \emph{only} to boundary data of intermediate smoothness.

In \cite{Agr07}, Agranovich investigated the inhomogeneous Dirichlet problem, in Lipschitz domains, for such operators $L$ that are elliptic (in the pointwise sense of \eqref{eqn:varelliptic1}, and not the more general condition \eqref{eqn:varelliptic}) and whose coefficients $a_{\alpha\beta}^{jk}$ are Lipschitz continuous in~$\Omega$.

He showed that if $h\in L^p_{-m-1+1/p+s}(\Omega)$ and $\dot f\in W\!A^{p}_{m-1+s}(\partial\Omega)$, for some $0<s<1$, and if $\abs{p-2}$ is small enough, then the Dirichlet problem
\begin{equation}
\label{eqn:inhom:L}
\left\{\begin{aligned}
Lu &=h && \text{in }\Omega,\\
\Tr \partial^\alpha u &= f_\alpha && \text{for all } 0\leq \abs{\alpha}\leq m-1
\end{aligned}\right.
\end{equation}
has a unique solution $u$ that satisfies the estimate
\begin{equation}
\label{eqn:estimate:Agr07}
\doublebar{u}_{L^p_{m-1+s+1/p}(\Omega)}
\leq
C  \doublebar{h}_{L^p_{-m-1+1/p+s}(\partial\Omega)}
+
C  \doublebar{\dot f}_{W\!A^p_{m-1+s}(\partial\Omega)}.
\end{equation}

Agranovich also considered the Neumann problem for such operators. As we discussed in Section~\ref{sec:neumann}, defining the Neumann problem is a delicate matter. In the context of zero boundary data, the situation is a little simpler as one can take a formal functional analytic point of view and avoid to some extent the discussion of estimates at the boundary.
We say that $u$ solves the Neumann problem for $L$, with homogeneous boundary data, if the equation \eqref{eqn:weaksoln2} in the weak formulation of $L$ is valid for all test functions $\varphi$ compactly supported in $\RR^n$ (but not necessarily in~$\Omega$.)
Agranovich showed that, if $h\in \mathring L^p_{-m-1+1/p+s}(\Omega)$, then there exists a unique function $u\in L^p_{m-1+1/p+s}(\Omega)$ that solves this Neumann problem with homogeneous boundary data, under the same conditions on $p$, $s$, $L$ as for his results for the Dirichlet problem.
Here $h\in\mathring L^p_{\alpha}(\Omega)$ if $h=g\big\vert_{\Omega}$ for some $g\in L^p_\alpha(\RR^n)$ that in addition is supported in~$\bar\Omega$.

In \cite{MazMS10}, Maz'ya, M.~Mitrea and Shaposhnikova considered the Dirichlet problem, again with boundary data in intermediate Besov spaces, for much rougher coefficients.
They showed that if $f\in W\!A^p_{m-1+s}(\partial\Omega)$, for some $0<s<1$ and some $1<p<\infty$, if $h$ lies in an appropriate space, and if $L$ is a divergence-form operator of order~$2m$ (as defined by \eqref{eqn:weaksoln2}), then under some conditions, there is a unique function $u$ that satisfies the Dirichlet problem
\eqref{eqn:inhom:L} subject to
the estimate
\begin{equation}
\label{eqn:MMSsquare}
\doublebar{u}_{W^p_{m,1-s-1/p}} =
\biggl(\sum_{\abs{\alpha}\leq m} \int_\Omega \abs{\partial^\alpha u(X)}^p \dist(X,\partial\Omega)^{p-ps-1}\,dX\biggr)^{1/p}<\infty.
\end{equation}
See \cite[Theorem~8.1]{MazMS10}.
The inhomogeneous data $h$ is required to lie in the space $V^p_{-m,1-s-1/p}(\Omega)$, the dual space to $V^{q}_{m,s+1/p-1}(\Omega)$, where
\begin{equation}
\label{eqn:MMSdual}
\doublebar{w}_{V^p_{m,a}}=\biggl(\sum_{\abs{\alpha}\leq m} \int_\Omega \abs{\partial^\alpha u(X)}^p \dist(X,\partial\Omega)^{pa+p\abs{\alpha}-pm} \,dX\biggr)^{1/p}.\end{equation}
Notice that $w\in V^p_{m,a}$ if and only if $w\in W^p_{m,a}$ and $\partial^\alpha w=0$ on $\partial\Omega$ for all $0\leq \abs{\alpha}\leq m-1$.

The conditions are that the coefficients $a_{\alpha\beta}^{jk}$ satisfy the weak ellipticity condition \eqref{eqn:varelliptic} considered in the theory of the Kato problem, that $\Omega$ be a Lipschitz domain whose normal vector $\nu$ lies in $VMO(\partial\Omega)$, and that the coefficients $a_{\alpha\beta}^{ij}$ lie in $L^\infty(\RR^n)$ and in $VMO(\RR^n)$. Recall that this condition on $\Omega$ has also arisen in \cite{MitM11} (it ensures the validity of  formula~\eqref{eqn:MM11}). Notice that the $L^\infty$ bound on the coefficients is a stronger condition than the bound \ref{eqn:weakbdd} of \cite{AusHMT01}, and the requirement that the coefficients lie in $VMO(\RR^n)$ is a regularity requirement that is weaker than the requirement of \cite{Agr07} that the coefficients be Lipschitz continuous.

In fact, \cite{MazMS10} provides a more intricate  result, allowing one to deduce a well-posedness range of $s$ and~$p$, given information about the oscillation of the coefficients~$a_{\alpha\beta}^{jk}$ and the normal to the domain~$\nu$. In the extreme case, when the oscillations for both are vanishing, the allowable range expands to $0<s<1$, $1<p<\infty$, as stated above.

The construction of solutions to the Dirichlet problem may be simplified using trace and extension theorems. In \cite[Proposition~7.3]{MazMS10}, the authors showed that if $\dot f\in W\!A^p_{m-1+s}(\partial\Omega)$, then there exists a function $F\in W^p_{m,a}$ such that $\partial^\alpha F=f_\alpha$ on $\partial\Omega$.
It is easy to see that if $F\in W^p_{m,a}$, and the coefficients $a_{\alpha\beta}^{jk}$ of $L$ are bounded pointwise, then $LF\in V_{-m,1-s-1/p}$. Thus, the Dirichlet problem
\begin{equation*}
\left\{\begin{aligned}
Lu&=h \text{ in }\Omega,
\\
\partial^\alpha u\big\vert_{\partial\Omega}&=f_\alpha, \text{ for all $\abs{\alpha}\leq m-1$},
\\
\doublebar{u}_{W^p_{m,1-s-1/p}}&\leq  \doublebar{h}_{V^p_{-m,1-s-1/p}}+\doublebar{\dot f}_{W\!A^p_{m-1+s}(\partial\Omega)}
\end{aligned}\right.
\end{equation*}
may be solved by solving the Dirichlet problem with homogeneous boundary data
\begin{equation*}
\left\{\begin{aligned}
Lw&=h-LF \text{ in }\Omega,
\\
\partial^\alpha w\big\vert_{\partial\Omega}&=0 \text{ for all $\abs{\alpha}\leq m-1$},
\\
\doublebar{w}_{W^p_{m,1-s-1/p}}&\leq \doublebar{h}_{V^p_{-m,1-s-1/p}}+\doublebar{LF}_{V^p_{-m,1-s-1/p}}
\end{aligned}\right.
\end{equation*}
for some extension~$F$ and then letting $u=w+F$.

Some limited results are available in the case where the coefficients $a_{\alpha\beta}^{jk}(x)$ satisfy no smoothness assumptions whatsoever. By \cite[Theorem~5.1]{BreM13}, the Dirichlet problem \eqref{eqn:inhom:L}, in a Lipschitz domain~$\Omega$ (whose unit outward normal need not be in~$VMO$), with data in appropriate spaces, has a unique solution $u$ that satisfies the estimate
\begin{equation}\label{eqn:BreM13}\doublebar{u}_{W^p_{m,1-s-1/p}}\leq C \doublebar{\arr f}_{W\!A^p_{m-1+s}(\partial\Omega)} + C \doublebar{\arr h}_{V^p_{-m,1-s-1/p}(\Omega)}\end{equation}
provided that $\abs{p-2}$ and $\abs{s-1/2}$ are small enough.

We comment on the estimate \eqref{eqn:MMSsquare}. First, by \cite[Propositon~S]{AdoP98} (listed above as formula~\eqref{eqn:equivestimates}), if $u$ is biharmonic then the estimate \eqref{eqn:MMSsquare} is equivalent to the estimate \eqref{eqn:estimate:Agr07} of \cite{Agr07}. Second, by \eqref{eqn:squarehigher}, if the coefficients $a_{\alpha\beta}^{jk}$ are constant,  one can draw connections between \eqref{eqn:MMSsquare} for $s=0$,~$1$ and the nontangential maximal estimates of the Dirichlet or regularity problems \eqref{eqn:systemdirichlet} or~\eqref{eqn:systemregularity}. However, as we pointed out earlier, this  endpoint case, corresponding to the true $L^p$-Dirichlet and regularity problems, has not been achieved.

\subsection{The Dirichlet problem for operators in composition form}
\label{sec:variablecomp}

Let us now discuss variable-coefficient fourth-order operators in composition form. Recall that this particular form arises naturally when considering the transformation of the bilaplacian under a pull-back from a Lipschitz domain (cf.\ \eqref{eqn:prodform}). The authors of the present paper have shown the well-posedness, for a class of such operators, of the Dirichlet problem with boundary data in~$L^2$, thus establishing the first results concerning the $L^p$-Dirichlet problem for variable-coefficient higher-order operators.

Consider the Dirichlet problem
\begin{equation}\label{eqn:Bdirichlet}
\left\{\begin{aligned}
L^*(aLu)&=0&& \text{in }\Omega,\\
 u&=f&&\text{on }\partial\Omega,\\
\nu\cdot A\nabla u&=g&&\text{on }\partial\Omega,\\
\doublebar{\tilde N(\nabla u)}_{L^2(\partial\Omega)}
&\leq C\doublebar{\nabla f}_{L^2(\partial\Omega)} + C\doublebar{g}_{L^2(\partial\Omega)}.
\negphantom{\text{on }\partial\Omega.}
\end{aligned}\right.
\end{equation}
Here $L$ is a \emph{second-order} divergence form differential operator $L=-\Div A(X)\nabla$, and $a$ is a scalar-valued function.
(For rough coefficients $A$, the exact weak definition of $L^*(aLu)=0$ is somewhat delicate, and so we refer the reader to \cite{BarM13}.)
The domain $\Omega$ is taken to be the domain above a Lipschitz graph, that is, $\Omega=\{(x,t):x\in\RR^{n-1},\>t>\varphi(x)\}$ for some function $\varphi$ with $\nabla\varphi\in L^\infty(\RR^{n-1})$.
As pointed out above, the class of equations $L^*(aLu)=0$ is preserved by a change of variables, and so well-posedness of the Dirichlet problem~\eqref{eqn:Bdirichlet} in such domains follows from well-posedness in upper half-spaces $\RR^n_+$. Hence, in the remainder of this section, $\Omega=\RR^n_+$.

The appropriate ellipticity condition is then
\begin{equation}
\label{eqn:Belliptic}
\lambda\leq a(X)\leq \Lambda,
\qquad \lambda\abs{\eta}^2\leq \re\overline{\eta}^t A(X)\eta,\quad \abs{A(X)}\leq\Lambda
\end{equation}
for all $X\in\RR^n$ and all $\eta\in\CC^n$, for some constants $\Lambda>\lambda>0$. The modified nontangential maximal function $\tilde N(\nabla u)$, defined by
\begin{equation*}\tilde N(\nabla u)(Q)=\sup_{X\in\Gamma(Q)} \biggl(\fint_{B(X,\dist(X,\partial\Omega)/2)} \abs{\nabla u}^2\biggr)^{1/2},\end{equation*}
is taken from \cite{KenP93} and is fairly common in the study of variable-coefficient elliptic operators.

In this case, we say that $u\big\vert_{\partial\Omega}=f$ and $\nu\cdot A\nabla u=g$ if
\begin{align*}
\lim_{t\to 0^+} \doublebar{u(\,\cdot\,+t\e)-f}_{W^2_1(\partial\Omega)}&=0,
\\
\lim_{t\to 0^+} \doublebar{\nu\cdot A\nabla u(\,\cdot\,+t\e)-g}_{L^2(\partial\Omega)}&=0\end{align*}
where $\e=\e_{n}$ is the unit vector in the vertical direction. Notice that by the restriction on the domain~$\Omega$, $\e$ is transverse to the boundary at all points. We usually refer to the vertical direction as the $t$-direction, and if some function depends only on the first $n-1$ coordinates, we say that function is $t$-independent.

In \cite{BarM13}, the authors of the present paper have shown that if $n\geq 3$, and if $a$ and $A$ satisfy \eqref{eqn:Belliptic} and are $t$-independent, then for every $f\in W^{2}_1(\partial\Omega)$ and every $g\in L^2(\partial\Omega)$, there exists a $u$ that satisfies \eqref{eqn:Bdirichlet}, provided that the second order operator $L=\Div A\nabla$ is good from the point of view of the second order theory.

Without going into the details, we mention that there are certain restrictions on the coefficients $A$ necessary to ensure the well-posedness even of the of the corresponding second-order boundary value problems; see \cite{CafFK81}. The key issues are good behavior in the direction transverse to the boundary, and symmetry. See \cite{JerK81A,KenP93} for results for symmetric $t$-independent coefficients,
\cite{KenKPT00,KenR09,Rul07,HofKMP12p, HofKMP13p} for well-posedness results and important counterexamples for non-symmetric coefficents, and
\cite{AusAH08,AusAM10,AlfAAHK11} for perturbation results for $t$-independent coefficients.

In particular, using the results of \cite{AusAH08,AusAM10,AlfAAHK11}, we have established that the $L^2$-Dirichlet problem \eqref{eqn:Bdirichlet} in the upper half-space is well-posed, provided the coefficients $a$ and $A$ satisfy \eqref{eqn:Belliptic} and are $t$-independent, if in addition one of the following conditions holds:
\begin{enumerate}
\item \label{BcondA} The matrix $A$ is real and symmetric,
\item \label{BcondB} The matrix $A$ is constant,
\item \label{BcondC} The matrix $A$ is in block form (see Section~\ref{sec:open})
and the Schwartz kernel $W_{t}(X,Y)$ of the operator $e^{-tL}$ satisfies certain pointwise bounds, or
\item There is some matrix $A_0$, satisfying \eqref{BcondA}, \eqref{BcondB} or~\eqref{BcondC}, that again satisfies \eqref{eqn:Belliptic} and is $t$-independent, such that $\doublebar{A-A_0}_{L^\infty(\RR^{n-1})}$ is small enough (depending only on the constants $\lambda$, $\Lambda$ in \eqref{eqn:Belliptic}).
\end{enumerate}

The solutions to \eqref{eqn:Bdirichlet} take the following form. Inspired by formula~\eqref{eqn:DKV} (taken from \cite{DahKV86}), and a similar representation in \cite{PipV92}, we let
\begin{equation}
\label{eqn:BarM14}
\mathcal{E} h=\int_{\Omega}F(X,Y)\frac{1}{a(Y)}\partial_n^2 \s_* h(Y)\,dY\end{equation}
for $h$ defined on~$\partial\Omega$, where $\s_*$ is the (second-order) single layer potential associated to $L^*$ and $F$ is the fundamental solution associated to~$L$. Then $a(X) \, L(\mathcal{E} h)(X)=\partial_n^2 \s_* f(X)$ in~$\Omega$ (and is zero in its complement); if $A^*$ is $t$-independent, then $L^*(\partial_n^2 \s_* h)=\partial_n^2 L^*(\s_* h)=0$. Thus $u=w+\mathcal{E} h$ is a solution to \eqref{eqn:Bdirichlet}, for any solution $w$ to $Lw=0$.
The estimate
$\doublebar{\tilde N(\nabla \mathcal{E} h)}_{L^2(\partial\Omega)}\leq
\doublebar{h}_{L^2(\partial\Omega)}$ must then be established. In the case of biharmonic functions (considered in \cite{PipV92}), this estimate follows from the boundedness of the Cauchy integral; in the case of \eqref{eqn:Bdirichlet}, this is the most delicate part of the construction, as the operators involved are far from being Calder\'on-Zygmund kernels.
Once this estimate has been established, it can be shown, by an argument that precisely parallels that of \cite{PipV92}, that there exists a $w$ and $h$ such that $Lw=0$ and $u=w+\mathcal{E}h$ solves~\eqref{eqn:Bdirichlet}.

\subsection{The fundamental solution}
\label{sec:fundamental}

A set of important tools, and interesting objects of study in their own right, are the fundamental solutions and Green's functions of differential operators in various domains. To mention some applications presented in this survey, recall from Sections~\ref{sec:pointwise} and~\ref{sec:Green} that bounds on Green's functions $G$ are closely tied to maximum principle estimates, and from Sections~\ref{sec:wiener} and~\ref{sec:wiener2} that the fundamental solution $F$ is used to establish regularity of boundary points (that is, the Wiener criterion). See in particular Theorem~\ref{thm:posweightF}.

Furthermore, fundamental solutions and Green's functions are often crucial elements of the construction of solutions to boundary-value problems. In the case of boundary-value problems in divergence form, the fundamental solution or Green's function for the corresponding higher-order operators are often useful; see the constructions in \cite{SelS81,CohG83,PipV95A,PipV95B,Ver96,MazMS10}, or in formulas \eqref{eqn:SelS81} and \eqref{eqn:CG83} above.
In the case of operators in composition form, it is often more appropriate to use the fundamental solution for the lower-order components; see, for example, formulas~\eqref{eqn:DKV}, \eqref{eqn:PipV93} and~\eqref{eqn:BarM14}, or the paper \cite{Ver90}, which makes extensive use of the Green's function for $(-\Delta)^m$ to solve boundary-value problems for $(-\Delta)^{m+1}$.

We now discuss some constructions of the fundamental solution.
In the case of the biharmonic equation, and more generally in the case of constant-coefficient equations, the fundamental solution may be found in a fairly straightforward fashion, for example, by use of the Fourier transform; see, for example, formulas~\eqref{eqp8.2} and \eqref{eqp8.3} above, \cite{Sha45,Mor54,Joh55,OrtW96,Hor03} (the relevant results of which are summarized as \cite[Theorem~4.2]{MitM13A}),
or \cite{Dal13, DalMM13}.
In the case of variable-coefficient second-order operators, the fundamental solution has been constructed in \cite{LitSW63,GruW82,KenN85,Fuc86,DolM95,HofK07,Ros13} under progressively weaker assumptions on the operators. The most recent of these papers, \cite{Ros13}, constructs the fundamental solution $F$ for a second-order operator~$L$ under the assumption that if $Lu=0$ in some ball $B(X,r)$, then we have the local boundedness estimate
\begin{equation}\label{eqn:Moser}\abs{u(X)}\leq C\biggl(\frac{1}{r^n} \int_{B(X,r)}\abs{u}^2\biggr)^{1/2}\end{equation}
for some constant~$C$ depending only on~$L$ and not on $u$, $X$ or~$r$.
This assumption is true if $L$ is a scalar second-order operator with real coefficients (see \cite{Mos61}) but is not necessarily true for more general elliptic operators (see \cite{Fre08}).

If $Lu=0$ in $B(X,r)$ for some elliptic operator of order $2m$, where $2m>n$, then the local boundedness estimate \eqref{eqn:Moser} follows from the Poincar\'e inequality, the Caccioppoli inequality
\begin{equation}\label{eqn:Caccioppoli}\int_{B(X,r/2)} \abs{\nabla^m u}^2 \leq \frac{C}{r^{2m}}\int_{B(X,r)} \abs{u}^2\end{equation}
and Morrey's inequality
\begin{equation*}\abs{u(X)}\leq C\sum_{j=0}^{N} r^j\biggl(\frac{1}{r^n}\int_{B(X,r/2)}\abs{\nabla^j u}^2\biggr)^{1/2} \text{ whenever $N/2>n$}.\end{equation*}
Weaker versions of the Caccioppoli inequality \eqref{eqn:Caccioppoli} (that is, bounds with higher-order derivatives appearing on the right-hand side) were established in \cite{Cam80} and \cite{AusQ00}.
In \cite{Bar14p}, the first author of the present paper established the full Caccioppoli inequality \eqref{eqn:Caccioppoli}, thus establishing that if $2m>n$ then solutions to $Lu=0$ satisfy the estimate \eqref{eqn:Moser}. (Compare the results of Section~\ref{sec:pointwise}, in which solutions to $(-\Delta)^m u=f$ are shown to be pointwise bounded only if $2m>n-2$; as observed in that section, Morrey's inequality yields one fewer degree of smoothness but was still adequate for the purpose of \cite{Bar14p}.)

Working much as in the second-order papers listed above, Barton then constructed the fundamental solution for divergence-form differential operators $L$ with order $2m>n$ and with bounded coefficients satisfying the ellipticity condition~\eqref{eqn:varelliptic}. In the case of operators $L$ with $2m\leq n$, she then constructed an auxiliary operator $\widetilde L$ with $2m>n$ and used the fundamental solution for $\widetilde L$ to construct the fundamental solution for~$L$; this technique was also used in \cite{AusHMT01} to pass from operators of high order to operators of arbitrary order, and for similar reasons (i.e., to exploit pointwise bounds present only in the case of operators of very high order).

This technique allowed the proof of the following theorem, the main result of \cite{Bar14p}.
\begin{theorem}\label{thm:fundamental:high:2}
Let $L$ be a divergence-form operator of order~$2m$, acting on functions defined on~$\R^n$, that satisfies the ellipticity condition~\eqref{eqn:varelliptic} and whose coefficients~$A$ are pointwise bounded.
Then there exists an array of functions $F^L_{j,k}(x,y)$ with the following properties.

Let $q$ and $s$ be two integers that satisfy $q+s<n$ and the bounds $0\leq q\leq \min(m,n/2)$, $0\leq s\leq \min(m,n/2)$.

Then there is some $\varepsilon>0$ such that if $x_0\in\R^n$, if $0<4r<R$, if $A(x_0,R)=B(x_0,2R)\setminus B(x_0,R)$, and if $q<n/2$ then
\begin{equation}
\label{eqn:fundamental:far:low:2}
\int_{y\in B(x_0,r)}\int_{x\in A(x_0,R)} \abs{\nabla^{m-s}_x \nabla^{m-q}_y \mat F^L(x,y)}^2\,dx\,dy \leq C r^{2q} R^{2s} \biggl(\frac{r}{R}\biggr)^\varepsilon
.\end{equation}
If $q=n/2$ then we instead have the bound
\begin{equation}
\label{eqn:fundamental:far:lowest:2}
\int_{y\in B(x_0,r)}\int_{x\in A(x_0,R)} \abs{\nabla^{m-s}_x \nabla^{m-q}_y \mat F^L(x,y)}^2\,dx\,dy \leq C(\delta)\, r^{2q} R^{2s} \biggl(\frac{R}{r}\biggr)^\delta
\end{equation}
for all $\delta>0$ and some constant $C(\delta)$ depending on~$\delta$.

We also have the symmetry property
\begin{equation}
\label{eqn:fundamental:symmetric:low}
\partial_x^\gamma\partial_y^\delta F^L_{j,k}(x,y) = \overline{\partial_x^\gamma\partial_y^\delta F^{L^*}_{k,j}(y,x)}
\end{equation}
as locally $L^2$ functions, for all multiindices $\gamma$, $\delta$ with $\abs{\gamma}=m-q$ and $\abs{\delta}=m-s$.

If in addition $q+s>0$, then for all $p$ with $1\leq p\leq 2$ and $p<n/(n-(q+s))$, we have that
\begin{equation}
\label{eqn:fundamental:near}
\int_{B(x_0,r)}\int_{B(x_0,r)} \abs{\nabla^{m-s}_x \nabla^{m-q}_y \mat F^L(x,y)}^p\,dx\,dy \leq C(p)\, r^{2n+p(s+q-n)}\end{equation}
for all $x_0\in\R^n$ and all $r>0$.

Finally, there is some $\varepsilon>0$ such that if $2-\varepsilon<p<2+\varepsilon$ then $\nabla^m\Pi^L$ extends to a bounded operator $L^p(\R^n)\mapsto L^p(\R^n)$. If $\gamma$ satisfies $m-n/p<\abs{\gamma}\leq m-1$ for some such~$p$, then
\begin{equation}
\label{eqn:fundamental:low}
\partial_x^\gamma
\Pi^L_j\arr h(x)
	= \sum_{k=1}^N \sum_{\abs{\beta}=m} \int_{\R^n} 	\partial_x^\gamma\partial_y^\beta F^L_{j,k}(x,y)\,h_{k,\beta}(y)\,dy
	\quad\text{for a.e.\ $x\in\R^n$}
\end{equation}
for all $\arr h\in L^p(\R^n)$ that are also locally in $L^{P}(\R^n)$, for some $P>n/(m-\abs{\gamma})$. In the case of $\abs{\alpha}=m$, we still have that
\begin{equation}
\label{eqn:fundamental:2}
\partial^\alpha\Pi^L_j\arr h(x)
	= \sum_{k=1}^N \sum_{\abs{\beta}=m} \int_{\R^n} 	\partial_x^\alpha\partial_y^\beta F^L_{j,k}(x,y)\,h_{k,\beta}(y)\,dy
	\quad\text{for a.e.\ $x\notin\supp \arr h$}
\end{equation}
for all $\arr h\in L^2(\R^n)$ whose support is not all of $\R^n$.
\end{theorem}
Here, if $\arr h\in L^2(\R^n)$, then $\Pi^L\arr h$ is the unique function in $\dot W^2_m(\R^n)$ that satisfies
\begin{equation*} \sum_{j=1}^\ell \sum_{\abs{\alpha}=m} \int_\Omega \partial^\alpha\varphi_j \,h_{j,\alpha}
=
(-1)^m\sum_{j,k=1}^\ell \sum_{\abs{\alpha}=\abs{\beta}=m} \int_\Omega \partial^\alpha \varphi_j \,a_{\alpha\beta}^{jk}\, \partial^\beta (\Pi^L\arr h)_k
 \end{equation*}
for all $\varphi\in \dot W^2_m(\R^n)$. That is, $u=\Pi^L\arr h$ is the solution to $Lu=\Div_m\arr h$. The formulas~\eqref{eqn:fundamental:low} and \eqref{eqn:fundamental:2} represent the statement that $F^L$ is the fundamental solution for~$L$, that is, that $L_x F^L(x,y)=\delta_y(x)$ in some sense.

Thus, \cite{Bar14p} contains a construction of the fundamental solution for divergence-form operators $L$ of arbitrary order, with no smoothness assumptions on the coefficients of~$L$ or on solutions to $Lu=0$ beyond boundedness, measurability and ellipticity. These results are new even in the second-order case, as there exist second-order operators $L=-\Div A\nabla$ whose solutions do not satisfy the local boundedness estimate~\eqref{eqn:Moser} (see \cite{MazNP91}, \cite{Fre08}) and thus whose fundamental solution cannot be constructed as in \cite{Ros13}.

\subsection{Formulation of Neumann boundary data}
\label{sec:neumann:var}

Recall from Section~\ref{sec:neumann} that even defining the Neumann problem is a delicate matter. In the case of higher-order divergence form operators with variable coefficients, the Neumann problem has thus received little study.

As discussed in Section~\ref{sec:variablediv}, Agranovich has established some well-posedness results for the inhomogeneous problem $Lu=h$ with homogeneous Neumann boundary data. He has also provided a formulation of inhomogeneous Neumann boundary values; see  \cite[Section~5.2]{Agr07}.

This formulation is as follows. Observe that if the test function $\varphi$ does not have zero boundary data, then formula~\eqref{eqn:weaksoln2} becomes
\begin{align}
\label{eqn:neumann:IBP}
\sum_{j=1}^\ell\int_{\Omega} (Lu)_j\, \varphi_j
&=
	(-1)^m \sum_{j,k=1}^\ell \sum_{\abs{\alpha}=\abs{\beta}=m}
	\int_\Omega \partial^\alpha \varphi_j(X) \,a_{\alpha\beta}^{jk}(X)\, \partial^\beta u_k(X)\,dX
	\\\nonumber&\qquad
	+\sum_{i=0}^{m-1}\sum_{j=1}^\ell
	\int_{\partial\Omega} B_{m-1-i}^j u \,\partial_\nu^i \varphi_j\,d\sigma
\end{align}
where $B_i u$ is an appropriate linear combination of the functions $\partial^\alpha u$ where $\abs{\alpha}=m+i$. The expressions $B_{i} u$ may then be regarded as the Neumann data for~$u$. Notice that if $L$ is a fourth-order constant-coefficient scalar operator, then $B_0=-M_A$ and $B_1=K_A$, where $K_A$, $M_A$ are given by~\eqref{eqn:generalDN}.
Agranovich provided some brief discussion of the conditions needed to resolve the Neumann problem with this notion of inhomogeneous boundary data. Essentially the same notion of Neumann boundary data was used in \cite{MitM13A} (a book considering only the case of constant coefficients); an explicit formula for $B_i u$ in this case may be found in \cite[Proposition~4.3]{MitM13A}.

However, there are several major problems with this notion of Neumann boundary data. These difficulties arise from the fact that
the different components $B_i u$ may have different degrees of smoothness. For example, in the case of the biharmonic $L^p$-Neumann problem of Section~\ref{sec:neumann}, the term $M_\rho u= -B_0 u$ is taken in the space $L^p(\partial\Omega)$, while the term $K_\rho u= B_1 u$ is taken in the negative smoothness space $W^p_{-1}(\partial\Omega)$.

If $\Omega$ is a Lipschitz domain, then the space $W^q_1(\partial\Omega)$ of functions with one degree of smoothness on the boundary is meaningful, and so we may define $W^p_{-1}(\partial\Omega)$, $1/p+1/q=1$, as its dual space. However, higher degrees of smoothness on the boundary and thus more negative smoothness spaces $W^p_{-k}(\partial\Omega)$ are not meaningful, and so this notion of boundary data is difficult to formulate on Lipschitz domains. (This difficulty may in some sense be circumvented by viewing the Neumann boundary data as lying in the dual space to $W\!A^p_{m-1+s}(\partial\Omega)$; see \cite{MitM13B,MitM13A}. However, this approach has some limits; for example, there is a rich theory of boundary value problems with boundary data in Hardy or Besov spaces which do not arise as dual spaces (i.e., with $p<1$) and which is thus unavailable in this context.)

Furthermore, observe that as we discussed in Section~\ref{sec:BVP}, a core result needed to approach Neumann and regularity problems is a Rellich identity-type estimate, that is, an equivalence of norms of the Neumann and regularity boundary data of a solution; in the second-order case this may be stated as
\begin{equation*}\doublebar{\nabla_\tau u}_{L^2(\partial\Omega)}\approx \doublebar{\nu\cdot A\nabla u}_{L^2(\partial\Omega)}\end{equation*}
whenever $\Div A\nabla u=0$ in~$\Omega$, for at least some domains~$\Omega$ and classes of coefficients~$A$. In Section~\ref{sec:open}, we will discuss some possible approaches and preliminary results concerning higher-order boundary value problems, in which a higher-order generalization of the Rellich identity is crucial; thus, it will be highly convenient to have notions of regularity and Neumann boundary values that can both be reasonably expected to lie in the space~$L^2$.

Thus, it is often convenient to formulate Neumann boundary values in the following way. Observe that, if $Lu=0$ in~$\Omega$, and $\partial\Omega$ is connected, then for all nice test functions $\varphi$, the quantity
\begin{equation*}
	\sum_{j,k=1}^\ell \sum_{\abs{\alpha}=\abs{\beta}=m}
	\int_\Omega \partial^\alpha \varphi_j \,a_{\alpha\beta}^{jk}\, \partial^\beta u_k
\end{equation*}
depends only on the values of $\nabla^{m-1}\varphi$ on~$\partial\Omega$; thus, there exist functions $M_A^{j,\gamma} u$ such that
\begin{equation}
\label{eqn:neumann:variable}
	\sum_{j,k=1}^\ell \sum_{\abs{\alpha}=\abs{\beta}=m}
	\int_\Omega \partial^\alpha \varphi_j \,a_{\alpha\beta}^{jk}\, \partial^\beta u_k
=
	\sum_{\abs{\gamma}=m-1}\sum_{j=1}^\ell
	\int_{\partial\Omega} M_A^{j,\gamma} u \,\partial^\gamma \varphi_j\,d\sigma
.\end{equation}
We may then consider the array of functions $\arr M_A u$ to be the Neumann boundary values of~$u$. This formulation only requires dealing with a single order of smoothness, but is somewhat less intuitive as there is no explicit formula for $\arr M_A u(X)$ in terms of the derivatives of $u$ evaluated at~$X$.

Also observe that, if we adopt this notion of Neumann boundary data, it is more natural to view the Dirichlet boundary values of~$u$ as the array $\{\partial^\gamma u\big\vert_{\partial\Omega}:\abs{\gamma}=m-1\}$, and not $\{\partial^\gamma u\big\vert_{\partial\Omega}:\abs{\gamma}\leq m-1\}$, as was done in Sections~\ref{sec:BVP} and~\ref{sec:variablediv}. The natural notion of regularity boundary values is then $\{\nabla_\tau\partial^\gamma u\big\vert_{\partial\Omega}:\abs{\gamma}=m-1\}$; again, all components conveniently may then be expected to have the same degree of smoothness.

\subsection{Open questions and preliminary results}
\label{sec:open}

The well-posed\-ness results of Section~\ref{sec:variablediv} cover only a few classes of elliptic differential operators and some special boundary-value problems; the theory of boundary-value problems for higher-order divergence-form operators currently contains many open questions.

Some efforts are underway to investigate these questions. Recall from \cite{BreM13} (see \eqref{eqn:BreM13} above) that the Dirichlet problem
\begin{equation}\label{eqn:open:Dirichlet}
Lu=0 \text{ in }\Omega, \qquad \partial^\alpha u\big\vert_{\partial\Omega}=f_\alpha,
\qquad
\doublebar{u}_{W^p_{m,1-s-1/p}}\leq C \doublebar{\arr f}_{W\!A^p_{m-1+s}(\partial\Omega)}\end{equation}
and the Poisson problem
\begin{equation}\label{eqn:open:Poisson}
Lu=h \text{ in }\Omega, \qquad \partial^\alpha u\big\vert_{\partial\Omega}=0,\qquad
\doublebar{u}_{W^p_{m,1-s-1/p}}\leq C \doublebar{h}_{V^p_{-m,1-s-1/p}}\end{equation}
are well-posed whenever $\abs{p-2}$ and $\abs{s-1/2}$ are small enough.
A similar argument should establish validity of the Neumann problem
\begin{equation}\label{eqn:open:Neumann}
Lu=0 \text{ in }\Omega, \qquad
M_A^{j,\gamma} u = g_{j,\gamma},
\qquad
\doublebar{u}_{W^p_{m,1-s-1/p}}\leq C \doublebar{\arr g}_{(W\!A^q_{m-1+s}(\partial\Omega))^*}\end{equation}
or Poisson problem with homogeneous Neumann data
\begin{equation}\label{eqn:open:Poisson:Neumann}
Lu=h \text{ in }\Omega, \qquad \arr M_A u=0,\qquad
\doublebar{u}_{W^p_{m,1-s-1/p}}\leq C \doublebar{h}_{V^p_{-m,1-s-1/p}}.\end{equation}

Turning to a broader range of exponents $p$ and~$s$, we observe that perturbative results for the Poisson problems \eqref{eqn:open:Poisson} and~\eqref{eqn:open:Poisson:Neumann} are often fairly straightforward to establish. That is, with some modifications to the relevant function spaces, it is possible to show that if \eqref{eqn:open:Poisson} or \eqref{eqn:open:Poisson:Neumann} is well-posed in some bounded domain~$\Omega$, for some operator~$L_0$ and for some $1<p<\infty$, $0<s<1$, and certain technical assumptions are satisfied, then the same problem must also be well-posed for any operator $L_1$ whose coefficients are sufficiently close to those of $L_0$ (in the $L^\infty$ norm).

Recall from Section~\ref{sec:variablediv} that for any $1<p<\infty$ and $0<s<1$, the Dirichlet problem \eqref{eqn:open:Dirichlet}
for boundary data $\arr f$ in the fractional smoothness space $W\!A^p_{m-1+s}(\partial\Omega)$, can be reduced to well-posedness of the Poisson problem \eqref{eqn:open:Poisson}. (Some results are also available at the endpoint $p=\infty$ and in the case $p\leq 1$; the integer smoothness endpoints $s=0$ and $s=1$ generally must be studied using entirely different approaches.)
A similar argument shows that well-posedness of the Neumann problem \eqref{eqn:open:Neumann} follows from well-posedness of the Poisson problem \eqref{eqn:open:Poisson:Neumann} for $1<p\leq \infty$. (In the case of the Neumann problem results for $p\leq 1$ are somewhat more involved.)
A paper \cite{Bar15p} containing these perturbative results is currently in preparation by the first author of the present paper.

A key component in the construction of solutions of \cite{Bar15p} are appropriate layer potentials, specifically, the Newton potential and the double and single layer potentials given by
\def\1{\mathbf{1}}
\allowdisplaybreaks
\begin{align*}
(\Pi^L \arr h)_j (x) &= \sum_{k=1}^\ell\sum_{\abs{\alpha}=m} \int_{\RR^n} \partial_y^\alpha F^L_{j,k}(x,y)\,h_{k,\alpha}(y)\,dy
,\\
(\D^A_\Omega \arr f)_i (x) &=
\1_{\bar\Omega^C}(x) \tilde f_i(x)
-\sum_{\abs{\alpha}=\abs{\beta}=m}\sum_{j,k=1}^\ell \int_{\bar\Omega^C} \partial_y^\alpha F^L_{i,j}(x,y)\,a_{\alpha\beta}^{jk}(y)\,\partial^\beta \!\tilde f_k(y)\,dy
,\\
(\s^A_\Omega \arr g)_j (x) &= \sum_{k=1}^\ell\sum_{\abs{\gamma}=m-1} \int_{\partial\Omega} \partial_y^\gamma F^L_{j,k}(x,y)\,g_{k,\gamma}(y)\,d\sigma(y)
\end{align*}
where $F^L$ denotes the fundamental solution discussed in Section~\ref{sec:fundamental} and where $\tilde f$ is any function that satisfies $\partial^\gamma \tilde f_k=f_{k,\gamma}$ on~$\partial\Omega$. We remark that these are very natural generalizations of layer potentials in the second-order case, and also of various potential operators used in the theory of constant coefficient higher order differential equations; see in particular \cite{MitM13A}.

In particular, to solve the Dirichlet or Neumann problems \eqref{eqn:open:Dirichlet} or \eqref{eqn:open:Neumann}, it was necessary to establish the bounds on layer potentials
\begin{equation*}\doublebar{\D^A_\Omega \arr f}_{W^p_{m,1-s-1/p}}\leq C\doublebar{\arr f}_{W\!A^p_{m-1+s}(\partial\Omega)}
,\qquad
\doublebar{\s^A_\Omega \arr g}_{W^p_{m,1-s-1/p}}\leq C\doublebar{\arr g}_{B^{p,p}{s-1}(\partial\Omega)}\end{equation*}
In \cite{Bar15p}, these bounds are derived from the bound
\begin{equation*}\doublebar{\Pi^L \arr h}_{W^p_{m,1-s-1/p}}\leq C\doublebar{\arr h}_{W^p_{0,1-s-1/p}}.\end{equation*}
This bound is the technical assumption mentioned above; we remark that it is stable under perturbation and is always valid if $p=2$ and $s=1/2$.

Analogy with the second-order case suggests that, in order to establish well-posedness of the $L^2$-Dirichlet, $L^2$-Neumann and $L^2$-regularity problems, a good first step would be to establish the estimates
\begin{equation*}\doublebar{\D^A_\Omega\arr f}_{\mathfrak{X}} \leq C \doublebar{\arr f}_{L^2(\partial\Omega)},
\quad
\doublebar{\D^A_\Omega\arr f}_{\mathfrak{Y}} \leq C \doublebar{\nabla_\tau \arr f}_{L^2(\partial\Omega)},
\quad
\doublebar{\s^A_\Omega\arr g}_{\mathfrak{Y}} \leq C \doublebar{\arr g}_{L^2(\partial\Omega)}
\end{equation*}
for some spaces $\mathfrak{X}$ and~$\mathfrak{Y}$. (We remark that the corresponding bounds for constant coefficient operators are Theorem~4.7 and Proposition~5.2 in \cite{MitM13A}, and therein were used to establish well-posedness results.) In the paper \cite{BarHM15p}, now in the later stages of preparation, Steve Hofmann together with the authors of the present paper have established the bounds
\begin{align}
\label{eqn:S:square}
\int_{\RR^n}\int_0^\infty \abs{\nabla^m \partial_t\s^A \arr g(x,t)}^2\,{t}\,dt\,dx
	& \leq C \doublebar{\arr g}_{L^2(\partial\Omega)}^2
,\\
\label{eqn:D:square}
\int_{\RR^n}\int_0^\infty \abs{\nabla^m \partial_t \D^A \!\arr f(x,t)}^2\,{t}\,dt\,dx
	& \leq C \doublebar{\nabla_x\smash{\arr f}}_{L^2(\partial\Omega)}^2
\end{align}
where $\D^A=\D^A_{\RR^n_+}$, $\s^A=\s^A_{\RR^n_+}$, for scalar operators~$L$, provided that the coefficients~$a_{\alpha\beta}$ are pointwise bounded, elliptic in the sense of \eqref{eqn:varelliptic}, and are constant in the $t$-direction, that is, the direction transverse to the boundary of $\RR^n_+$. As discussed in Section~\ref{sec:variablecomp}, this assumption of $t$-independent coefficients is very common in the theory of second-order differential equations.

We hope that in a future paper we may be able to extend this result to domains of the form $\Omega=\{(x,t):t>\varphi(x)\}$ for some Lipschitz function~$\varphi$; in the second-order case, this generalization may be obtained automatically via a change of variables, but in the higher-order case this technique is only available for equations in composition form (Section~\ref{sec:variablecomp}) and not in divergence form.

In the case of second-order equations $-\Div A\nabla u=0$, where the matrix $A$ of coefficients is real (or self-adjoint) and $t$-independent, a straightforward argument involving Green's theorem establishes the Rellich identity
\begin{equation*}\doublebar{\nabla_x u(\,\cdot\,,0)}_{L^2(\partial\RR^n_+)}
\approx \doublebar{-\vec e\cdot A\nabla u}_{L^2(\partial\RR^n_+)}\end{equation*}
where $-\vec e\,$ is the unit outward normal to~$\RR^n_+$; that is, we have an equivalence of norms between the regularity and Neumann boundary values of a solution~$u$ to $\Div A\nabla u=0$.
Together with boundedness and certain other properties of layer potentials, this estimate leads to well-posedness of the $L^2$-regularity and $L^2$-Neumann problems; it is then a straightforward argument to derive well-posedness of the Dirichlet problem. See \cite{AlfAAHK11}, \cite{AusA11}, \cite{HofKMP13p}, \cite{BarM13p}, and others.

We have hopes that a similar argument will yield the higher-order Rellich identity
\begin{equation*}\doublebar{\nabla_x\nabla^{m-1} u(\,\cdot\,,0)}_{L^2(\partial\RR^n_+)}
\approx \doublebar{\arr M_A u}_{L^2(\partial\RR^n_+)}\end{equation*}
where the Neumann boundary values $\arr M_A u$ are as in Section~\ref{sec:neumann:var},
for solutions $u$ to divergence-form equations with $t$-independent and self-adjoint coefficients (that is, coefficents that satisfy $a_{\alpha\beta}=\overline{a_{\beta\alpha}}$), and that a similar argument will imply well-posedness of higher-order $L^2$ boundary value problems.

The results of Section~\ref{sec:kato} may also lead to well-posedness of $L^2$ boundary-value problems for a different class of operators, namely, operators of block type. Again, this argument would proceed by establishing a Rellich-type identity.

Let us review the theory of second-order divergence-form operators
${\mathbb L}=-\Div {\mathbb A} \nabla$ in $\RR^{n+1}$, where ${\mathbb A}$ is an $(n+1)\times (n+1)$, $t$-independent matrix in block form; that is,
${\mathbb A}_{j,n+1}={\mathbb A}_{n+1,j}=0$ for $1\leq j\leq n$, and ${\mathbb A}_{n+1,n+1}=1$. It is fairly easy to see that one can formally realize the solution to ${\mathbb L}u=0$ in $\RR^{n+1}_+$, $u\big\vert_{\RR^n}=f$, as the Poisson semigroup $u(x,t)=e^{-t\sqrt L} f(x)$, $(x,t)\in \RR^{n+1}_+$. Then the Kato estimate~\eqref{eqn:Kato} essentially provides an analogue of the Rellich identity-type estimate for the block operator $\mathbb L$, that is, the $L^2$-equivalence between normal and tangential derivatives of the solution on the boundary
\begin{equation*}\doublebar{\partial_tu(\,\cdot\,, 0)}_{L^2(\RR^n)}\approx \doublebar{\nabla_x u(\,\cdot\,, 0)}_{L^2(\RR^n)}.\end{equation*}
Boundedness of layer potentials for block matrices also follows from the Kato estimate.

Following the same line of reasoning, one can build a higher order ``block-type'' operator ${\mathbb L}$, for which the Kato estimate \eqref{eqn:Kato} of Section~\ref{sec:kato} would imply a certain comparison between normal and tangential derivatives on the boundary \begin{equation*}\doublebar{\partial_t^mu(\,\cdot\,, 0)}_{L^2(\RR^n)}\approx \doublebar{\nabla^m_xu(\,\cdot\,, 0)}_{L^2(\RR^n)}.\end{equation*}
It remains to be  seen whether these bounds lead to standard well-posedness results. However, we would like to emphasize that  such a result would be restricted to very special, block-type, operators.



\providecommand{\bysame}{\leavevmode\hbox to3em{\hrulefill}\thinspace}
\providecommand{\MR}{\relax\ifhmode\unskip\space\fi MR }
\providecommand{\MRhref}[2]{%
  \href{http://www.ams.org/mathscinet-getitem?mr=#1}{#2}
}
\providecommand{\href}[2]{#2}

\end{document}